\pdfoutput=1
\documentclass[11pt]{amsart}
\usepackage[left=3.2cm, right=3.2cm, top=3.2cm, bottom=3.2cm]{geometry}
\usepackage[T2A]{fontenc} %encoding of output characters, for cyrillic and other stuff
\usepackage[utf8]{inputenc} % encoding of input characters
\usepackage{yfonts} % old German fonts
\usepackage{dirtytalk} % for quotes: \say{dsfsd}
\usepackage{needspace}
\usepackage{graphics}
\usepackage{epsfig}
\usepackage{wrapfig}

\usepackage{csquotes}

\usepackage{tikz} % for pictures
\usepackage{tikz-cd} % for commutative diagrams
\usetikzlibrary{decorations.pathmorphing,shapes,arrows,calc,shapes,decorations.pathreplacing}
\usetikzlibrary{arrows,arrows.meta}
\tikzcdset{arrow style=tikz, diagrams={>=latex}}
\tikzset{>=latex}

\usepackage{mathtools} % for math stuff
\usepackage{verbatim}

\DeclareFontFamily{U}{mathx}{\hyphenchar\font45}
\DeclareFontShape{U}{mathx}{m}{n}{
      <5> <6> <7> <8> <9> <10>
      <10.95> <12> <14.4> <17.28> <20.74> <24.88>
      mathx10
      }{}
\DeclareSymbolFont{mathx}{U}{mathx}{m}{n}
\DeclareFontSubstitution{U}{mathx}{m}{n}
\DeclareMathAccent{\widecheck}{0}{mathx}{"71}
\DeclareMathAccent{\wideparen}{0}{mathx}{"75}

\usepackage[bookmarks=true,pdfencoding=auto,bookmarksdepth=3]{hyperref} % for references inside paper to work
\usepackage{bookmark}
\hypersetup{colorlinks=true,citecolor=black,linkcolor=black}
\usepackage{leftidx} % for easy access to left indxex
\usepackage{amsfonts,amssymb,amsmath,mathrsfs,amsthm,amscd}
\usepackage{subcaption} % for subcaptions in figures
\usepackage{float} % for positioning of figures, see https://en.wikibooks.org/wiki/LaTeX/Floats,_Figures_and_Captions#Figures

\theoremstyle{plain}
\newtheorem{theorem}{Theorem}[section]
\newtheorem{proposition}[theorem]{Proposition}
\newtheorem{lemma}[theorem]{Lemma}
\newtheorem{corollary}[theorem]{Corollary}
\newtheorem*{conjecture*}{Conjecture}
\newtheorem{conjecture}[theorem]{Conjecture}

\theoremstyle{definition}
\newtheorem{definition}[theorem]{Definition}
\newtheorem{example}[theorem]{Example}
\newtheorem{computation}[theorem]{Computation}
\newtheorem{construction}[theorem]{Construction}

\newtheorem{remark}[theorem]{Remark}

\numberwithin{equation}{section}

\newcommand{\ov}[1]{\overline{#1}}
\renewcommand{\tilde}[1]{\widetilde{#1}}
\newcommand{\dtilde}[1]{%
  \tilde{\tilde{#1}}%
}
\newcommand\x{\mathbf x}
\newcommand\y{\mathbf y}

\def\co{\colon\thinspace\relax}

\renewcommand\C{\mathbb{C}}
\newcommand\id{{\text{id}}}
\newcommand\rk{{\text{rk}}}
\renewcommand\Im{{\text{Im}}}
\renewcommand\Re{{\text{Re}}}
\renewcommand\det{{\text{det}}}
\newcommand\B{{\mathcal{B}}}
\newcommand\F{{\mathcal{F}}}
\newcommand\Zpmc{{\mathcal{Z}}}
\newcommand\A{{\mathcal{A}}}
\renewcommand\H{{\mathcal{H}}}
\newcommand\R{\mathbb{R}}

\renewcommand\b{\boldsymbol}
\renewcommand\ge{\geqslant}
\renewcommand\geq{\geqslant}
\renewcommand\le{\leqslant}
\renewcommand\leq{\leqslant}
\renewcommand\o{\otimes}
\newcommand{\lefttorightarrow}{\, \rotatebox[origin=c]{-90}{   $  \circlearrowleft   $ } \, }

\author{Artem Kotelskiy}
\email{artofkot@iu.edu}
\address{Department of Mathematics \\ Princeton University}
\curraddr{Department of Mathematics \\ Indiana University}
\urladdr{http://artofkot.github.io/}

\begin{document}

\title[Explicit]{Comparing homological invariants for mapping classes of surfaces}

\begin{abstract}
    We compare two different types of mapping class invariants: the Hochschild homology of an $A_\infty$ bimodule coming from bordered Heegaard Floer homology, and fixed point Floer cohomology. We first compute the bimodule invariants and their Hochschild homology in the genus two case. We then compare the resulting computations to fixed point Floer cohomology, and make a conjecture that the two invariants are isomorphic. We also discuss a construction of a map potentially giving the isomorphism. It comes as an open-closed map in the context of a surface being viewed as a $0$-dimensional Lefschetz fibration over the complex plane.
\end{abstract}

\maketitle
    
\section{Introduction}
    % Original goal is to find a computable invariant for smooth closed 4-manifolds, which resembles Seiberg-Witten invariants. Following TQFT approach, we can try to define first an invariant for $3$-manifolds, and then maps for 4-dim. cobordisms between them. Let us narrow the class of 3 manifolds we consider down to $3$-manifolds which fiber over the circle, because those appear as preimages of circles in Lefschetz and Broken Lefschetz fibrations (any 4-manifold admits BLF). These 3 manifolds are determined by the monodromy map around the circle, which is a mapping class of a surface. Thus we first want to have a homology theory for mapping classes of surfaces, i.e. fibered $3$-manifolds, and then we want to have maps between homology theories for 4-dim. cobordisms. Here we mainly concentrate on the first task.

    Denote by $\Sigma$ a compact oriented genus $g$ surface, possibly with boundary. We will be studying elements of a strongly based mapping class group $MCG_0(\Sigma,\partial\Sigma=S^1)$, which consists of orientation preserving self-diffeomorphisms $\phi\co (\Sigma,\partial\Sigma) \rightarrow (\Sigma,\partial\Sigma)$ fixing the boundary, up to isotopy.
    \subsection{Overview}
    Suppose we are given a mapping class $\phi \in MCG_0(\Sigma,\partial\Sigma=S^1)$. We will be studying two homological invariants associated to $\phi$: 

    % the <объект>!
    % a <понятие>!

    \begin{enumerate}
        \item  An $A_\infty$ bimodule $N(\phi)$ (or, more precisely, its $A_\infty$ homotopy equivalence class), and its Hochschild homology $HH_*(N(\phi))$, which is a $\mathbb Z_2$-graded vector space over ${\mathbb F}_2$. The bimodule $N(\phi)$ comes from the bordered Heegaard Floer theory:  see~\cite{LOT-mcg}, where the bimodule is also denoted by $N(\phi)$, and the original paper~\cite{LOT-bim}, where the bimodule is denoted by $\widehat{CFDA}(\phi,-g+1)$. In Section~\ref{sec:bimodule from bordered} we cover the original construction of the bimodule $N(\phi)$, whereas in Section~\ref{sec:bimodule_via_Fukaya} we also describe an equivalent construction of this bimodule, using the partially wrapped Fukaya category of the surface. This construction will be useful in understanding the connection with the next invariant. 

        \item Suppose for a moment that $\phi$ is a mapping class of a closed surface, and pick a generic area-preserving representative $\phi$ in that mapping class. Then we consider $HF^*(\phi)$, a $\mathbb Z_2$-graded vector space over ${\mathbb F}_2$ defined as a cohomology of chain complex $CF^*(\phi)$. The generators of $CF^*(\phi)$ are non-degenerate constant sections of the mapping torus $T_\phi \rightarrow S^1$ (i.e. non-degenerate fixed points), and the differentials are pseudo-holomorphic cylinder sections of $T_\phi \times {\mathbb R} \rightarrow S^1 \times  {\mathbb R} $. The same theory can be set up using fixed points as generators and pseudo-holomorphic discs in the Lagrangian Floer cohomology of graphs of $\id$ and $\phi$ as differentials, but we will use the sections and cylinders approach. The invariant $HF^*(\phi)$ is called \emph{fixed point Floer cohomology}, or \emph{symplectic Floer cohomology}. It is $\mathbb Z_2$-graded by the sign of $\det(d\phi - \id )$ at the fixed points of $\phi$.

        In order to generalize this construction to mapping classes fixing the boundary, we have to specify in which direction to twist the boundary slightly to eliminate degenerate fixed points. There are two choices (we call them $+$ and $-$) for each boundary, see Figure~\ref{fig:perturbation_conventions} for the conventions. 

        In our case of a mapping class $\phi \in MCG_0(\Sigma,\partial\Sigma=S^1=U_1)$, we actually consider the induced mapping class $\tilde{\phi}\co (\tilde{\Sigma},\partial \tilde{\Sigma}=U_1 \cup U_2) \rightarrow (\tilde{\Sigma},\partial \tilde{\Sigma}=U_1 \cup U_2)$, where the surface $\tilde{\Sigma}=\Sigma \setminus D^2$ is obtained by removing a disc in the small enough neighborhood of the boundary $U_1$, such that $\phi$ is identity on that neighborhood. We then consider fixed point Floer cohomology $HF^{*}(\tilde{\phi}\:;\:  U_2^+,U_1^-)$ with two different perturbation twists on the two boundaries.

    \end{enumerate}

    The bimodule invariant $N(\phi)$ was computed for mapping classes of the genus one surface in~\cite[Section~10]{LOT-bim}. In Section~\ref{subsec:bim_comps}, we compute $N(\phi)$ in the genus two case; we do so by explicitly describing the bimodules associated to Dehn twists $\tau_l$, which generate the mapping class group. For that we write down the bimodules based on a holomorphic curve count, and then use the description of arc-slide type $DD$ bimodules from~\cite{LOT-arc} to prove that the bimodules $N(\tau_l)$ are the correct ones.
     
    We also describe how to compute Hochschild homology of the bimodule in the genus two case. There is an obstacle we encountered in computing Hochschild homology: none of the smallest models of bimodules $N(\phi)$ for the Dehn twists are bounded, and thus their Hochschild complex is infinitely generated. We write down a certain bounded identity bimodule $[\mathbb{I}]^b$ in the genus two case, so that the bimodule $[\mathbb{I}]^b \boxtimes N(\phi) \boxtimes [\mathbb{I}]^b$ is bounded and belongs to the same $A_\infty$ homotopy equivalence class as $N(\phi)$. Thus by replacing $N(\phi)$ with $[\mathbb{I}]^b \boxtimes N(\phi) \boxtimes [\mathbb{I}]^b$ we solve the problem of $N(\phi)$ being not bounded.

    Based on the computations of Hochschild homology $HH_*(N(\phi))$ in the genus two case, and the corresponding computations of the fixed point Floer cohomology, we make the following conjecture.

    \begin{conjecture} \label{conj:in_intro} For every mapping class $\phi \in MCG_0(\Sigma,\partial\Sigma=S^1=U_1)$ there is an isomorphism of $\mathbb{Z}_2$-graded vector spaces
    $$
    HH_*(N(\phi^{-1})) \cong HF^{*+1}(\tilde{\phi}\:;\:  U_2^+,U_1^-).$$
    \end{conjecture}

    \begin{theorem}\label{thm:int}
    The conjecture above is true in the following two cases of mapping classes of genus two surface $\Sigma_2$: 
    \begin{itemize}
        \item The identity mapping class $\phi=\id \lefttorightarrow \Sigma_2$.
        \item Single Dehn twist along an arbitrary curve $\phi=\tau  \lefttorightarrow \Sigma_2$.
    \end{itemize}
    \end{theorem}

    \noindent The result above is based on computer calculations using~\cite{Pyt}. In fact, we confirmed Conjecture~\ref{conj:in_intro} in many other cases; in Section~\ref{sec:conj_iso} we prove Theorem~\ref{thm:int}, and list some of the other relevant computations.

    We further explain where an isomorphism $HH_*(N(\phi^{-1})) \cong HF^{*+1}(\tilde{\phi}\:;\:  U_2^+,U_1^-)$ may come from. In Section~\ref{sec:bordered<->Fuk}, we describe the spectacular connection, due to Auroux, between bordered Heegaard Floer theory and the partially wrapped Fukaya categories of punctured surfaces and their symmetric products. In the light of this connection the bimodule $N(\phi)$ becomes a graph bimodule $N_{\F_z}(\phi)$, associated to an automorphism of the partially wrapped Fukaya category of a surface $\F_z(\Sigma)$ induced by $\phi$. Moreover, there is a natural $open$-$closed$ map from the Hochschild homology of the graph bimodule into fixed point Floer cohomology, provided that we consider the same kind of Hamiltonian perturbations for both invariants. Based on these Fukaya categorical structures, in Section~\ref{sec:LF} we obtain theoretical evidence for Conjecture~\ref{conj:in_intro}.  Namely, we show how the double basepoint version of Conjecture~\ref{conj:in_intro} can be viewed as an instance of a more general conjecture of Seidel~\cite[Conjecture 7.18]{Sei3II}, which states that the open-closed map in the Fukaya-Seidel category of a Lefschetz fibration is an isomorphism. 

    Below we drew a diagram which concisely describes the main invariants and relationships between them, which we study in this paper:
    $$
    \begin{tikzcd}[row sep=50pt,column sep=25pt]
        &  
        \text{Mapping class }\phi \lefttorightarrow (\Sigma, \partial \Sigma )
        \arrow[dl, rightsquigarrow,  "\text{Section~\ref{alt_bim}}" left ]
        \arrow[d, rightsquigarrow, "\text{Section~\ref{sec:bimodule from bordered}}" left]
        \arrow[ddr, rightsquigarrow, "\text{Section~\ref{sec:fixed point Floer}}"]
        &
        \\
        N_{\F_z}(\phi)
        \arrow[d, rightsquigarrow, "\text{Section~\ref{sec:hf_bim}}" right] 
        &
        N(\phi)
        \arrow[l, leftrightarrow, "\simeq" above, "\text{Prop.~\ref{prop:bims_are_eq}}" below] 
        \arrow[d, rightsquigarrow, "\text{Section~\ref{sec:hf_bim}}" left] 
        &
        \\
        HH_*(N_{\F_z}(\phi))
        \arrow[rr, bend right=25, "\text{Section~\ref{subsec:O-C map}}" above, "\text{(in the double basepoint case)}" below]
        & 
        HH_*(N(\phi))
        \arrow[l, leftrightarrow, "\cong" above, "\text{Corollary~\ref{cor}}" below]
        \arrow[r, leftrightarrow, "\overset{?}{\cong}" above, "\text{Conjecture~\ref{conj}}" below] 
        & 
        HF^*(\tilde{\phi}\:;\ U_2^+,U_1^-)
    \end{tikzcd}
    $$

    In addition to stating the conjecture in Section~\ref{sec:statement}, performing lots of computations supporting the conjecture in Section~\ref{supporting_comps}, and discussing theoretical evidence for the conjecture in Section~\ref{sec:LF}, some of the main contributions of the present paper are the new bimodule and Hochschild homology computations, made possible by a software package~\cite{Pyt}. From the point of view of low-dimensional topology, in Section~\ref{subsec:bim_comps} we compute the bordered Heegaard Floer bimodule invariants for mapping classes of the genus two surface (before they were only computed for the torus). From the point of view of symplectic geometry, the computed invariants are the graph bimodules for the partially wrapped Fukaya category of the genus two surface. These bimodule computations, while tedious, follow from the paper~\cite{LOT-arc} 
    \footnote{In~\cite{LOT-arc} the authors compute $\widehat{CFDD}$ bimodule for all arc-slides, and so one can tensor those bimodules with $\widehat{CFAA}(\id)$ (which was studied in~\cite{Boh2}), and obtain $\widehat{CFDA}$ bimodules for arc-slides. In Section~\ref{subsec:bim_comps} we took a slightly different, reverse approach: we guessed (based on holomorphic theory) the $\widehat{CFDA}$ bimodules for arc-slides, and then proved that those are the right ones by tensoring with $\widehat{CFDD}(\id)$.}
    . More importantly, having computed those bimodules, in Section~\ref{alg_hoh} we introduce an algorithm for computing their Hochschild homologies. Computations with the bimodules and Hochschild homologies quickly get too complicated to do by hand, and so we decided to develop a python package~\cite{Pyt}. This program is the only software which computes the $\widehat{CFDA}$ bimodules in the second lowest $\text{Spin}^c$ structure, although only for the genus one and two cases\footnote{Zhan's python package~\cite{Pyt-Zha} allows one to work with bimodules $\widehat{CFDA}(\phi,0)$, which are the middle $\text{Spin}^c$ summands.}. Moreover, it also computes the Hochschild homology of bimodules, and thus makes new computations of knot Floer homology possible. Namely, given the factorization of a mapping class $\phi$ into Dehn twists, the program can compute knot Floer homology of the binding of the open book corresponding to monodromy $\phi$, in the second lowest Alexander grading; see~\cite[Showcase~1]{Pyt}. This is based on the fact that Hochschild homology is equal to knot Floer homology. 

    Also, assuming Conjecture~\ref{conj:in_intro} is true, our computational methods allow to effectively compute the number of fixed points of a mapping class $\tilde{\phi}$ by simply running a program, even in the pseudo-Anosov case. For example, taking a mapping class $\psi=\tau_A \tau_B^{-1}$ (see Figure~\ref{fig:Sigma_2}) as an input, the program~\cite[Showcase~2]{Pyt} outputs:
    \begin{center}
    \begin{tabular}{ |c|c|c| } 
        \hline
        Automorphism & Number of fixed points \\ 
        \hline
        $\psi$ & 5   \\ 
        $\psi^2$ & 9  \\
        $\psi^3$ & 20  \\
        $\psi^4$ & 49 \\
        $\psi^5$ & 125 \\ 
        \hline
    \end{tabular}
    \end{center}
    As a byproduct of simplicity of such computations, we can determine if the mapping class is periodic, reducible with all components periodic, or pseudo-Anosov: the rank of $HH_*(N(\phi^{n}))$ is respectively bounded, grows linearly, or grows exponentially (see~\cite[Corollary 4.2]{LOT-mcg},~\cite[Corollary 1.7]{CC}).

    \subsection{Context}
    Here we discuss possible generalizations of Conjecture~\ref{conj:in_intro}, and other closely related constructions.

    First, let us explain why we work only with the $1$-strand moving summand $N(\phi)=\widehat{CFDA}(\phi,-g+1)$ of the full bimodule $\bigoplus_{0 \leq k \leq 2g} \widehat{CFDA}(\phi;-g+k)$. There is a Fukaya categorical interpretation of $\widehat{CFDA}(\phi;-g+1)$ as a graph bimodule of $\phi$, which we discuss in~Section~\ref{sec:bimodule_via_Fukaya}. This suggests that the Fukaya categorical interpretation of $\widehat{CFDA}(\phi;-g+k)$ is a graph bimodule of the induced symplectomorphism $Sym^{k}(\phi)\co Sym^{k}(\Sigma) \rightarrow Sym^{k}(\Sigma)$ (see~\cite[Section~1.1]{Perutz_LMII} for how to obtain a canonical symplectomorphism $Sym^{k}(\phi)$). Thus, in the light of the existence of the open-closed map (see Section~\ref{subsec:O-C map} for the $k=1$ case), the Hochschild homology $HH_*(\widehat{CFDA}(\phi, -g+k))$ should be isomorphic to some version of fixed point Floer cohomology of $Sym^{k}(\phi)$. While this invariant can be defined, if $k\ge 2$ there are no effective computational tools to manage it. However, in the $k=1$ case, there are a wealth of methods to compute fixed point Floer cohomology for mapping classes of surfaces (see Section~\ref{subsec:exist_comp_methods}). This allows us to compare (a specific version of) fixed point Floer cohomology of $\phi \lefttorightarrow \Sigma$ to $HH_*(N(\phi))$. 
    
    Second, there is another isomorphism of invariants of mapping classes of surfaces, which is directly related to our work. Consider a $3$-manifold $Y^3_\phi$, where the subscript $\phi$ indicates that the manifold fibers over the circle $Y^3_\phi \rightarrow S^1$ with the monodromy $ \phi \lefttorightarrow  \Sigma $, where $\Sigma$ is a closed surface. The statement says that the Heegaard Floer homology in the second lowest $\text{Spin}^c$ structures (evaluating to $-2g+4$ on the fiber) is equal to the fixed point Floer cohomology of the corresponding monodromy:  $HF^+_*(Y^3_\phi\:;\:  -2g+4) \cong HF^*(\phi)$. Assuming the genus is greater than 2, this isomorphism follows from the following chain of isomorphisms:
    
    \begin{itemize}

    \item $HF^*(\phi) \cong  HP^*_{degree=1}(\phi)$, \\ where $HF^*(\phi)$ is the fixed point Floer cohomology of mapping class $\phi$ of a closed surface, see~\cite{Sei2}, and $HP^*_{degree=1}(\phi)$ is periodic Floer cohomology, see~\cite{HS}. This isomorphism is addressed in~\cite[Appendix B]{LT}.

    \item $HP^*_{degree=1}(\phi) \cong \widecheck{HM}_*(Y^3_\phi,c_+ \:;\:  -2g+4)$, \\ where $Y^3_\phi$ is the $3$-manifold fibered over the circle with a monodromy $\phi$, and $\widecheck{HM}_*(Y^3_\phi,c_+ \:;\:  -2g+4)$ denotes an invariant of a $3$-manifold called monopole Floer homology, defined in~\cite{KM}. The $c_+$ indicates the version of $\widecheck{HM}_*$ with a monotone positive perturbation. The isomorphism was proved in~\cite{LT}. The positivity of the perturbation comes from the genus being greater than 2. 

    \item $\widecheck{HM}_*(Y^3_\phi,c_+ \:;\:  -2g+4) \cong \widecheck{HM}_\bullet(Y^3_\phi,c_+ \:;\:  -2g+4)$, \\
    where $\bullet$ indicates the negative completion of the coefficient ring. The definition of $\widecheck{CM}_*$ does not depend on the negative completion, so the above groups are isomorphic, see~\cite[p. 606]{KM}.
    
    \item $\widecheck{HM}_\bullet(Y^3_\phi,c_+ \:;\:  -2g+4) \cong \widecheck{HM}_\bullet(Y^3_\phi \:;\:  -2g+4)$, \\
    where the absence of $c_+$ indicates exactness of the perturbation in the definition of $\widecheck{HM}_\bullet$. The isomorphism is proved in~\cite[Theorems 31.1.2]{KM}.
    
    \item $\widecheck{HM}_\bullet(Y^3_\phi \:;\:  -2g+4) \cong \widecheck{HM}_\bullet(Y^3_\phi,c_b \:;\:  -2g+4)$, \\
    where $c_b$ indicates the balanced perturbation in the definition of $\widecheck{HM}_\bullet$. The isomorphism is proved in~\cite[Theorems 31.1.1]{KM}.
    
    \item $\widecheck{HM}_\bullet(Y^3_\phi,c_b \:;\:  -2g+4) \cong \widecheck{HM}_*(Y^3_\phi,c_b \:;\:  -2g+4) $, \\
    which again follows from the fact that the negative completion of the coefficient ring does not affect $\widecheck{HM}_*$.

    \item $\widecheck{HM}_*(Y^3_\phi,c_b \:;\:  -2g+4) \cong HF^+_*(Y_\phi^3\:;\:  -2g+4)$, \\
    which is a deep and very difficult theorem, despite the fact that the definition of Heegaard Floer homology was inspired by monopole Floer homology type constructions. It was proved via passing through another invariant, called embedded contact homology (ECH), in~\cite{KLT,KLT2,KLT3,KLT4,KLT5}, and~\cite{CGH,CGH1,CGH2,CGH3}.

    \end{itemize}

    % Such isomorphism was initially expected because of the following: Heegaard Floer homology and monopole Floer homology count solutions to Seiberg-Witten equations, and in a special $\text{Spin}^c$ structures its dimension reduction to a surface gives the surface itself as a space of solutions (in other $\text{Spin}^c$ structures solutions are points of higher symmetric products $Sym^k(\Sigma)$). So on a mapping torus $Y^3_\phi$ one would expect that the global solutions to Seiberg-Witten equations are related to fixed points of the monodromy. 
    
    Our Conjecture~\ref{conj:in_intro} is analogous to the proved isomorphism $HF^+(Y^3_\phi\:;\:  -2g+4) \cong HF(\phi)$. We work in a slightly different $3$-manifold. Suppose we fix a lift of $\phi$ from the mapping class group of a closed surface $MCG(\Sigma)$ to the strongly based mapping class group $MCG_0(\Sigma,\partial \Sigma=U_1)$. Then, instead of the fibered manifold $Y^3_\phi$, we consider the open book corresponding to $\phi$; denote this open book by $M_\phi^\circ$ and its binding by $K$. Manifolds $Y^3_\phi$ and $M_\phi^\circ$ are related: $M_\phi^\circ$ is obtained from $Y^3_\phi$ by $0$-surgery on the constant section of $Y^3_\phi \rightarrow S^1$, which comes from the lift of $\phi$;  $Y^3_\phi$ is obtained from $M_\phi^\circ$ by $0$-surgery on $K$. Instead of $HF^+(Y^3_\phi\:;\:  -2g+4)$ we consider the knot Floer homology of the binding (in the second lowest Alexander grading) $\widehat{HFK}(M_\phi^\circ,K \:;\:  -g+1)$. It is equal to the Hochschild homology $HH_*(N(\phi))$ (see~\cite[Theorem~14]{LOT-bim}), with which we are actually working in this paper. The relevant version of fixed point Floer cohomology turns out to be $HF^{*}(\tilde{\phi}\:;\:  U_2^+,U_1^-)$. It is possible that the Conjecture~\ref{conj:in_intro} can be deduced from  the proved isomorphism $HF^+(Y^3_\phi\:;\:  -2g+4) \cong HF(\phi)$.
    
    It is also interesting to compare our results to the work of Spano~\cite{Sp}. He develops the full version of embedded contact knot homology
    %(the \say{hat} version was developed in [todo] via sutured embedded contact homology)
    $ECK(Y,K, \alpha)$, and conjectures it to be isomorphic to $HFK^-(Y,K)$. In~\cite[Section~3.3.1]{Sp}, the connection to symplectic Floer homology is explained. Namely, in case of the knot being a binding of an open book, the embedded contact knot homology is equal to a certain periodic Floer homology, see~\cite[Theorem~3.19]{Sp}. In the degree one case, it follows that if $ECK(Y, K, \alpha)\cong HFK^-(Y,K)$, then  $HF(\phi, U_1+)\cong HFK^-(M_\phi^\circ,K \:;\:  -g+1).$ Thus our Conjecture~\ref{conj:in_intro} can be viewed as the \say{hat} version of the conjecture of Spano.

    Working in the \say{hat} version allows us to consider Hochschild homology of the bimodule $HH(N(\phi))$ instead of the knot Floer homology $\widehat{HFK}(M_\phi^\circ,K \:;\:  -g+1)$. This transition is quite powerful, because two things become possible: computations using bordered Floer theory, and the connection to the partially wrapped Fukaya category of a surface, specifically to twisted open-closed maps there. The latter connection provides hope that the Conjecture~\ref{conj:in_intro} can be proved by more algebraic methods, using the structure of the Fukaya category. In this direction see~\cite{Gan}, where it is proved that the 
    untwisted open-closed map is an isomorphism for the non-degenerate %i.e. there is an essential collection, i.e. Lagrangians which generate the category. By Ab criterion its the same as requiring 1 to be hit by O-C.
    wrapped Fukaya category, in the exact setting. %Ultimately, it is possible, that after establishing bordered theory for $HF^-$ version of Heegaard Floer homology, and its relationship to relative Fukaya category, Spano's conjecture can also be proved as an isomorphism induced by an open-closed map.
    
    \subsection{Outline of the paper}
    This paper is organized as follows. For clarity comments on the novelty of the content are added.
        \begin{description}
            \item[Section 2.1] we cover the background on bordered Heegaard Floer theory, which is relevant to the construction of the bimodule invariant $N(\phi)$. This is an exposition of results from~\cite{LOT-bim}.
            \item[Section 2.2] we perform computations of the bimodule invariant $N(\phi)$ in the genus two case. These computations are new.
            \item[Section 2.3] we describe a method to compute Hochschild homology $HH_*(N(\phi))$ of bimodules in the genus two case, which we implemented in the program~\cite{Pyt}. The method and the computations, made possible by~\cite{Pyt}, are new.
            \item[Section 3.1] we sketch the definition of fixed point Floer cohomology $HF^*(\phi)$, in both cases $\partial \Sigma = \emptyset$ and $\partial \Sigma \neq \emptyset$. This is an exposition of results from~\cite{Flo,DS,Sei2}, based on other expositions from~\cite{Gau,CC,Sei3II,Ul}.
            \item[Section 3.2] we describe the existing methods to compute fixed point Floer cohomology. This is an exposition of results from~\cite{Sei1,Gau,Eft,CC}.
            \item[Section 4.1] we state our main conjecture. This is new.
            \item[Section 4.2] we perform various computations in the genus two case, and notice the equality of ranks of the groups $HF^{*}(\tilde{\phi}\:;\:  U_2^+,U_1^-)$ and $HH_*(N(\phi))$, which confirms the conjecture. We also prove Theorem~\ref{thm:int}. The results are new.
            \item[Section 5.1] we describe Auroux's symplectic geometric interpretation of bordered Heegaard Floer homology in terms the partially wrapped Fukaya categories. This is an exposition of~\cite{Aur1,Aur2}.
            \item[Section 5.2] we explain the Fukaya categorical interpretation of the bimodule $N(\phi)$. This interpretation was known, see~\cite[Lemma~4.2]{AGW14}.
            \item[Section 6.1] we explain the generalization of the bimodule $N(\phi)$ to the double basepoint version $N^\text{2bp}(\phi)$. This is an exposition of known constructions, which are based on~\cite{Zarev}.
            \item[Section 6.2] we introduce the double basepoint generalization $HF^{*}(\dtilde{\phi}\:; \ U_3^+,U_2^+,U_1^-)$ of the fixed point Floer cohomology $HF^{*}(\tilde{\phi}\:;\:  U_2^+,U_1^-)$, and state the double basepoint version of our main conjecture. This material is new. 
            \item[Section 6.3] we describe the background material on Lefschetz fibrations and the Fukaya-Seidel category. This is an expositions of results from~\cite{Sei4,Sei5,Sei3I,Sei3II}.
            \item[Section 6.4] we show how our conjecture is a special case of the conjecture of Seidel~\cite[Conjecture 7.18]{Sei3II}, which states that the open-closed map in the Fukaya-Seidel category is an isomorphism. While the  discussion of the open-closed map is  an exposition of~\cite[Section~7]{Sei3II}, the connection between our work and~\cite{Sei3II} is new.
        \end{description}

    \vspace{0.2cm}

    {\setlength{\parindent}{0pt} \bf Assumptions and conventions.}
    \begin{itemize}
        \item By $\phi$ we will usually denote not only the diffeomorphism, but also the mapping class which it represents.
        \item We will use the convention $\omega(X_H,\cdot)=-dH$ for Hamiltonian vector fields.
        \item Every homological invariant we consider will be defined over the field ${\mathbb F}_2$.
        \item We will be working with the fixed point Floer $co$homology, rather than $ho$mology.
    \end{itemize}

    \medskip

    {\setlength{\parindent}{0pt} \bf Acknowledgements.}
    I am thankful to my adviser Zoltán Szabó for his guidance and support during this project. I thank Nick Sheridan for suggesting I read one of the key references~\cite{Sei3II}. I also would like to thank Denis Auroux, Sheel Ganatra, Peter Ozsváth, Paul Seidel, András Stipsicz, Mehdi Yazdi, and Bohua Zhan for helpful conversations. Finally, I thank the referees for their helpful comments and suggestions.

\Needspace{8\baselineskip}
\section{Bimodule invariant coming from bordered Heegaard Floer homology}
    Everything in this section is based on the bordered Heegaard Floer theory. It was developed by  Lipshitz, Ozsváth and Thurston in~\cite{LOT-main} and~\cite{LOT-bim}. We refer to those papers for the theory of $A_\infty$ algebras, modules and bimodules, for how such objects arise in the Heegaard Floer theory, and for the proofs of the propositions we state and use along the way.
    \subsection{Background: bimodules and their Hochschild homology}\label{sec:bimodule from bordered}
        \subsubsection{Pointed matched circles}\label{sec:pmc}
            We will be considering surfaces with one boundary component. Moreover, it is useful to consider parameterized surfaces, i.e. surfaces with a specified $1$-handle decomposition. Thus let us start with the following definition.

            \begin{definition}
            A \emph{pointed matched circle} is an oriented circle $\Zpmc$, equipped with a basepoint $z$ on it, and additional $4g$ points coming in pairs (distinct from each other and $z$) such that performing surgery on all $2g$ pairs results in one circle.
            \end{definition}

\begin{figure}[h]
            \includegraphics[width=0.5\textwidth]{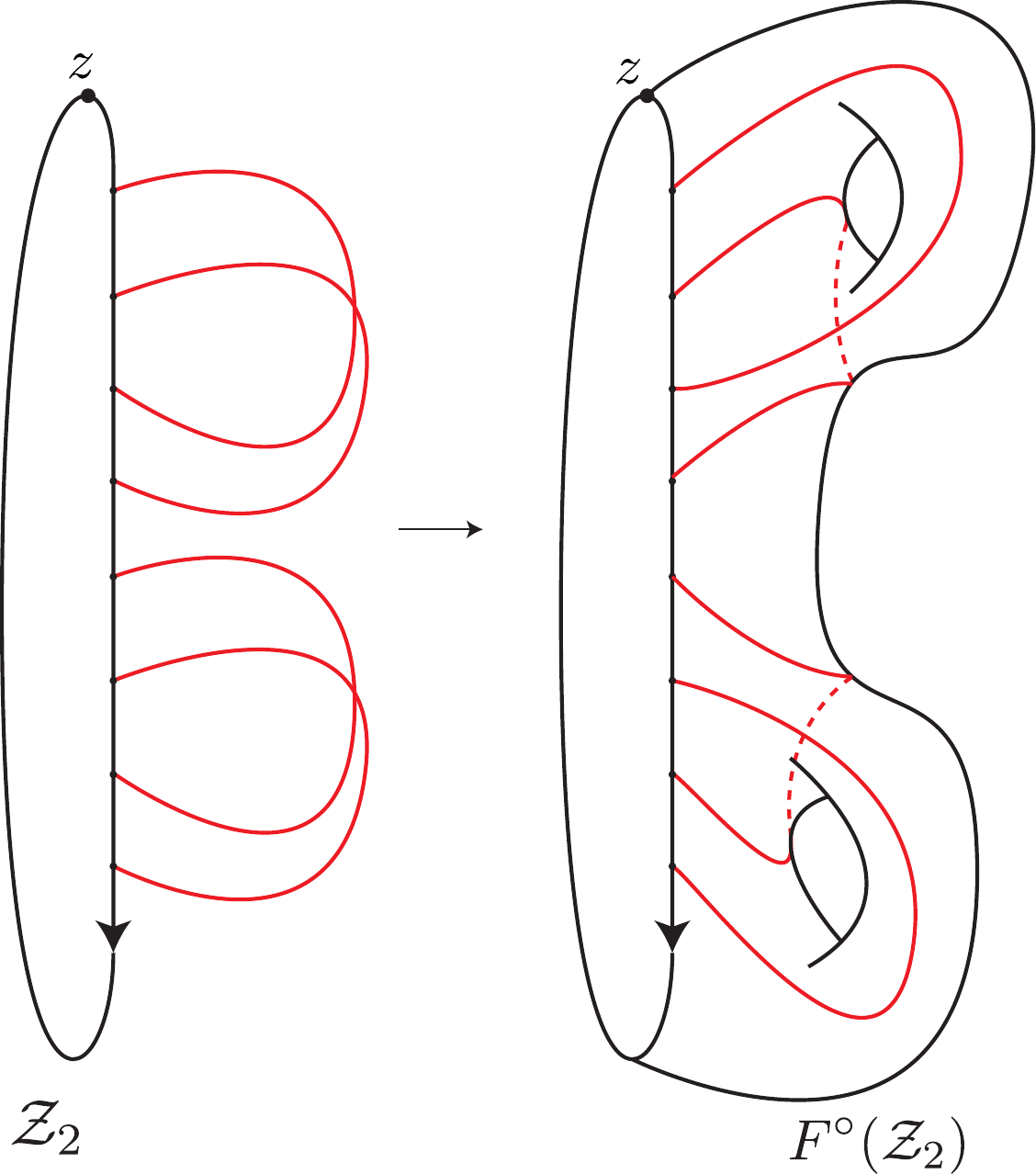}
            \caption{An example of a pointed matched circle for $g=2$, and its corresponding genus two surface.}
            \label{fig:pointed_matched_circle}
            \end{figure}

            \begin{construction}[The surface associated to a pointed matched circle]
            Given a pointed matched circle $\Zpmc$, we can associate a surface whose boundary is a circle $\Zpmc$, viewing the $2g$ pairs of points as feet of $1$-handles. Specifically, we thicken $\Zpmc$ into a band $\Zpmc \times [0,1]$, then glue the $1$-handles to $\Zpmc \times \{1\}$, and then cap off the boundary component which is not $\Zpmc \times \{0\}$ (see below Figure~\ref{fig:pointed_matched_circle} and~\cite[Figure 1.1]{LOT-main}). We denote this surface by $F^{\circ}(\Zpmc)$, and the orientation on it is induced from the boundary via the usual rule \say{outward normal first}. Let $F(\Zpmc)$ denote the result of filling in a disc $D_\Zpmc$ to the boundary component of $F^{\circ}(\Zpmc)$ (so mapping classes of $F^{\circ}(\Zpmc)$ fixing the boundary naturally correspond to mapping classes of $F(\Zpmc)$ fixing the disc $D_\Zpmc$). Note that any two surfaces specified by the same pointed matched circle are homeomorphic, via a homeomorphism which is uniquely determined up to isotopy.
            \end{construction}

            \begin{example}
            In Figure~\ref{fig:pointed_matched_circle}, we provide an example of a pointed matched circle in the $g=2$ case, and its corresponding surface of genus two. In our computations of the mapping class invariant we will be using this pointed matched circle, which we denote by $\Zpmc_2$. Notice that there are other pointed matched circles for a genus two surface, not isomorphic to $\Zpmc_2$. We could have used them. Thus here we make a particular choice which can be understood as a choice of a parameterization of the surface by specifying a $0$-handle (the preferred disc) and $1$-handles.
            \end{example}

            Having a pointed matched circle $\Zpmc$, by $-\Zpmc$ we denote the same circle but with the reversed orientation. The corresponding surfaces also have opposite orientations: $F^{\circ}(-\Zpmc)= -F^{\circ}(\Zpmc)$.

            Consider  now \emph{the genus $g$ strongly based mapping class groupoid}, which is a category where the objects are pointed matched circles with $4g$ points, and the morphism sets are
            $$MCG_0(\Zpmc_L,\Zpmc_R)=\{\phi\co F^{\circ}(\Zpmc_L)\xrightarrow{\cong} F^{\circ}(\Zpmc_R) \ | \  \phi(z_L)=z_R\} /\sim,$$
            i.e. orientation preserving diffeomorphisms respecting the boundary and the basepoint, considered up to isotopy. For any pointed matched circle $\Zpmc$ with $4g$ points the corresponding group of self-diffeomorphisms $MCG_0(\Zpmc,\Zpmc) \cong MCG_0(\Sigma,\partial \Sigma)$ is the mapping class group of the genus $g$ surface with one boundary component.

            Our goal for the rest of Section~\ref{sec:bimodule from bordered} is to explain how to associate a (homotopy equivalence class of) type $DA$ bimodule to an element $\phi \in MCG_0(\Zpmc_L,\Zpmc_R)$. This uses a number of correspondences, which are outlined in Figure~\ref{fig:plan_2}. If we take $\Zpmc_L=\Zpmc_R$, then we will produce an invariant of a mapping class of a surface. We now proceed to explaining the different pieces of the diagram.

            \begin{figure}[h]
            \tikzstyle{block} = [rectangle, rounded corners, minimum width=3cm, minimum height=1cm, text centered, text width=8cm]
            \tikzstyle{arrow} = [thick,->,>=stealth]
            \begin{tikzpicture}[node distance = 2.4cm]
                \node [block] (mcg)  {$MCG_0(\Zpmc_L,\Zpmc_R)$};
                
                \node [block, below of=mcg] (mapping_cylinders) {Mapping cylinders up to isomorphisms};
                
                \node [block, below of=mapping_cylinders] (heegard_diagrams) {Heegaard diagrams of mapping cylinders up to certain moves};
                
                \node [block, below of=heegard_diagrams](bimodules) {$A_\infty$ bimodules ${}^{\A(\Zpmc_L)}N_{\A(\Zpmc_R)}$ up to homotopy equivalence};

                \draw[decorate,decoration={brace,mirror,raise=3pt,amplitude=5pt}, thick]
                    (mapping_cylinders.north west)--(mapping_cylinders.south west) ;
                \draw[decorate,decoration={brace,raise=3pt,amplitude=5pt}, thick]
                    (mapping_cylinders.north east)--(mapping_cylinders.south east); 
                \draw[decorate,decoration={brace,mirror,raise=3pt,amplitude=5pt}, thick]
                    (heegard_diagrams.north west)--(heegard_diagrams.south west) ;
                \draw[decorate,decoration={brace,raise=3pt,amplitude=5pt}, thick]
                    (heegard_diagrams.north east)--(heegard_diagrams.south east); 
                % \draw[decorate,decoration={brace,mirror,raise=3pt,amplitude=5pt}, thick]
                %     (mcg.north west)--(mcg.south west) ;
                % \draw[decorate,decoration={brace,raise=3pt,amplitude=5pt}, thick]
                %     (mcg.north east)--(mcg.south east); 
                \draw[decorate,decoration={brace,mirror,raise=3pt,amplitude=5pt}, thick]
                    (bimodules.north west)--(bimodules.south west) ;
                \draw[decorate,decoration={brace,raise=3pt,amplitude=5pt}, thick]
                    (bimodules.north east)--(bimodules.south east); 

                \draw [arrow,<->,shorten >=0.2cm,shorten <=0.2cm] (mcg) -- node[anchor=west]{1-1} (mapping_cylinders);
                \draw [arrow,<->,shorten >=0.2cm,shorten <=0.2cm] (mapping_cylinders) -- node[anchor=west]{1-1} (heegard_diagrams);
                \draw [arrow,shorten >=0.2cm,shorten <=0.2cm] (heegard_diagrams) -- (bimodules);
            \end{tikzpicture}
            \caption{}
            \label{fig:plan_2}
            \end{figure}
        \subsubsection{Mapping cylinders}
            We need the notion of a \emph{strongly bordered $3$-manifold} with two boundary components. It consists of the following data (following~\cite[Definition 5.1]{LOT-bim}):
            \begin{enumerate}
                \item an oriented $3$-manifold $Y$ with two boundary components $\partial_1 Y$ and $\partial_2 Y$;
                \item a preferred disc and a basepoint (on the boundary of that disc) in each boundary component;
                \item parameterizations of each boundary components by some fixed surfaces $\psi_i\co(F_i, D_i, z_i)\rightarrow\partial_i Y$ respecting distinguished discs and basepoints;
                \item a framed arc connecting the basepoints such that the framing on the boundaries points into the distinguished discs.
            \end{enumerate}
            
            Given a parameterization of the boundaries of $Y$ by surfaces $F_1$ and $F_2$, then there is a natural notion of isomorphism of strongly bordered $3$-manifolds --- it is a diffeomorphism of the corresponding $3$-manifolds, which respects every piece of the additional data, i.e. parameterizations of boundaries, arcs connecting the basepoints, and their framings. 

            Given a strongly based mapping class we wish to form a strongly bordered $3$-manifold, which is called \emph{the corresponding mapping cylinder}.

\begin{figure}[h]
            \centering
            \includegraphics[width=1\textwidth]{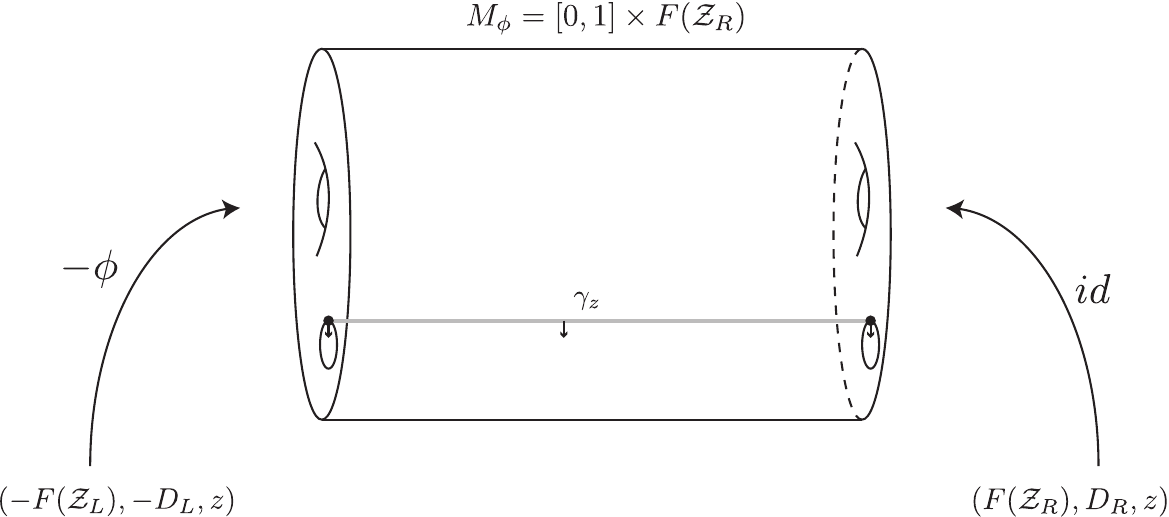}
            \caption{Mapping cylinder of $\phi\co (F(\Zpmc_L),D_L,z_L)\rightarrow (F(\Zpmc_R),D_R,z_R)$.}
            \label{fig:mapping_cylinder}
            \end{figure}

            \begin{construction}[Mapping cylinder]\label{mapping cylinder}
            Fix pointed matched circles $\Zpmc_L$, $\Zpmc_R$ and a mapping class $\phi\co (F(\Zpmc_L),D_L,z_L)\rightarrow (F(\Zpmc_R),D_R,z_R)$. We can form a \emph{mapping cylinder } \\ $M_\phi={}_\phi([0, 1]\times F(\Zpmc_R))_{\id }$, which is a strongly bordered $3$-manifold  with the following data: 
            \begin{enumerate} 
                \item The $3$-manifold and its boundary components are 
                $$Y = [0, 1] \times F(\Zpmc_R), \qquad  \partial_L Y =  \{0\} \times -F(\Zpmc_R), \qquad \partial_R Y =\{1\} \times F(\Zpmc_R)$$
                \item a parametrization of its boundary given by 
                $$\psi_L=-\phi\co -F(\Zpmc_L) \rightarrow \partial_L Y, \qquad \psi_R=\id \co F(\Zpmc_R) \rightarrow \partial_R Y ;$$
                \item two distinguished discs $\{0\}\times D_R$ in $\partial_L Y$ and $\{1\}\times D_R$ in $\partial_R Y$;
                \item a framed path $\gamma_z=[0,1]\times \{z_R\}$ between $z_L \in \partial_L Y$ and $z_R \in \partial_R Y$ such that the framing points into the distinguished discs $D_R$ at every fiber $\{t\} \times F(\Zpmc_R)$. See Figure~\ref{fig:mapping_cylinder}.
            \end{enumerate}
            \end{construction}

            The following lemma allows us to talk about mapping cylinders instead of mapping classes (and vice versa), i.e. it explains the first correspondence in Figure~\ref{fig:plan_2}.

            \begin{lemma}[\text{\cite[Lemma~5.29]{LOT-bim}}]
            Fix pointed matched circles $\Zpmc_L$ and $\Zpmc_R$. Then any strongly bordered $3$-manifold Y, whose boundary is parameterized by $F(\Zpmc_L)$ and $F(\Zpmc_R)$, and whose underlying space can be identified with a product of a surface with an interval (so that arc $\gamma_z$ is identified with the product of a point with the interval, respecting the framing) is of the form $M_\phi$ for some choice of strongly based mapping class $\phi\co F^\circ(\Zpmc_L) \rightarrow F^\circ(\Zpmc_R)$. Moreover, two such strongly bordered three-manifolds are isomorphic if and only if they represent the same strongly based mapping class. 
            \end{lemma}
        \subsubsection{Heegaard diagrams}\label{sec:hd}
            Now, having constructed the mapping cylinder $M_\phi$, we would like to have a $2$-dimensional presentation of it.

            \begin{definition}
            An \emph{arced bordered Heegaard diagram with two boundary components} is a quadruple $(\ov{\Sigma},\ov{\b\alpha},\b{\beta},\b{z})$ where
            \begin{itemize}
                \item $\ov{\Sigma}$ is an oriented compact surface of genus $g$ with two boundary components, $\partial_L \ov{\Sigma}$ and $\partial_R \ov{\Sigma}$;

                \item $\ov{\b{\alpha}}=\{\ov{\alpha}^\text{arc,left}_1,\dots,\ov{\alpha}^\text{arc,left}_{2l},\ov{\alpha}^\text{arc,right}_1,\dots,\ov{\alpha}^\text{arc,right}_{2r},\dots,\alpha_1^\text{curve},\dots,\alpha_{g-l-r}^\text{curve}\}$ is a collection of pairwise disjoint $2l$ embedded arcs with boundaries on $\partial_L \ov{\Sigma}$, $2r$ embedded arcs with boundaries on $\partial_R \ov{\Sigma}$, and $g-l-r$ circles in the interior (in particular $g \geq l+r$);

                \item $\b{\beta}=\{\beta_1,\cdots,\beta_g\}$ is a $g$-tuple of pairwise disjoint curves in the interior of $\ov{\Sigma}$;

                \item $\b z$ is a path in $\ov{\Sigma} \setminus (\ov{\b \alpha} \cup \b \beta )$ between $\partial_L \ov{\Sigma}$ and $\partial_R \ov{\Sigma}$;
            \end{itemize}
            These are required to satisfy:
            \begin{itemize}
                \item $\ov{\Sigma} \setminus \ov{\b \alpha}$ and $ \ov{\Sigma} \setminus  \b \beta$ are connected;
                
                \item $\ov{\b \alpha}$ intersect $\b \beta$ transversely.
            \end{itemize}
            \end{definition}

            \begin{figure}[h]
            \centering
            \includegraphics[width=0.4\textwidth]{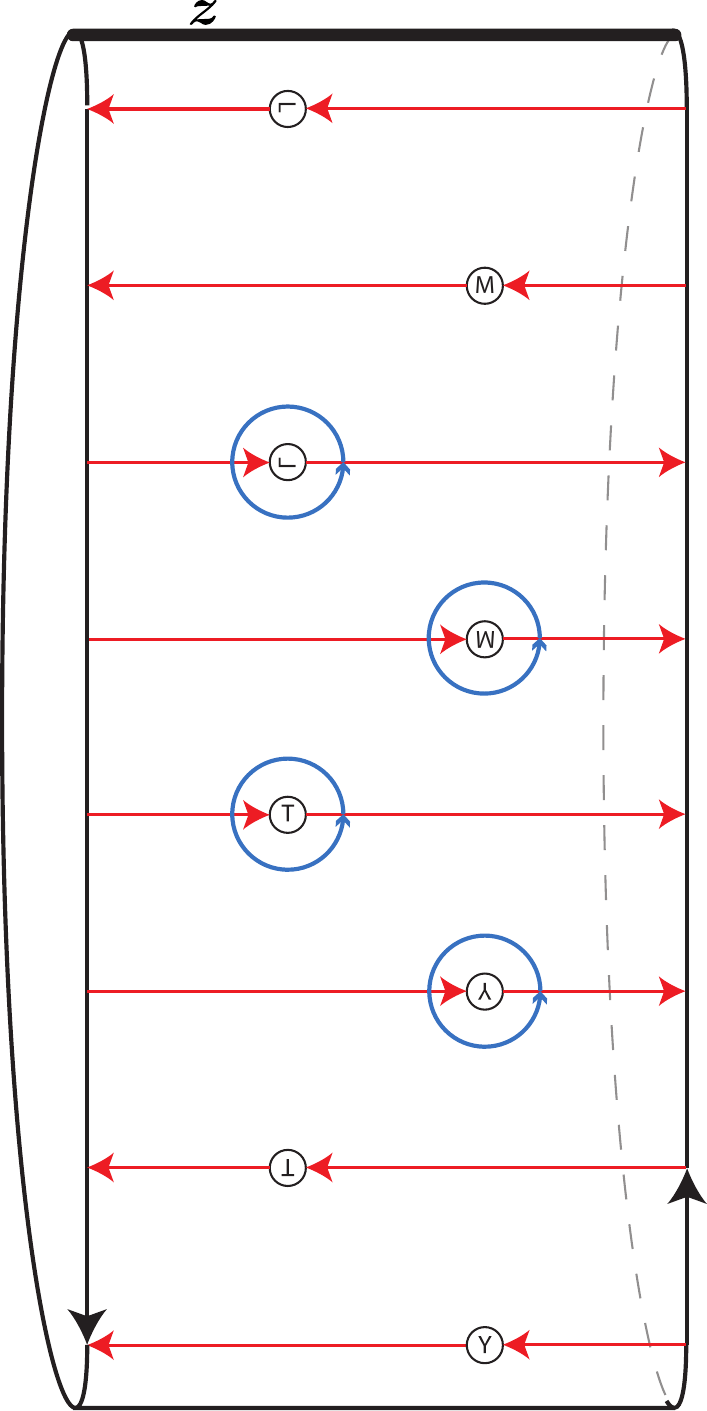}
            \caption{A Heegaard diagram for $M_{\id}$, where $\id\co F^\circ(\Zpmc_2) \rightarrow F^\circ(\Zpmc_2)$ is the identity mapping class of the genus two surface. The pointed matched circle $\Zpmc_2$ is the one from Figure~\ref{fig:pointed_matched_circle}. We also indicate here the orientations of the $\alpha$ and $\beta$ curves, because later they will give $\mathbb{Z}_2$-grading on Hochschild homology.}
            \label{fig:heegard_diagram_id}
            \end{figure}

            Notice that two boundaries of any arced bordered Heegaard diagram specify two pointed matched circles. In Figure~\ref{fig:heegard_diagram_id} we provide an example of a Heegaard diagram of the mapping cylinder of $\id\co\Sigma_2 \rightarrow \Sigma_2$ in the genus two case.

            The following proposition provides the second correspondence in Figure~\ref{fig:plan_2}.

            \begin{proposition}[\text{\cite[Construction~5.6,Propositions~5.10,~5.11]{LOT-bim}}] \label{moves_thm}
            Any arced bordered Heegaard diagram with two boundary components gives rise to a strongly bordered $3$-manifold. For the other direction, suppose a strongly bordered $3$-manifold has a boundary parameterized by $F(\Zpmc_1)$ and $F(\Zpmc_2)$. Then this $3$-manifold has an arced bordered Heegaard diagram with boundary pointed matched circles $\Zpmc_1$ and $\Zpmc_2$. The diffeomorphism type of this Heegaard diagram is unique up to certain moves that transform one diagram to another.
            % \begin{itemize}
            %     \item Isotopies of $\alpha$- and $\beta$- curves and arcs;

            %     \item Handleslides among the $\alpha$-circles and among the $\beta$-circles;

            %     \item Handleslides of an $\alpha$-arc over an $\alpha$-circle;

            %     \item Stabilization of the diagram.
            % \end{itemize}
            \end{proposition}
            \noindent
            In one direction, the procedure of getting a strongly bordered $3$-manifold from an arced bordered Heegaard diagram consists of thickening the surface, attaching $2$-handles along the circles (along the $\beta$-circles from one side, and along the $\alpha$-circles from the other side), and then carefully analyzing what happens on the $\alpha$ side of the boundary. There we have two surfaces of genera $l$ and $r$ with parameterizations coming from the Heegaard diagram, which are connected by an annulus, with the path $\b z$ in it. This annulus, along with the path $\b z$ in it, specify the framed arc $\gamma_z$ by which the two boundary surfaces are connected in the definition of a strongly bordered $3$-manifold. In the other direction, the existence of a Heegaard diagram follows from Morse theory, see~\cite[Proposition~5.10]{LOT-main}. To obtain uniqueness, Heegaard diagrams connected by certain moves are set to be equivalent, and we refer the reader to~\cite[Proposition~5.11]{LOT-bim} for the list of those moves and the proof.
        \subsubsection{Bimodules} \label{sec:hf_bim}
            Here the original definition of the bimodule invariant is described, following~\cite{LOT-bim}. Let us note that there is an alternative, more elementary approach, studied in~\cite{LOT-mcg} and~\cite{Sie}. It is similar in spirit to the Fukaya categorical interpretation of the bimodule invariant, which is described in Section~\ref{sec:bimodule_via_Fukaya}.

            As a prerequisite to the subsequent material, we refer the reader to the excellent sources~\cite[Section~2]{LOT-main} and~\cite[Section~2]{LOT-bim} for the algebraic theory of \emph{$dg$ algebras}, \emph{$A_\infty$ modules}, and \emph{$A_\infty$ bimodules}. In particular, the most important algebraic object for us will be a \emph{type $DA$ bimodule over $dg$ algebras}, see \cite[Definition~2.2.43]{LOT-bim}.  Following~\cite{LOT-bim}, we will use subscripts and superscripts to indicate the algebras of type $DA$ bimodules: if algebra $\mathcal A_1$ is on the $D$-side and algebra $\mathcal A_2$ is on the $A$-side of type $DA$ bimodule $N$, then we will write ${}^{\mathcal A_1}N_{\mathcal A_2}$. We will also  occasionally make use of \emph{type $AA$} and \emph{type $DD$ bimodules}, using notations ${}_{\mathcal A_1}N_{\mathcal A_2}$ and ${}^{\mathcal A_1}N^{\mathcal A_2}$ respectively. 

\begin{figure}[ht]
            \centering
                \includegraphics[width=0.7\textwidth]{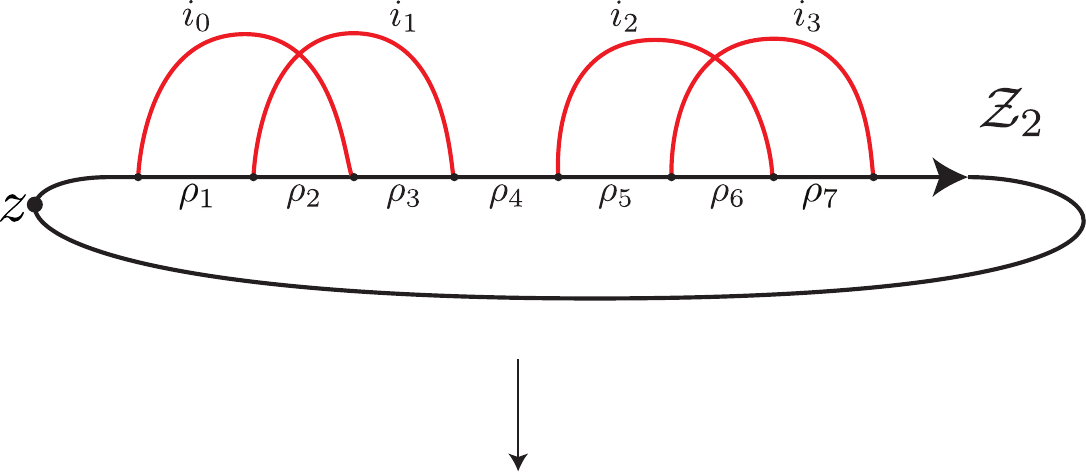}
                \tikzstyle{arrow} = [thick,->,>=stealth]

                \begin{tikzpicture}[scale=1,baseline=1.5cm]
                \node[circle](i_0) at (0,0){$i_0$};
                \node[circle](i_1) at (3,0){$i_1$};
                \node[circle](i_2) at (5,0){$i_2$};
                \node[circle](i_3) at (8,0){$i_3$};

                \draw[arrow] (i_0) to[bend left] node[above, sloped]{$\rho_1$} (i_1);
                \draw[arrow] (i_1) to node[above, sloped]{$\rho_2$} (i_0);
                \draw[arrow] (i_0) to[bend right] node[above, sloped]{$\rho_3$} (i_1);
                \draw[arrow] (i_2) to[bend left] node[above, sloped]{$\rho_5$} (i_3);
                \draw[arrow] (i_3) to node[above, sloped]{$\rho_6$} (i_2);
                \draw[arrow] (i_2) to[bend right] node[above, sloped]{$\rho_7$} (i_3);
                \draw[arrow] (i_1) to node[above, sloped]{$\rho_4$} (i_2);
                
                \draw (0,1) to[bend right] node[left,text width=3cm, text centered]{Path algebra \\ over $\mathbb{F}_2$} (0,-1);
                \draw (8,1) to[bend left] (8,-1);
                \draw (-2.5,-1.4) to (8,-1.4);
                \node at (3,-2){Relations $\rho_1 \rho_4= \rho_4 \rho_7= \rho_2 \rho_1=\rho_3 \rho_2=\rho_6 \rho_5=\rho_7 \rho_6=0$};
                \node[scale=1.2] at (-4.5,-1.4){$\B(\Zpmc_2)=$};

                \end{tikzpicture}
            \caption{A genus two example of how to associate a $dg$ algebra to a pointed matched circle.}
            \label{fig:Z_2--B-Z_2-}
            \end{figure}

            The first step in constructing any bordered Heegaard Floer invariants is always specifying the algebra.
            \begin{construction}[The $dg$ algebra associated to a pointed matched circle]\label{construction:dg-algebra}
            To a pointed matched circle $\Zpmc$ we associate a $dg$ algebra $\B(\Zpmc)$, which is $\A(\Zpmc,-g+1)$ in the notation of~\cite{LOT-bim} (i.e. we have only \say{one strand moving}). First, given a pointed matched circle $\Zpmc$, we construct a directed graph (\emph{quiver} for short) $\Gamma(\Zpmc)$. The vertices of $\Gamma(\Zpmc)$ are the matched pairs of points in the pointed matched circle, and the edges are the length one chords between the points, which do not cross the basepoint $z$. We illustrate this in Figure~\ref{fig:Z_2--B-Z_2-}, where we draw an example of a pointed matched circle $\Zpmc_2$ and the corresponding quiver $\Gamma(\Zpmc_2)$ below, inside the parenthesis.

            To the quiver $\Gamma(\Zpmc)$ we associate a path algebra ${\mathbb F}_2(\Gamma(\Zpmc))$. The generators of the underlying $\mathbb{F}_2$-vector space are all the paths of the quiver (including the constant paths). If edges are denoted by letters $\rho_j$, we use the notation $\rho_{j_1 \dots j_k}$ for a path $(\rho_{j_1}, \ldots, \rho_{j_k})$. The multiplication in the algebra is given by concatenating the paths, if possible: for example, in ${\mathbb F}_2(\Gamma(\Zpmc))$ from Figure~\ref{fig:Z_2--B-Z_2-}, we have $i_1 \cdot \rho_4 \cdot \rho_5=\rho_{45}$. If paths do not concatenate, we declare the product to be zero. 

            At last, to obtain the algebra $\B(\Zpmc)$ we quotient the algebra ${\mathbb F}_2(\Gamma(\Zpmc))$ by setting certain paths to be zero. Those are the paths in which not consecutive chords are concatenated. For example, in ${\mathbb F}_2(\Gamma(\Zpmc))$ from Figure~\ref{fig:Z_2--B-Z_2-}, despite the fact that $\rho_{32}$ is a valid path, $\rho_2$ does not go right after $\rho_3$ in the pointed matched circle $\Zpmc_2$, and thus we set $\rho_{32}=0$. The full set of relations for $\B(\Zpmc_2)$ is given in the bottom of Figure~\ref{fig:Z_2--B-Z_2-}. The differential in $dg$ algebras $\B(\Zpmc)$ is always set to be trivial. All the constant paths of $\Gamma(\Zpmc)$ are idempotents in the algebra $\B(\Zpmc)$, and they correspond to the matched pairs of points in $\Zpmc$, i.e. $1$-handles of $F^\circ(\Zpmc)$. We usually denote these idempotents by $i_k$; there are four of them in $\B(\Zpmc_2)$ from Figure~\ref{fig:Z_2--B-Z_2-}. The sum of all the idempotents is a unit in $\B(\Zpmc)$ .
            \end{construction}

            Now we turn to the construction of type $DA$ bimodule invariant associated to a mapping class.
            \begin{construction}[The type $DA$ bimodule associated to a mapping class]\label{constr_bim}
            Fix a surface with one boundary component $(\Sigma,\partial \Sigma = S^1)=F^\circ(\Zpmc)$. A mapping class $\phi \in MCG_0(\Sigma,\partial \Sigma)= MCG_0(\Zpmc,\Zpmc)$ gives rise to a strongly bordered $3$-manifold $M_\phi$, the mapping cylinder of $\phi$. After this we can consider a Heegaard diagram $\H(M_\phi)$ representing $M_\phi$, with pointed matched circles $\Zpmc$ and $\Zpmc$ on the boundaries. To such a Heegaard diagram we associate a type $DA$ bimodule ${}^{\B(\Zpmc)}N(\phi)_{\B(\Zpmc)}$, over the algebra $\B(\Zpmc)$ from both sides. 

            The generators of the underlying $\mathbb{F}_2$-vector space of the bimodule consist of tuples $\x$ of intersection points between $\alpha$- and $\beta$- curves such that every $\alpha$ and $\beta$ circle gets one point, only one $\alpha$-arc on the right boundary gets a point (denote this $\alpha$-arc by $\alpha^\text{arc,right}_{\x}$), and all except one $\alpha$-arc on the left boundary get one point (denote this $\alpha$-arc by $\alpha^\text{arc,left}_{\x}$). See an example in Figure~\ref{fig:heeg_diag_dehn_twist}, where we marked all the generators on a Heegaard diagram.

            The idempotent subalgebra of $\B(\Zpmc)$ acts on these generators in the following way. By Construction~\ref{construction:dg-algebra} for an $\alpha$-arc there is an associated idempotent of $\B(\Zpmc)$, which we denote by $i(\alpha^{arc})$. For a generator $\x$, we have actions $i(\alpha^\text{arc,left}_{\x}) \cdot \x= \x $, $\x \cdot i(\alpha^\text{arc,right}_{\x})= \x $, and the other idempotent actions are zero.

            Higher $A_\infty$ type $DA$ actions on these generators are defined by counting of pseudo-holomorphic curves; for this, analytic in nature, definition we refer the reader to~\cite[Section~6.3]{LOT-bim} and references therein. Let us note an important point: along the way of the definition, an analytic choice of a family of almost complex structures on some space has to be made. If one makes a different analytic choice, or chooses a different Heegaard diagram for $M_\phi$, the resulting type $DA$ bimodule will be the same up to $A_\infty$ homotopy equivalence of type $DA$ bimodules. 
            \end{construction}

            This finishes the explanation of how, to a mapping class $\phi \in MCG_0(\Sigma,\partial \Sigma)=MCG_0(\Zpmc,\Zpmc)$, we can associate an $A_\infty$ homotopy equivalence class of type $DA$ bimodules ${}^{\B(\Zpmc)}N(\phi)_{\B(\Zpmc)}$. An important feature of this mapping class invariant is that there is an operation  on bimodules which corresponds to multiplication in the mapping class group. This operation is called \emph{box tensor product}~\cite[Definitions~2.3.2 and~2.3.9]{LOT-bim}, and is denoted by $\boxtimes$.
            \begin{theorem}[\text{\cite[Theorem~12]{LOT-bim}}]\label{DA pairing}
            Suppose $\phi_1$ and $\phi_2$ are two elements in the mapping class group $MCG_0(\Zpmc,\Zpmc)$. Then we have the following homotopy equivalence of bimodules: 
            $${}^{\B(\Zpmc)}N(\phi_1 \phi_2)_{\B(\Zpmc)} \simeq {}^{\B(\Zpmc)}N(\phi_2 )_{\B(\Zpmc)}      \boxtimes {}^{\B(\Zpmc)}N(\phi_1)_{\B(\Zpmc)}.$$
            \end{theorem}
            % \noindent The high level strategy for the proof is as follows. It is clear that the mapping cylinder of the product $\phi_1 \phi_2$ can be obtained by gluing the two mapping cylinders: $M_{\phi_1 \phi_2} \cong M_{\phi_2}  {}_{\partial_R}{\cup}_{\partial_L}  M_{\phi_1}$. The same is true for Heegaard diagrams: $\H_{\phi_1 \phi_2} \cong \H_{\phi_2}  {}_{\partial_R}{\cup}_{\partial_L}  \H_{\phi_1}$. Thus, the content of~\cite[Theorem~12]{LOT-bim} is that gluging Heegaard diagrams corresponds to taking the box tensor product on the bimodule side; the proof involves deep analysys of pseudo-holomorphic curves.
            % Let us also point out that the box tensor product ooperation is, in fact, quite natural: it turns out that $\boxtimes$ is a model for the $A_\infty$ tensor product in the case of pairing the $D$ side to the $A$ side, see~\cite[Sections~2.2 and~2.4]{LOT-main} for the explanation of that.
            
            Along the way of associating the bimodule invariant to $\phi \in MCG_0(\Sigma,\partial \Sigma)=MCG_0(\Zpmc,\Zpmc)$, we made a choice of a particular pointed matched circle $\Zpmc$, and an identification of $F^\circ(\Zpmc)$ with $(\Sigma, \partial \Sigma)$. It turns out that for us these choices are not important. Given two different pointed matched circles $\Zpmc, \ \Zpmc'$, and two different parameterizations of the surface $(\Sigma, \partial \Sigma) \cong_1 F^\circ (\Zpmc)  $ and $ (\Sigma, \partial \Sigma) \cong_2 F^\circ (\Zpmc
            ') $, 
            % (note that if pointed matched circles are equal, it does not mean that the identifications are the same)
            the two mapping class groups $MCG_0(F^\circ (\Zpmc),F^\circ (\Zpmc
            ))$ and $MCG_0(F^\circ (\Zpmc'),F^\circ (\Zpmc
            '))$ can be bijectively identified via conjugation by some element $a \in MCG_0(F^\circ (\Zpmc'),F^\circ (\Zpmc
            ))$: 
            \begin{align*}
            MCG_0(F^\circ (\Zpmc'),F^\circ (\Zpmc
            ')) &\rightarrow MCG_0(F^\circ (\Zpmc),F^\circ (\Zpmc
            )) \\
            \psi&\mapsto a \psi a^{-1}
            \end{align*}
            Conveniently, the bimodules are also related: if $\phi=a \psi a^{-1}$, then, according to Theorem~\ref{DA pairing}, we have
            $$N(\phi) \simeq N(a^{-1}) \boxtimes N(\psi) \boxtimes N(a).$$ 
            More importantly, the main invariant for us is going to be the Hochschild homology of a bimodule, which is invariant with respect to conjugation of a mapping class.

            Lastly, we quote~\cite[Theorem~4]{LOT-bim}, which says that ${}^{\B(\Zpmc)}N(\id)_{\B(\Zpmc)} \simeq {}^{\B(\Zpmc)}[\mathbb{I}]_{\B(\Zpmc)}$ is the \emph{identity type $DA$ bimodule}, i.e. it has generators $_{i_k}{i_k}_{i_k}$ for every idempotent of $\B(\Zpmc)$, and actions $({}_{i_k}{i_k}_{i_k},a) \rightarrow a \o {}_{i_l}{i_l}_{i_l}$ for every element $i_k  a i_l = a$ in the algebra $\B(\Zpmc)$; see~\cite[Definition 2.2.48]{LOT-bim}. From the definition it follows that the operations $P\boxtimes{}^{\B(\Zpmc)}[\mathbb{I}]_{\B(\Zpmc)}$ and ${}^{\B(\Zpmc)}[\mathbb{I}]_{\B(\Zpmc)} \boxtimes P$ do not affect the bimodule $P$, up to isomorphism; this will be important later.

            \subsubsection{Hochschild homology}\label{sec:hoch}
            One of the main goals of this paper is to relate the bimodule $N(\phi)$ to fixed point Floer cohomology (defined in Section~\ref{sec:fixed point Floer}). In turns out that for that we need to apply a certain algebraic operation to the bimodule, which is called \emph{Hochschild homology} and denoted by $HH_*(N(\phi))$. In general, it is a homology theory associated to type $AA$ bimodules ${}_{\mathcal A_1}N_{\mathcal A_2}$ and type $DA$ bimodules ${}^{\mathcal A_1}N_{\mathcal A_2}$ that satisfy $\mathcal A_1=\mathcal A_2$. Hochschild homology should be though of as self-tensoring of the bimodule, pairing the left side to the right side. We refer to~\cite[Section~2.3.5]{LOT-bim} for the algebraic definitions, basic properties, and a way to compute Hochschild homology for bounded type $DA$ bimodules.

            There are two important points about this algebraic structure. First, Hochschild homology depends only on the $A_\infty$ homotopy equivalence class of the bimodule. Thus, the Hochschild homology of the bimodule $HH_*(^{\B(\Zpmc)}N(\phi)_{\B(\Zpmc)})$ gives us an invariant of a mapping class. Second, $HH_*(^{\B(\Zpmc)}N(\phi)_{\B(\Zpmc)})$ can be identified with knot Floer homology.
            \begin{theorem}[\text{\cite[Theorem~7]{LOT-bim}}]
            Given a mapping class $\phi \in MCG(\Sigma,\partial \Sigma)$, denote by $M_\phi^\circ$ the open book with monodromy $\phi\co (\Sigma,\partial \Sigma) \rightarrow (\Sigma, \partial \Sigma)$ and binding $K$. Then the following two vector spaces are isomorphic:
            $$HH_*(N(\phi))\cong \widehat{HFK}(M_\phi^\circ,K \:;\:  -g+1),$$
            where the latter is the Alexander grading $(-g+1)$ summand of knot Floer homology of $K$.
            \end{theorem}
            % Let us explain the reasons behind this theorem. The open book $M_\phi^\circ$  is obtained from the mapping cylinder $M_\phi$ (Construction~\ref{mapping cylinder}) by (1) identifying the two boundaries (turning the arc  $\gamma_z$ into a circle); and (2) doing a surgery on $\gamma_z$ with respect to the framing in the definition of $M_\phi$. 

            % The construction of a Heegaard diagram for this open book is as follows. First we take an arced bordered Heegaard diagram $\H (M_\phi)$ for the mapping cylinder $M_\phi$, and glue the left boundary to the right one (on the algebraic level this corresponds to self-tensoring the bimodule). Then we do a surgery on the arc $\b z$ (which became a circle after self-gluing). Note that in order to block those discs which we did not count, we need to place two basepoints on the two sides after surgery --- these two basepoints specify a knot, which is the binding of the open book. 

            % Up to the self-pairing procedure, this explains why $\widehat{HFK}(M_\phi^\circ,K \:;\:  -g+1)\cong HH_*(N(\phi))$, see~\cite[Theorem~14]{LOT-bim} for the details. Because we work in the second lowest $\text{Spin}^c$ structure, we only get the Alexander grading $(-g+1)$ summand.

            \begin{remark}\label{rmk:inv_conj}
            The open book $M_\phi^\circ$ is invariant with respect to conjugation of $\phi$, which implies that $HH_*(N(\phi))$ is invariant with respect to conjugation of $\phi$.
            \end{remark}

            There is a $\mathbb{Z}_2$-grading on $\widehat{HFK}(M_\phi^\circ,K \:;\:  -g+1)$, coming from the sign of intersections of tori $\mathbb{T}_\alpha$ and $\mathbb{T}_\beta$ in the definition of generators of knot Floer homology.  In conjunction with the theorem above, this endows Hochschild homology with a $\mathbb{Z}_2$-grading. 

            % Because the Heegaard diagram for the open book is constructed via self-gluing, the $\mathbb{Z}_2$-grading, i.e. the choice of orientations on $\mathbb{T}_\alpha$, $\mathbb{T}_\beta$, amounts to the choice of orientations on $\beta$ and $\alpha$ curves on the Heegaard diagram for $M_\phi$ (s.t. they are consistent, i.e. after gluing the orientations of the left arcs match with the orientations of the right arcs). Let us choose these orientations as in Figure~\ref{fig:heegard_diagram_id}. To obtain a consistent $\mathbb{Z}_2$-grading on $HH_*(N(\phi))$ in general case for $\phi\co MCG_0(\Sigma,\partial\Sigma)$, we make such choices of orientations in the standard Heegaard diagram for $M_\phi$ (see~\cite[Section~5.3]{LOT-bim} or the next Section~\ref{subsec:bim_comps} for a construction of the standard $\H ( M_\phi)$), that the gradings of all $2g$ generators of $HH_*(N(\id))$ are zero.
            \subsubsection{Cancellation} 
            Let us finish this section by describing a process called \emph{cancellation}; for details we refer the readers to~\cite[Section~2.6]{Levine} (where it is called the ‘‘edge reduction algorithm'') and~\cite[Section~3.1]{Boh}. Suppose there are two generators $_i x_j$ and $_i y_j$ (the subscripts indicate their left and right idempotents) in a type $DA$ bimodule $P$ such that the only action between them is $\delta^1_1(_i x_j)=i \o _i y _j$ (we will label such actions by $1$, see the arrow between two generators $x_2$ and $t_{12}$ in Figure~\ref{fig:C_inv_bimodule}). Then we can cancel these two generators, i.e. erase $x$ and $y$ and the arrows involving them from the bimodule, and then add some other arrows between the remaining generators in the bimodule, guided by a cancellation rule: for every \say{zigzag} configuration of arrows 
            $$z_2 \xleftarrow{b \o (d_1,\ldots,d_l)} x \xrightarrow{1} y  \xleftarrow{a \o (c_1,\ldots,c_k)} z_1$$
            in the initial bimodule $P$ we will add an arrow
            $$z_2 \xleftarrow{a\cdot b \o (c_1,\ldots,c_k,d_1,\ldots,d_l)} z_1$$
            The outcome is a new bimodule $P'$ with less generators, which is homotopy equivalent to the previous one $P'\simeq P$. 

    \subsection{New bimodule computations}\label{subsec:bim_comps}
        In this section, we compute the type $DA$ bimodules $N(\phi)$ for mapping classes of the genus two surface; the computations are based on the arc-slide bimodules studied in~\cite{LOT-arc}.

\begin{figure}[H]
        \centering
        \includegraphics[width=0.6\textwidth]{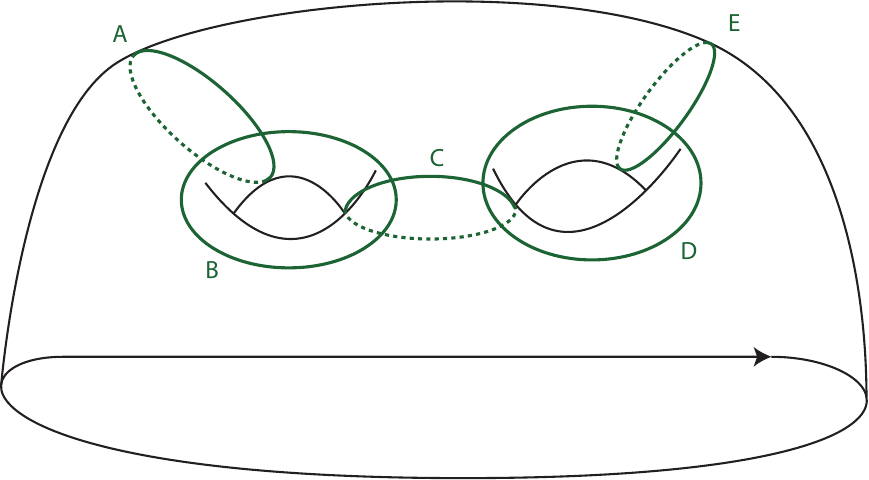}
        \caption{Dehn twists along these curves generate $MCG_0(\Sigma_2)$.}
        \label{fig:Sigma_2}
        \end{figure}

        First, fix a genus two surface $\Sigma_2$ with one boundary component, and a set of curves on it as in Figure~\ref{fig:Sigma_2}.
        The following is a presentation of the mapping class group of the genus two surface with one boundary component (see~\cite[Theorem~2]{W}): 
        \begin{equation}
        \begin{aligned}\label{g2_mcg_presentationg}
        MCG_0(\Sigma_2,\partial \Sigma_2)= \langle \tau_A,\tau_B,\tau_C,\tau_D,\tau_E | \ \ &\text{commuting  relations}, \ \ \text{braid relations}, \\
        &(\tau_E \tau_D \tau_C \tau_B)^5=\tau_A\tau_B\tau_C\tau_D\tau_E \tau_E\tau_D\tau_C\tau_B \tau_A \rangle,
        \end{aligned}
        \end{equation}
        where $\tau_l$ is a right handed Dehn twist along the curve $l$. By commuting relations we mean that, if curves $l_1$ and $l_2$ do not intersect, then $\tau_{l_1} \tau_{l_2} =\tau_{l_2}  \tau_{l_1}$. By braid relations we mean that, if curves $l_1$ and $l_2$ intersect at a single point transversely, then $\tau_{l_2} \tau_{l_1} \tau_{l_2} = \tau_{l_1} \tau_{l_2} \tau_{l_1}$. 

        Now we explain how to compute $N(\phi)$ for any $\phi \in MCG_0(\Sigma_2,\partial \Sigma_2)$. Every such mapping class can be represented by a product of Dehn twists $\tau_A,\tau_B,\tau_C,\tau_D,\tau_E$, or their inverses. Theorem~\ref{DA pairing} governs how bimodules $N(\phi)$ behave with respect to composition of mapping classes: the corresponding operation is the box tensor product. Thus it is enough to compute $N(\tau)$ only for these Dehn twists and their inverses:
        \begin{equation}\label{ten_bimodules}
        N(\tau_A),\ N(\tau_B),\ N(\tau_C),\ N(\tau_D),\ N(\tau_E), N(\tau^{-1}_A),\ N(\tau^{-1}_B),\ N(\tau^{-1}_C),\ N(\tau^{-1}_D),\ N(\tau^{-1}_E)
        \end{equation}

        \begin{remark}
        There is a very interesting strategy, pioneered by Zhan~\cite{Boh}, to assign Floer theoretic invariants to topological objects without referring to pseudo-holomorphic theory. Assuming the ten bimodules~\eqref{ten_bimodules} are computed (we will proceed to computing them after the remark), this strategy can be applied in our case: we can redefine the bimodule invariant $N(\phi)$ for $\phi \in MCG_0(\Sigma_2,\partial \Sigma_2)$ in the following combinatorial way. For a factorization of $\phi$ into ten Dehn twists~\eqref{ten_bimodules} --- let us assume it is $\phi=\tau_A \tau_C \tau^{-1}_E$, for example, --- we associate a tensor product of the corresponding bimodules $N(\tau^{-1}_E) \boxtimes N(\tau_C) \boxtimes N(\tau_A)$. To make sure it is an invariant of the mapping class group element up to conjugation, we need to check the mapping class group relations. For example for the relation $\tau_{A} \tau_{B} \tau_{A} = \tau_{B} \tau_{A} \tau_{B}$ we need to check the following homotopy equivalence: $N(\tau_A) \boxtimes N(\tau_B) \boxtimes N(\tau_A) \simeq N(\tau_B) \boxtimes N(\tau_A) \boxtimes N(\tau_B)$. Indeed, in our case, after the computation of ten bimodules~\eqref{ten_bimodules}, computer can show that all the relations in presentation~\eqref{g2_mcg_presentationg} are satisfied; see~\cite[Showcase~3]{Pyt} for an illustration. 

        In~\cite{Boh}, Zhan gives a combinatorial definition of $\widehat{CFDA}(\phi, 0)$; it is an analogue of our invariant $\widehat{CFDA}(\phi, -g+1)=N(\phi)$, where generators occupy $g$ arcs on the left and $g$ arcs on the right boundary of the Heegaard diagram. In the definition he uses arc-slides (as opposed to Dehn twists) as generators of the mapping class groupoid. Using this definition for $\widehat{CFDA}(\phi,0)$, he then gives a combinatorial definition of the \say{hat} version of Heegaard Floer homology of a $3$-manifold $\widehat{HF}(Y^3)$.
        \end{remark}

\begin{figure}[ht]
        \centering
        \includegraphics[width=0.8\textwidth]{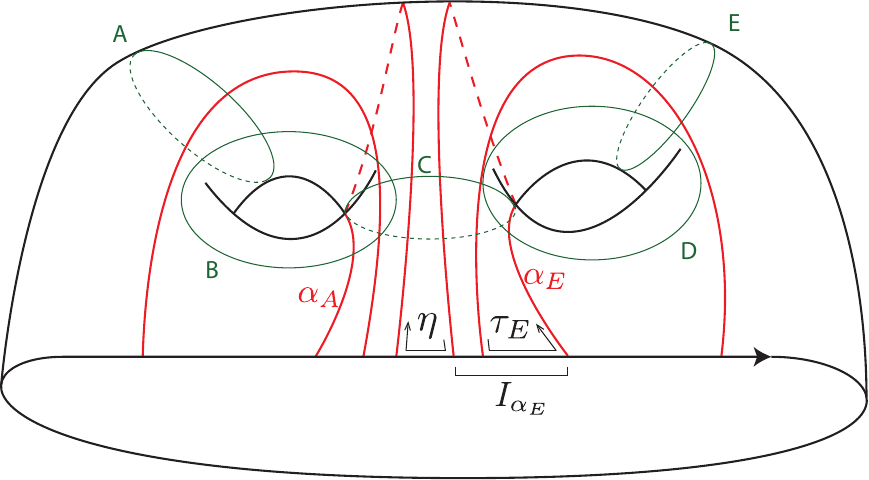}
        \caption{Parameterization of the surface $\Sigma_2 \cong F^\circ(\Zpmc_2)$.}
        \label{fig:Sigma_2_parameterized}
        \end{figure}

        We now compute the ten bimodules~\eqref{ten_bimodules}. First we need to fix a parameterization of our surface $(\Sigma_2,\partial \Sigma_2) \cong F^\circ(\Zpmc)$. This will specify a $dg$ algebra. We use the pointed matched circle $\Zpmc_2$ and its corresponding algebra $\B(\Zpmc_2)$ from Figure~\ref{fig:Z_2--B-Z_2-}. For an identification $(\Sigma_2,\partial \Sigma_2) \cong F^\circ(\Zpmc_2)$ see Figure~\ref{fig:Sigma_2_parameterized}.
        % first note that right shaded surface in Figure~\ref{fig:curves_on_heeg_diag} is naturally identified with $F^\circ(\Zpmc_2) \setminus D^2$ (identification is specified by red arcs, which are $1$-handles). Now use a map $(\Sigma_2,\partial \Sigma_2) \rightarrow F^\circ(\Zpmc_2)$ such that curves $A,B,C,D,E$ are mapped to $A_r,B_r,C_r,D_r,E_r$.

\begin{figure}[ht]
        \includegraphics[width=0.6\textwidth]{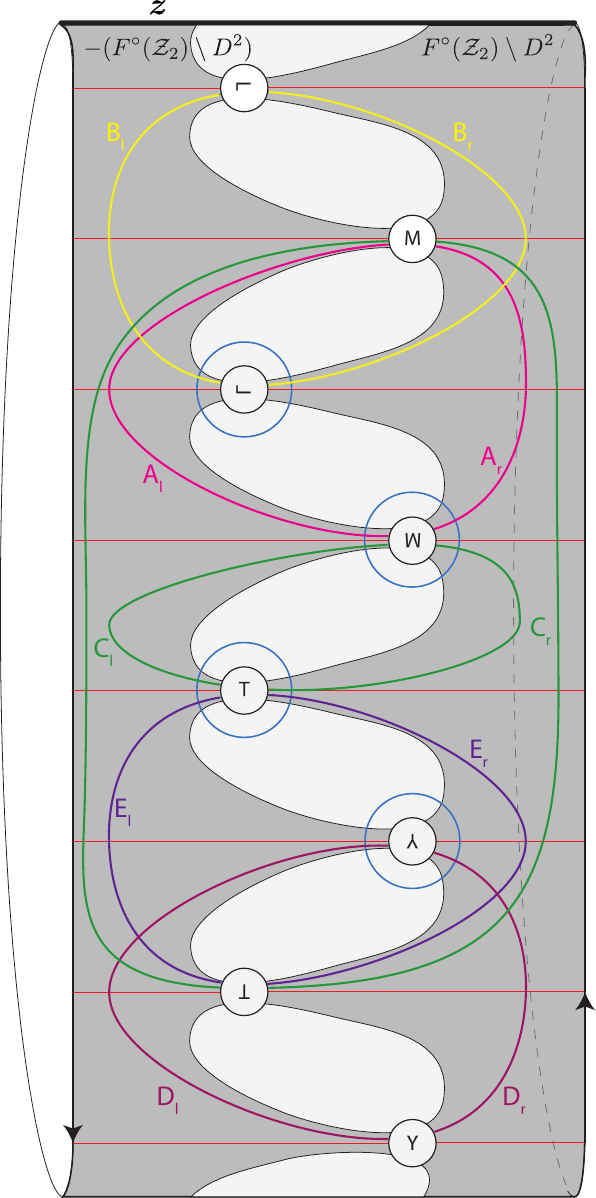}
        \caption{Heegaard diagram $\H(M_{\id})$ of identity mapping class with curves over which we do the Dehn twists.}
        \label{fig:curves_on_heeg_diag}
        \end{figure}

        For every Dehn twist $\tau_l \in MCG_0(\Sigma_2,\Sigma_2)$ we need to specify a Heegaard diagram for a mapping cylinder $M_{\tau_l}$. Following~\cite[Section~5.3]{LOT-bim}, consider first the standard Heegaard diagram $\H(M_{\id})$ for $\id\co F^\circ(\Zpmc_2) \rightarrow F^\circ(\Zpmc_2)$, see Figure~\ref{fig:curves_on_heeg_diag}. There is a shaded region of the diagram on the right that is identified with the right boundary $F^\circ(\Zpmc_2)\setminus D^2$ of the mapping cylinder. Analogously, there is a shaded region on the left part of the diagram which is identified with $-(F^\circ(\Zpmc_2)\setminus D^2)$. There are also curves $A_l,~B_l,~C_l,~D_l,~E_l$ on the left surface, and $A_r,~B_r,~C_r,~D_r,~E_r$ on the right surface, via the specified above identification $F^\circ(\Zpmc_2) \cong \Sigma_2$.

        % Important point about sign conventions: in the procedure of getting mapping cylinder from Heegaard diagram there is a technical orientation change of the boundary surfaces. Thus identification of shaded regions with boundaries are orientation reversing. This explains the next commented paragraph.

\begin{figure}[!ht]
            \begin{subfigure}{0.5\textwidth}
              \centering
              \includegraphics[width=0.58\linewidth]{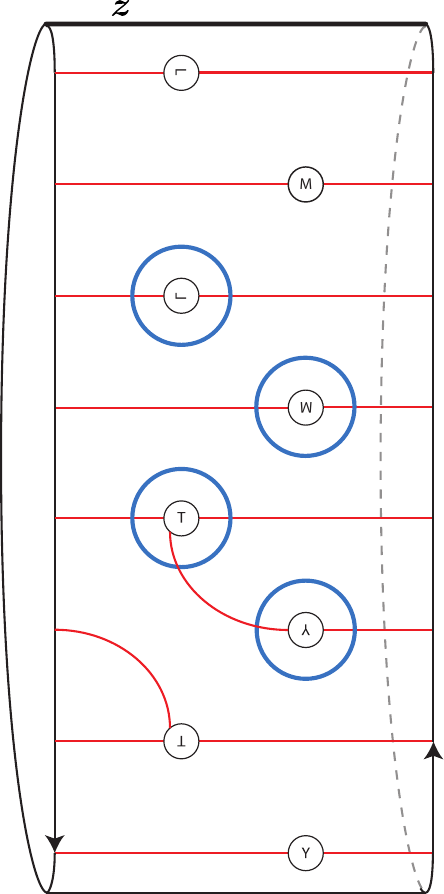}
              \caption{1st type.}
            \end{subfigure}%
            \begin{subfigure}{0.5\textwidth}
              \centering
              \includegraphics[width=0.58\linewidth]{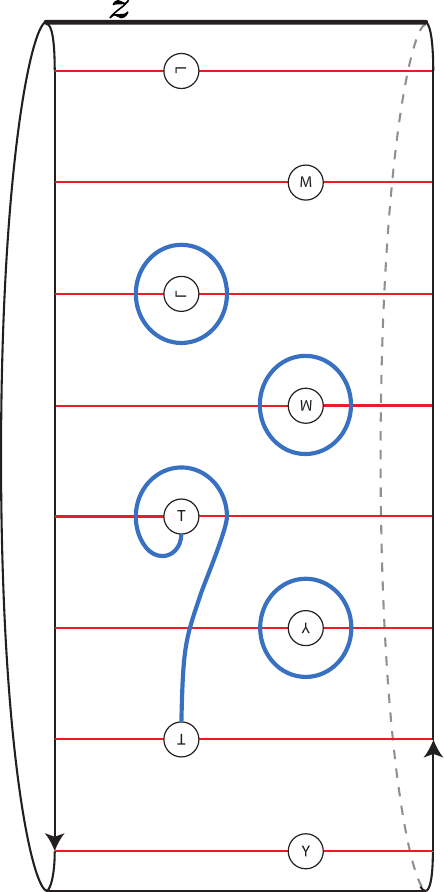}
              \caption{2nd type.}
            \end{subfigure}\\
            \vspace*{1cm}
            \begin{subfigure}{0.5\textwidth}
              \centering
              \includegraphics[width=0.58\linewidth]{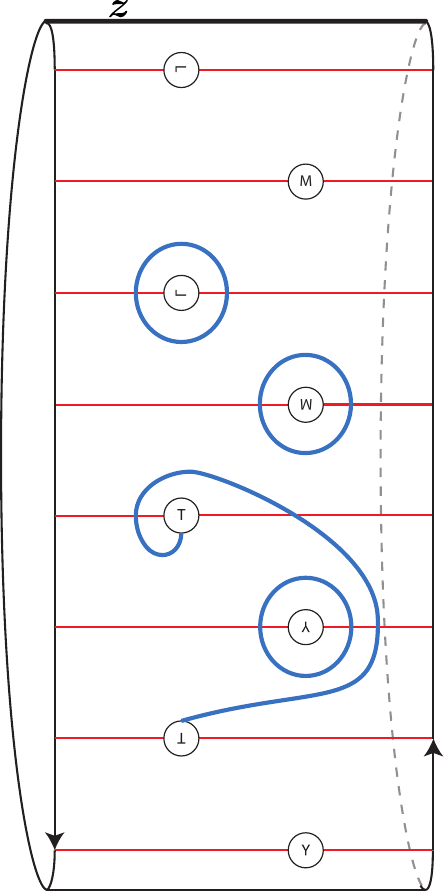}
              \caption{3rd type.}
            \end{subfigure}%
            \begin{subfigure}{0.5\textwidth}
              \centering
              \includegraphics[width=0.58\linewidth]{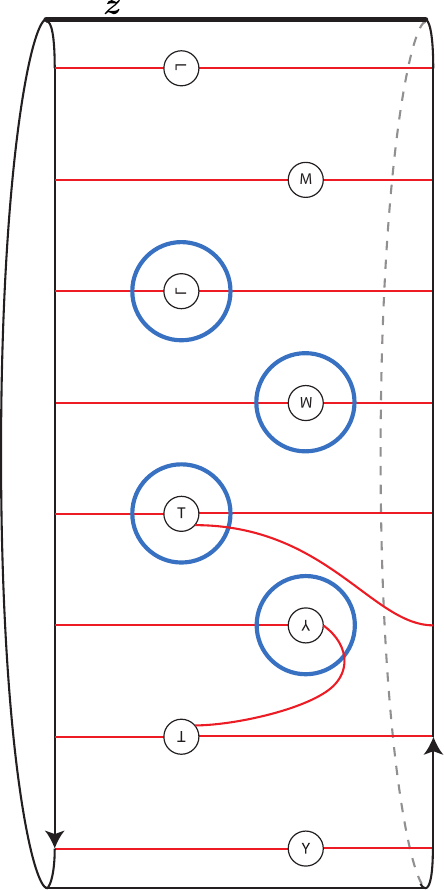}
              \caption{4th type.}
            \end{subfigure}%
        \caption{Four Heegaard diagrams for a Dehn twist along the curve E.}
        \label{fig:heeg_diag_dehn_twists_1-4}
        \end{figure}

        Now, suppose we want to draw a Heegaard diagram for $M_{\tau_E}$. Then it is enough to change all the red arcs of the left side of $\H(M_{\id})$ by applying $\tau_{E}$. This corresponds to a parameterization $-\tau_E:-F^\circ(\Zpmc_2) \rightarrow \partial_L M_{\tau_E}$. So we apply the right handed Dehn twist $\tau_{E_l}$ to the red arcs.  % with respect to the orientation $-F^\circ(\Zpmc_2)$, which coincides with the usual orientation of a plane we are used to, where $x$-coordinate is the first one in basis. 
        Alternatively we can apply $\tau_{E_l}^{-1}$ to all the blue curves (this corresponds to applying self-diffeomorphism $\tau_{E_l}^{-1}$ to the Heegaard diagram). We also could have applied $\tau_{E_r}$ to all the red arcs, or $\tau_{E_r}^{-1}$ to all the blue curves. All these possibilities are depicted in Figure~\ref{fig:heeg_diag_dehn_twists_1-4}.  All of the four diagrams are equivalent up to the equivalence moves from Proposition~\ref{moves_thm}  and self-diffeomorphisms applied to the diagrams. The resulting Heegaard diagrams here are analogous to the ones for the genus one case in~\cite[Section~10.2]{LOT-bim}.

        \begin{remark}
        The orientation convention (essentially the signs of the Dehn twists on the diagrams) is chosen so that the map $\phi$ goes \say{from left to right} on the mapping cylinder and the Heegaard diagram, see~\cite[Appendix A]{LOT-mor}. This ensures the desired behavior with respect to gluing: $\H_{\phi_1 \phi_2} \cong \H_{\phi_2}  {}_{\partial_R}{\cup}_{\partial_L}  \H_{\phi_1}$.
        \end{remark}

        The Heegaard diagrams for Dehn twists along the curves $E,~D,~A,~B$ are very similar in their complexity and structure. In contrast to this, as one can see from Figure~\ref{fig:two_Dehn}, the Heegaard diagrams for $\tau_E$ and $\tau_{C}^{-1}$ are very different, which results in the bimodule  $N(\tau_C^{-1})$ having much more generators than the bimodule $N(\tau_{E})$. Thus we have two fundamentally different computations of the bimodules to perform: for Dehn twists $\tau_E$ and $\tau_{C}$ (or their inverses). Below, we compute the bimodule $N(\tau_E)$ via the second type of Heegaard diagram for $\tau_{E}$, and the bimodule $N(\tau_C^{-1})$ via the third type of diagram for $\tau_{C}^{-1}$.  All other eight bimodules can be computed analogously: computations for $\tau_{E}^{-1},\tau_{D},\tau_{D}^{-1}$ are analogous to the one for $\tau_E$, while all the other five bimodule invariants for $\tau_{A}^{-1},\tau_{A},\tau_{B}^{-1},\tau_{B},\tau_{C}$ can be deduced from the previous five by reflection about the $x$-axis and changing the orientation of the Heegaard diagram. More precisely, to obtain $N(\tau_{A}^{-1})$ from $N(\tau_{E})$ for instance, one has to apply a certain involution to the algebra $\B(\Zpmc_2)$, reverse the direction of arrows, and reverse the order of algebra elements on the $A$ side. The involution is given by $i_k \leftrightarrow i_{3-k}$ and $\rho_k \leftrightarrow \rho_{8-k}$ (in particular, it does not preserve the multiplication). For example, the action 
        $${}_{i_0}(x_0)_{i_0} \o (\rho_{3},\rho_{234}) \rightarrow \rho_{34} \o {}_{i_2}(x_2)_{i_2}\quad \text{ in }N(\tau_{E})$$
         results in action 
        $${}_{i_1}(x_2)_{i_1} \o (\rho_{456}, \rho_{5}) \rightarrow \rho_{45} \o {}_{i_3}(x_0)_{i_3} \quad \text{ in } N(\tau_{A}^{-1}).$$
        All the ten bimodules~\eqref{ten_bimodules} are explicitly described in the program~\cite[Showcase~0]{Pyt}. For a general mapping class we first factorize it into Dehn twists, and then box tensor multiply all the bimodules associated to these Dehn twists.

        While the actual technical computations of bimodules $N(\tau_E),~N(\tau_C^{-1})$ are done in Computations~\ref{comp: dehn twist, the easier one} and~\ref{harder one}, here we first describe the resulting bimodules. Figure~\ref{fig:heeg_diag_dehn_twist} shows the Heegaard diagram $\H (M_{\tau_E})$ along with the marked generators of the bimodule, and also the idempotents corresponding to the arcs. 
        % Notice that we reversed the orientation of the boundary on the left, because it corresponds to the D side of the bimodule, which is over $\B(-(-\Zpmc_2))=\B(\Zpmc_2)$. 
        The type $DA$ bimodule $N(\tau_E)$ is depicted in Figure~\ref{fig:E_bimodule}. The subscripts of the generators represent the underlying left and right idempotents; the arrow between generators $\x$ and $\y$ with the label $a \o (b,c)$ means that there is a type $DA$ action $\x \o (b,  c) \rightarrow a \o \y$; if there is $1$ in the label on the right, it means that there are no incoming algebra elements in that action; if there is $1$ in the label on the left, it means that the corresponding action's outgoing algebra element is an idempotent. The $1$ notation refers to the fact that the sum of the idempotents in the algebra is a unit. 

        Figure~\ref{fig:heeg_diag_dehn_twist_connecting_curve} shows the Heegaard diagram $\H (M_{\tau_C^{-1}})$ along with the marked generators of the bimodule. Because there are many generators, for the generators $t_i$ we only mark an intersection point on the right side of the diagram; its corresponding $2g-1=3$ intersections on the left are uniquely determined. The bimodule $N(\tau_C^{-1})$ is depicted in Figure~\ref{fig:C_inv_bimodule}. On some arrows (which are a little lighter on the picture) we did not write the actions; those actions correspond to all the rectangles in the rectangle area with vertices $t_0,t_{12}$ and the right edge $\rho_{23456}$.

        We now proceed to proving that the bimodules from Figures~\ref{fig:E_bimodule} and~\ref{fig:C_inv_bimodule} are indeed the bimodules $N(\tau_E)$ and $N(\tau_C^{-1})$ coming from the holomorphic structure contained in Heegaard diagrams in Figures~\ref{fig:heeg_diag_dehn_twist} and~\ref{fig:heeg_diag_dehn_twist_connecting_curve}. 

        \begin{computation}[$N(\tau_E)$] \label{comp: dehn twist, the easier one}
            Let us denote for a moment our candidate bimodule from Figure~\ref{fig:E_bimodule} by $N'(\tau_E)$, and the bimodule $N(\tau_E)$ will be the one which corresponds to the Heegaard diagram in Figure~\ref{fig:heeg_diag_dehn_twist}. Thus we want to prove that $N'(\tau_E)\simeq N(\tau_E)$. There are two ways to do it. One is to use directly the definition of $A_\infty$ actions via pseudo-holomorphic curves. In the genus one case such computations were done in~\cite[Section~10.2]{LOT-bim}, and it is possible to generalize them to compute $N(\tau_E)$. However, it is more difficult to do this for $N(\tau_C^{-1})$, which is our next computation. Thus we choose another approach, which we will also use to compute $N(\tau_C^{-1})$.

            Before proceeding to the proof of $N'(\tau_E)\simeq N(\tau_E)$, let us describe the necessary background material. 

            The algebra $\B(\Zpmc)$ is a direct summand of a bigger algebra 
            $$\B(\Zpmc)=\A(\Zpmc, -g +1) \quad \subset  \quad \A(\Zpmc)=\bigoplus\limits_{-g \leq k \leq g}\A(\Zpmc, k),$$
            see~\cite[Section~3]{LOT-main} for the definition of $\A(\Zpmc)$. 
            % The algebra $\A(\Zpmc, g -1)$ is Koszul dual to $\B(\Zpmc)$, and in general we have Koszul duality between algebras $\A(\Zpmc,-k)$ and $\A(\Zpmc,k)$, see~\cite[Section~8]{LOT-mor} for that.

            To a Heegaard diagram of the mapping cylinder $\H(M_\phi)$ one can associate not only a 
            $$\text{type $DA$ bimodule}\quad {}^{\B(\Zpmc)}N(\phi)_{\B(\Zpmc)}={}^{\B(\Zpmc)}\widehat{CFDA}(\phi,-g+1)_{\B(\Zpmc)},$$ 
            but also 
            \begin{align*}
            \text{type $DD$ bimodule}\quad &{}^{\B(\Zpmc)}\widehat{CFDD}(\phi,-g+1)^{\A(\Zpmc, g -1)}, \text{ and} \\
            \text{type $AA$ bimodule}\quad &{}_{\A(\Zpmc, g -1)}\widehat{CFAA}(\phi,-g+1)_{\B(\Zpmc)},
            \end{align*}
            see~\cite[Section~6]{LOT-bim} for the definitions. The sets of generators for these three bimodules are the same, but the $A_\infty$ actions are different. Note that all the three bimodules are direct summands of more general ones 
            \begin{align*}
            {}^{\A(\Zpmc)}\widehat{CFDA}(\phi)_{\A(\Zpmc)}=&\bigoplus_{0\leq j\leq 2g} {}^{\A(\Zpmc, -g +j )}\widehat{CFDA}(\phi,-g+j)_{\A(\Zpmc, -g +j )}, \\ 
            ^{\A(\Zpmc)}\widehat{CFDD}(\phi)^{\A(\Zpmc)}=&\bigoplus_{0\leq j\leq 2g} {}^{\A(\Zpmc, -g +j)}\widehat{CFDD}(\phi,-g+j)^{\A(\Zpmc, g -j)}, \\ 
            _{\A(\Zpmc)}\widehat{CFAA}(\phi)_{\A(\Zpmc)}=&\bigoplus_{0\leq j\leq 2g} {}_{\A(\Zpmc, g -j)}\widehat{CFAA}(\phi,-g+j)_{\A(\Zpmc, -g +j)}.
            \end{align*}
            The $(-g+1)$-summands, which are relevant for us, are characterized by the number of arcs generators occupy on the left and the right side of a Heegaard diagram --- for us it is $2g-1$ on the left, and $1$ on the right. Equivalently, we can say that the $(-g+1)$-summands consist of generators in those $\text{Spin}^c$ structures of $M_\phi$ whose Chern class evaluates to $-2g+2$ on each of the boundaries of $M_\phi$. As we noted in the introduction, the reason why we focus our attention on the $(-g+1)$-summands is because those are the ones corresponding to the fixed point Floer theory of the surface. Because we will be only interested in the $(-g+1)$-summands, we will omit the index $-g+1$ for these bimodules below.

            To complement Theorem~\ref{DA pairing}, we describe below the relationship between the type $AA$, $DD$ and $DA$ bimodules, which follow from \cite[Theorem~12]{LOT-bim}:
            \begin{equation}\label{prrr}
            \begin{split}
            ^{\B(\Zpmc)}\widehat{CFDD}(\phi)^{\A(\Zpmc, g -1)} \boxtimes {}_{\A(\Zpmc, g -1)}\widehat{CFAA}(\psi)_{\B(\Zpmc)} \simeq {}^{\B(\Zpmc)}\widehat{CFDA}(\psi \phi)_{\B(\Zpmc)}, \\
            ^{\B(\Zpmc)}\widehat{CFDA}(\phi)_{\B(\Zpmc)} \boxtimes {}^{\B(\Zpmc)}\widehat{CFDD}(\psi)^{\A(\Zpmc, g -1)} \simeq {}^{\B(\Zpmc)}\widehat{CFDD}(\psi \phi)^{\A(\Zpmc, g -1)}.
            \end{split}
            \end{equation}
            (We use here the notation $\widehat{CFDA}(\phi)$ instead of $N(\phi)$ just to emphasize that $D$ sides are paired with $A$ sides.)

            \emph{Arc-slides} are generators of the mapping class groupoid, see~\cite[Figure 3]{LOT-arc} for a definition of an arc-slide. In particular, they generate Dehn twists. In practice, \cite[Lemma~2.1]{LOT-arc} gives a standard way to decompose a Dehn twist $\tau_l\co(\Sigma,\partial \Sigma) \rightarrow (\Sigma,\partial \Sigma)$ into a product of arc-slides; let us describe it. First, we pick a parameterization of a surface $(\Sigma,\partial \Sigma) \cong F^\circ(\Zpmc)$ such that $l$ is isotopic to an arc $\alpha \subset F^\circ(\Zpmc)$, whose ends are connected along the part of the boundary which does not contain a basepoint --- we denote this part by $I_\alpha$. Then we consider the composition of arc-slides, each of which moves a point on $I_\alpha$, in turn, once along the $\alpha$. This will be the desired Dehn twist. For example, see Figure~\ref{fig:Sigma_2_parameterized}, where the curve $E$ is in the correct position with respect to the arc $\alpha_E$. Thus the Dehn twist $\tau_E$ is equal to the single arc-slide over the arc $\alpha_E$, which is indicated on the picture. 
            There is also a standard way to produce a Heegaard diagram for an arc-slide. In our example of $\tau_E$, this is the 3rd type of the diagram in Figure~\ref{fig:heeg_diag_dehn_twists_1-4}. Based on these standard Heegaard diagrams, the type $DD$ bimodules for all arc-slides were explicitly computed in~\cite{LOT-arc}. 
            % Let us mention that, compared to type $DA$ and type $AA$ bimodules, the moduli spaces of pseudo-holomorphic curves involved in the definition of type $DD$ bimodules are simpler.

            Let us return to the proof of $N'(\tau_E) \simeq N(\tau_E)$. We first prove the following isomorphism:
            \begin{equation}\label{eqqqq}
            {}^{\B(\Zpmc)}N'(\tau_E)_{\B(\Zpmc)} \boxtimes {}^{\B(\Zpmc)}\widehat{CFDD}(\id)^{\A(\Zpmc, g -1)} \cong {}^{\B(\Zpmc)}\widehat{CFDD}(\tau_E)^{\A(\Zpmc, g -1)}.
            \end{equation}
            On the left-hand-side, the bimodule $N'(\tau_E)$ is our candidate bimodule, see Figure~\ref{fig:heeg_diag_dehn_twist}. The bimodule ${}^{\B(\Zpmc)}\widehat{CFDD}(\id)^{\A(\Zpmc, g -1)}$, according~\cite[Theorem~1]{LOT-arc}, is the identity type $DD$ bimodule described in~\cite[Definition~1.3]{LOT-arc}. It consists of four generators 
            $${}_{i_0}(i_0)_{\bar{i}_0}, {}_{i_1}(i_1)_{\bar{i}_1}, {}_{i_2}(i_2)_{\bar{i}_2}, {}_{_3}(i_3)_{\bar{i}_3},$$
            where ${\bar{i}_k}\in \A(\Zpmc, g -1)$ is the idempotent complementary to  ${i}_k\in \B(\Zpmc)$. The type $DD$ actions are 
            \begin{equation}\label{act}
            \partial(i)=\sum_{ \substack {(i \xrightarrow{\rho} j  )\in \B(\Zpmc) \\ i\neq j }} a(\rho)\big|_{\B(\Zpmc)} \o j \o a(\rho)\big|_{\A(\Zpmc, g -1)}, 
            \end{equation}
            where $a(\rho)$ is the notation from ~\cite[Definition~3.23]{LOT-main}, and $a(\rho)\big|_{\A(\Zpmc, g -1)}$ stands for the projection of $a(\rho)$ onto the summand  $\A(\Zpmc, g -1)$ (in general $a(\rho)\in \A(\Zpmc)$).
            Note that, in our case of $\Zpmc=\Zpmc_2$, the property $i\neq j$ implies that in the type $DD$ actions the chord $\rho$ cannot be equal to $\rho_{12}, \rho_{23}, \rho_{56}$ or $\rho_{67}$.

            Next, we describe in detail the right-hand-side of Isomorphism~\eqref{eqqqq}. According to~\cite[Theorem~2]{LOT-arc}, the bimodule ${}^{\B(\Zpmc)}\widehat{CFDD}(\tau_E)^{\A(\Zpmc, g -1)}$ is an arc-slide type $DD$ bimodule described in~\cite[Definition~1.7]{LOT-arc}. Note that our arc-slide $\tau_E$ is an \emph{under-slide} (see~\cite[Definition~4.2]{LOT-arc}), and thus the actions of the corresponding bimodule are described explicitly in~\cite[Definition~4.19]{LOT-arc}. Figure~\ref{fig:arc-slide_DD_bimodule} depicts the resulting bimodule; let us pause to describe the notation. An arrow between generators $\x$ and $\y$ with the label $a \o b$ means that there is a type $DD$ action $\x \rightarrow a \o \y \o b$. For elements of algebra $\A(\Zpmc, g -1)$ we use the strand diagram notation from~\cite[Section~3]{LOT-main}, with the difference that our strands are the horizontal lines numbered by $0$ to $7$, while in~\cite{LOT-main} the strands are numbered by $1$ to $8$. As an example of our notation, $|(0, 2), (1, 3), (4 \rightarrow 5)|$ represents an element of the algebra corresponding to the strand going from the 4th to the 5th line supplemented with idempotents $(0, 2)$ and $(1, 3)$ (this makes this an element of three-strands-moving algebra). In the notation of~\cite[Definition~3.23]{LOT-main}:
            $$|(0, 2), (1, 3), (4 \rightarrow 5)|=a(\rho_5)\big|_{\A(\Zpmc, g -1)}$$ 
            In the notation of~\cite[Definition~3.25]{LOT-main}:
            $$|(0, 2), (1, 3), (4 \rightarrow 5)|=\begin{bmatrix}4 & 0  &1  \\ 5\end{bmatrix}$$   
            Also, every action in arc-slide bimodules has its type, see~\cite[Definition 4.19]{LOT-arc} for the relevant here case of an under-slide. In Figure~\ref{fig:arc-slide_DD_bimodule}, we specify the types of the actions by the superscripts U-n.

            Given the complete description of all the bimodules involved, it is straightforward to check that Isomorphism~\eqref{eqqqq} holds. The actions in Figures~\ref{fig:E_bimodule} and~\ref{fig:arc-slide_DD_bimodule} are intentionally spaced in a similar way, to suggest how the isomorphism works. Taking an action from Figure~\ref{fig:E_bimodule}, and, if possible, box tensoring it with actions from ${}^{\A(\Zpmc, g -1)}\widehat{CFAA}(\id)^{\B(\Zpmc)}$ (described in Equation~\ref{act}), precisely results in the corresponding action from Figure~\ref{fig:arc-slide_DD_bimodule}. 

            Now, box tensoring both sides of Isomorphism~\eqref{eqqqq} by ${}_{\A(\Zpmc, g -1)}\widehat{CFAA}(\id)_{\B(\Zpmc)}$ gives 
            \begin{equation*}
            \begin{split}
            {}^{\B(\Zpmc)}N'(\tau_E)_{\B(\Zpmc)} &\boxtimes {}^{\B(\Zpmc)}\widehat{CFDD}(\id)^{\A(\Zpmc, g -1)} \boxtimes {}_{\A(\Zpmc, g -1)}\widehat{CFAA}(\id)_{\B(\Zpmc)} \simeq \\
            &\cong
            ^{\B(\Zpmc)}\widehat{CFDD}(\tau_E)^{\A(\Zpmc, g -1)} \boxtimes {}_{\A(\Zpmc, g -1)}\widehat{CFAA}(\id)_{\B(\Zpmc)}.
            \end{split}
            \end{equation*}
            In the view of pairings~\eqref{prrr}, the previous isomorphism can be transformed to
            \begin{equation*}
            {}^{\B(\Zpmc)}N'(\tau_E)_{\B(\Zpmc)} \boxtimes {}^{\B(\Zpmc)}\widehat{CFDA}(\id)_{\B(\Zpmc)} \simeq {}^{\B(\Zpmc)}\widehat{CFDA}(\tau_E)_{\B(\Zpmc)}.
            \end{equation*}
            The bimodule ${}^{\B(\Zpmc)}\widehat{CFDA}(\id)_{\B(\Zpmc)}$ is the identity bimodule (\cite[Theorem~4]{LOT-bim}), and thus tensoring with it has no effect. Thus, remembering ${}^{\B(\Zpmc)}N(\tau_E)_{\B(\Zpmc)}={}^{\B(\Zpmc)}\widehat{CFDA}(\tau_E)_{\B(\Zpmc)}$, the above homotopy equivalence transforms to
            \begin{equation*}
            {}^{\B(\Zpmc)}N'(\tau_E)_{\B(\Zpmc)} \simeq  {}^{\B(\Zpmc)}N(\tau_E)_{\B(\Zpmc)},
            \end{equation*}
            which finishes the proof.
        \end{computation}
        \begin{computation}[$N(\tau_C^{-1})$]\label{harder one}
            Let us denote for a moment our candidate bimodule from Figure~\ref{fig:C_inv_bimodule} by $N'(\tau_C^{-1})$, and the bimodule $N(\tau_C^{-1})$ will be the one which corresponds to the Heegaard diagram in Figure~\ref{fig:heeg_diag_dehn_twist_connecting_curve}
            . Thus we want to prove that $N'(\tau_C^{-1}) \simeq N(\tau_C^{-1})$. 

            First, we factorize the Dehn twist $\tau_C^{-1}$ into the product of arc-slides. In Figure~\ref{fig:Sigma_2_parameterized}, consider a slide of the arc $\alpha_E$ over the arc $\alpha_A$, and let us call this arc-slide $\eta$. Then we get a new parameterization of the surface, where instead of the arc $\alpha_E$ we have a new arc $ \alpha'_E$ which is isotopic to the curve $C' \sim \eta(C)$, if 
            % $C' \sim \eta(C) \sim C$
            we connect the ends of the arc $ \alpha'_E$. According to~\cite[Lemma~2.1]{LOT-arc}, the Dehn twist $\tau_{C'}^{-1}$ can be factorized into four arc-slides $ \mu_1,\mu_2,\mu_3,\mu_4$ along the arc $ \alpha'_E$, which we picture in Figure~\ref{fig:arcslides_sequence}. Now, by a mapping class group relation $f\tau_l f^{-1}=\tau_{f(l)}$ we obtain the desired factorization:
            $$
            \tau_C^{-1}=\tau_{\eta^{-1}(C')}^{-1}=\eta^{-1} \tau_{C'}^{-1} \eta =\eta^{-1} \mu_4\mu_3\mu_2\mu_1 \eta.
            $$
            \begin{figure}[H]
            \centering
            \includegraphics[width=0.6\textwidth]{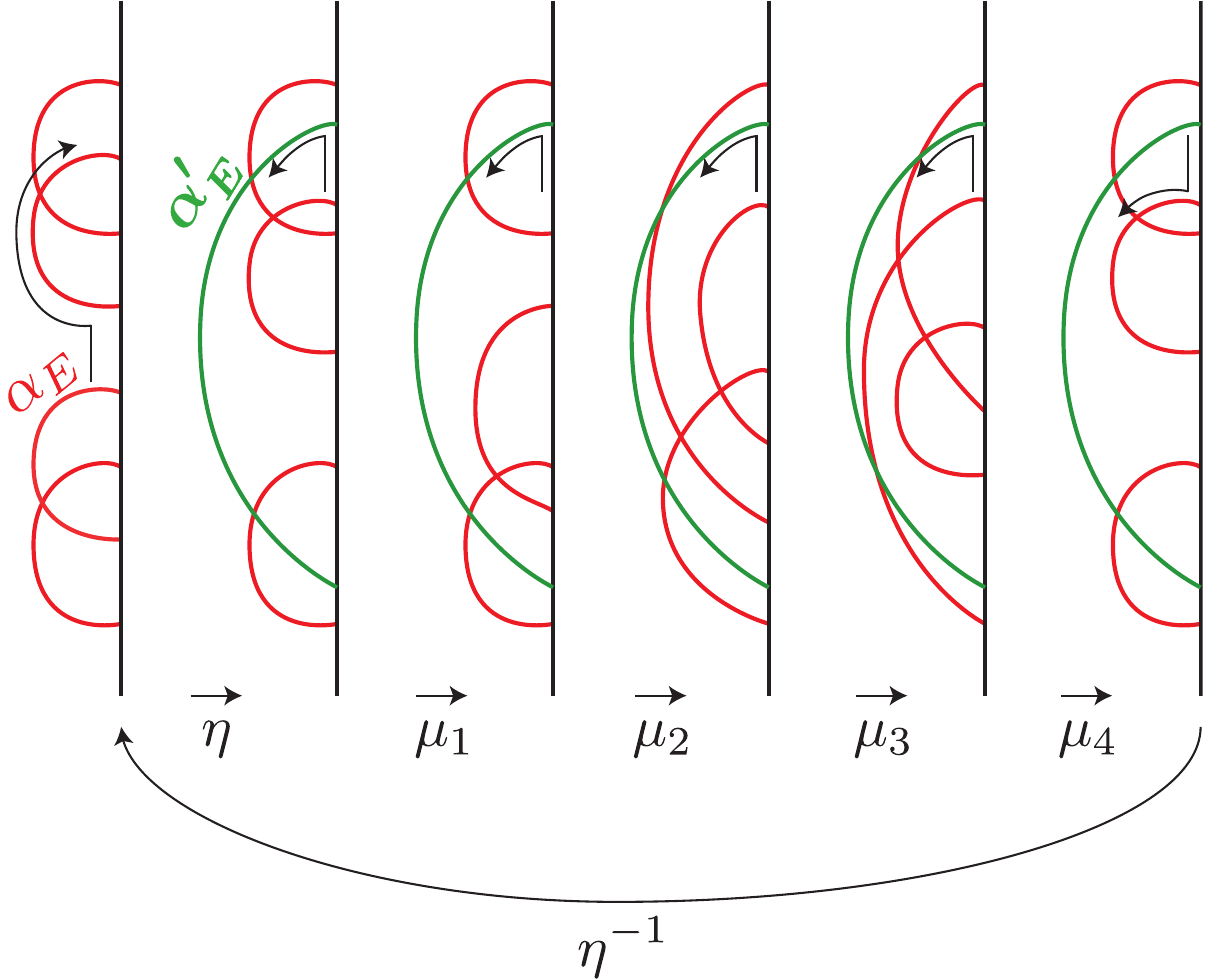}
            \caption{Composition of the above six arc-slides gives a left handed Dehn twist along the curve C in Figure~\ref{fig:Sigma_2_parameterized}, i.e. $\tau_C^{-1}=\eta^{-1} \mu_4\mu_3\mu_2\mu_1 \eta$.}
            \label{fig:arcslides_sequence}
            \end{figure}
            From this factorization we get
            $$
            N(\eta)\boxtimes N(\mu_1)  \boxtimes N(\mu_2) \boxtimes N(\mu_3) \boxtimes N(\mu_4) \boxtimes N(\eta^{-1}) \simeq N(\tau_{C'}^{-1}),
            $$
            and so to compute $N(\tau_{C'}^{-1})$ it is left to compute the bimodules for the six arc-slides. These are computed using exactly the same method we used in the previous computation of $N(\tau_{E})$; see \cite[Showcase~10]{Pyt} for explicit descriptions of all 6 bimodules. (The strand notation for algebra elements in the program is the same as the one used in Figure~\ref{fig:arc-slide_DD_bimodule}.)
            % (with the exception which is described in the remark at the end of Computation~\ref{comp: dehn twist, the easier one}).

            For computing $N(\eta)\boxtimes N(\mu_1)  \boxtimes N(\mu_2) \boxtimes N(\mu_3) \boxtimes N(\mu_4) \boxtimes N(\eta^{-1})$ we used the program~\cite[Showcase~10]{Pyt}: tensoring all the six arc-slide bimodules, and then doing all possible cancellations, results in a bimodule isomorphic to the one in Figure~\ref{fig:C_inv_bimodule} with canceled differential $x_2 \rightarrow t_{12}$. This proves the desired homotopy equivalence:
            $$
            N'(\tau_{C'}^{-1}) \simeq  N(\eta)\boxtimes N(\mu_1)  \boxtimes N(\mu_2) \boxtimes N(\mu_3) \boxtimes N(\mu_4) \boxtimes N(\eta^{-1}) \simeq N(\tau_{C'}^{-1}).
            $$
        \end{computation}

        \newpage
        \begin{figure}[ht]
        \begin{subfigure}{0.5\textwidth}
        \centering
        \includegraphics[width=1\textwidth]{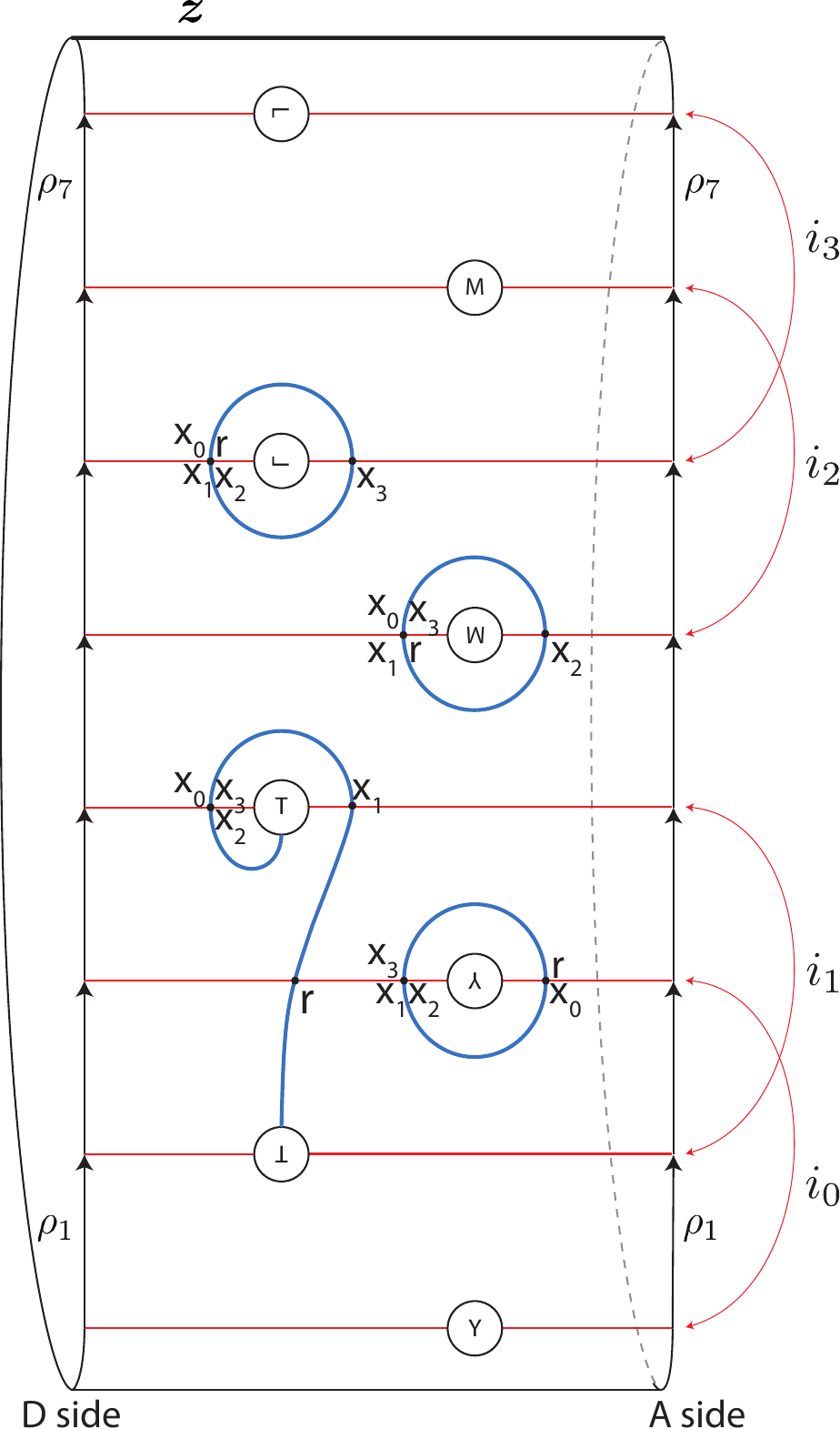}
        \caption{Heegaard diagram $\H (M_{\tau_E})$ for the right handed Dehn twist along the curve E.}
        \label{fig:heeg_diag_dehn_twist}
        \end{subfigure}
        \begin{subfigure}{0.446\textwidth}
        \centering
        \includegraphics[width=1\textwidth]{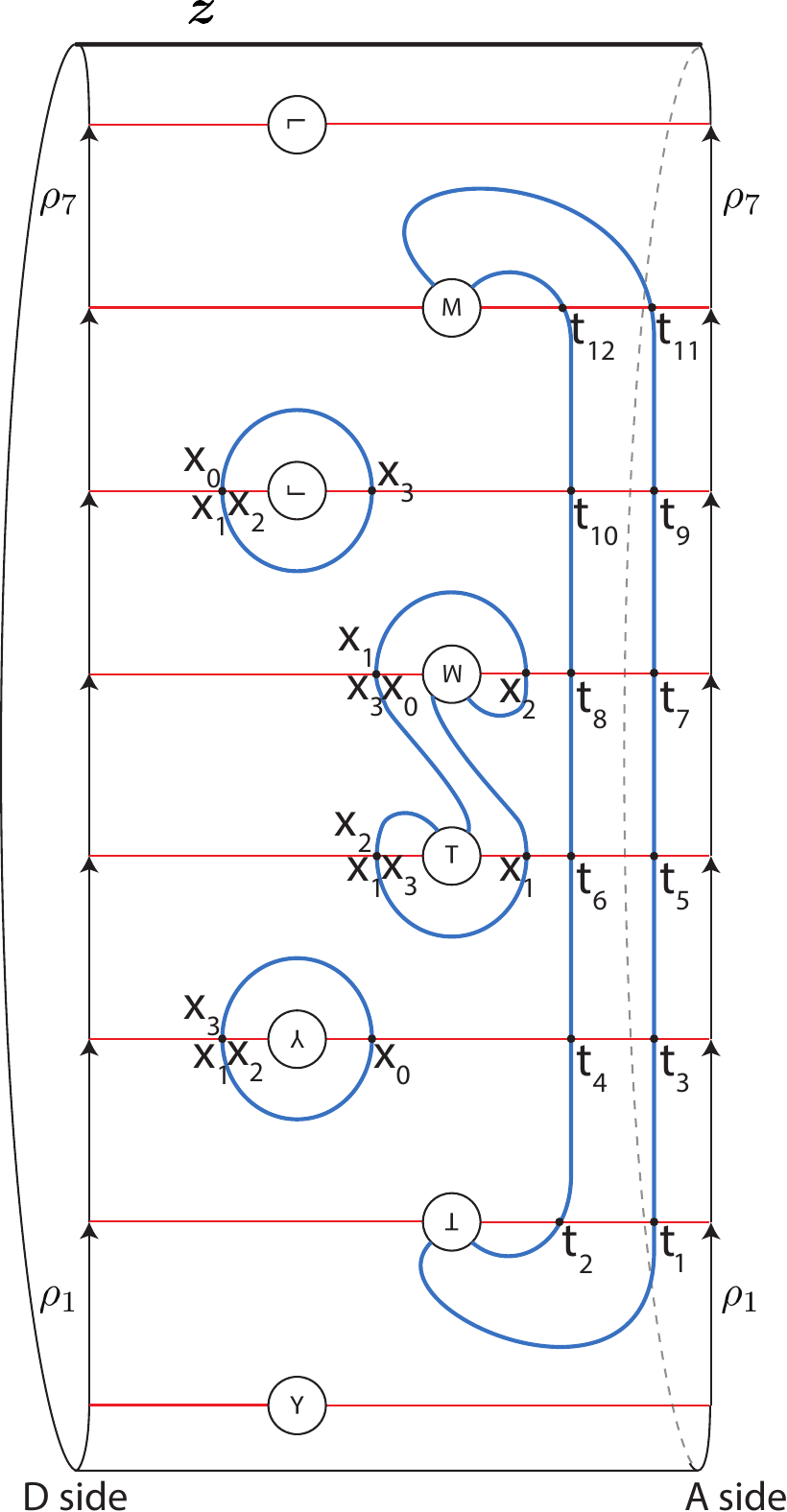}
        \caption{Heegaard diagram $\H (M_{\tau_C^{-1}})$ for the left handed Dehn twist along the curve C.}
        \label{fig:heeg_diag_dehn_twist_connecting_curve}
        \end{subfigure}
        \caption{}
        \label{fig:two_Dehn}
        \end{figure}
        \newpage
        \begin{figure}[H]
            \tikzstyle{arrow} = [thick,->,>=stealth]
            \begin{tikzpicture}[scale=1]
                \node[](x_2) at (-4,7){$_{i_2}{(x_2)}_{i_2}$};
                \node[](x_3) at (6,7){$_{i_3}{(x_3)}_{i_3}$};
                \node[](x_0) at (-4,0){$_{i_0}{(x_0)}_{i_0}$};
                \node[](x_1) at (6,0){$_{i_1}{(x_1)}_{i_1}$};
                \node[](r) at (1,-7){$_{i_1}{(r)}_{i_0}$};

                \draw[arrow] (x_2) to[loop above,looseness=15] node[above]{$\rho_{56} \o \rho_{56}$} (x_2);
                \draw[arrow] (x_3) to[loop above,looseness=15] node[above]{$\rho_{67} \o \rho_{67}$} (x_3);
                \draw[arrow] (x_1) to[loop,out=120,in=150,looseness=15] node[above,sloped]{$\rho_{23} \o \rho_{23}$} (x_1);
                \draw[arrow] (x_3) to[bend right=10] node[above]{$\rho_6 \o \rho_6$} (x_2);
                \draw[arrow] (x_2) to[bend left] node[above,text width=2.5cm,text centered]{$\rho_5 \o \rho_5$ \\  $+ \rho_7 \o \rho_7$ \\  $+ \rho_{567} \o \rho_{567}$} (x_3);

                \draw[arrow] (x_0) to[bend right] node[above,text width=3.3cm,text centered]{$\rho_1 \o \rho_1$ \\ $+\rho_3 \o (\rho_3,\rho_{23})$ \\  $+ \rho_{123} \o \rho_{123}$} (x_1);
                \draw[arrow] (r) to[bend left, in=180] node[above,sloped]{$1\o \rho_{3}$} (x_1);
                \draw[arrow] (x_1) to[bend left] node[above,sloped]{$\rho_{23}\o \rho_{2}$} (r);
                \draw[arrow] (r) to[bend left] node[above,sloped]{$ \rho_{2} \o 1$} (x_0);
                \draw[arrow] (x_0) to[bend left, out=0] node[below,sloped]{$\rho_{123}\o \rho_{12} + \rho_{3}\o (\rho_{3}, \rho_2)$} (r);

                \draw[arrow] (x_0) to node[above,sloped]{$\rho_{1234}\o \rho_{1234} + \rho_{3456}\o (\rho_{3},\rho_{23456})$} node[below,sloped]{$\rho_{34}\o (\rho_{3},\rho_{234}) + \rho_{123456}\o \rho_{123456}$} (x_2);

                \draw[arrow] (x_0) to[bend left, out=15] node[above,sloped,text width=4cm,text centered]{$\rho_{12345}\o \rho_{12345}$\\  $+ \rho_{345}\o (\rho_{3},\rho_{2345})$\\  $\ + \rho_{34567}\o (\rho_{3},\rho_{234567})$\\  $+ \rho_{1234567}\o \rho_{1234567}$} (x_3);

                \draw[arrow] (x_1) to[bend left=47, in=165] node[above,sloped,text width=4cm,text centered]{$\rho_{234}\o \rho_{234}$\\  $+ \rho_{4}\o\rho_{4}+ \rho_{456}\o\rho_{456}$\\  $+ \rho_{23456}\o \rho_{23456}$} (x_2);

                \draw[arrow] (x_1) to node[above,sloped]{$\rho_{4567}\o \rho_{4567} + \rho_{2345}\o \rho_{2345}$}  node[below,sloped]{$\rho_{45}\o \rho_{45}  + \rho_{234567}\o \rho_{234567}$} (x_3);

                \draw[arrow] (r) to[bend left=53] node[below,sloped]{$\rho_{4}\o \rho_{34} + \rho_{456}\o \rho_{3456}$} (x_2);

                \draw[arrow] (r) to[bend right=50] node[below,sloped,text width=4cm,text centered]{$\rho_{45}\o \rho_{345}$\\ $+ \rho_{4567}\o \rho_{34567}$} (x_3);
            \end{tikzpicture}
            \caption{Bimodule ${}^{\B(\Zpmc_2)}N(\tau_E)_{\B(\Zpmc_2)}$.}
            \label{fig:E_bimodule}
        \end{figure}
        \newpage
        \begin{figure}[H]
            \tikzstyle{arrow} = [thick,->,>=stealth]
            \begin{tikzpicture}[scale=0.95]
                \node(x_2) at (2,8.5){$_{i_2}{(x_2)}_{i_2}$};
                \node(x_3) at (10,8.5){$_{i_3}{(x_3)}_{i_3}$};
                \node(x_0) at (2,5){$_{i_0}{(x_0)}_{i_0}$};
                \node(x_1) at (10,5){$_{i_1}{(x_1)}_{i_1}$};
                \node(t_12) at (6,1){$_{i_2}{(t_{12})}_{i_2}$};
                \node(t_11) at (10,1){$_{i_1}{(t_{11})}_{i_2}$};
                \node(t_10) at (6,-1){$_{i_2}{(t_{10})}_{i_3}$};
                \node(t_9) at (10,-1){$_{i_1}{(t_9)}_{i_3}$};
                \node(t_8) at (6,-3){$_{i_2}{(t_8)}_{i_2}$};
                \node(t_7) at (10,-3){$_{i_1}{(t_7)}_{i_2}$};
                \node(t_6) at (6,-5){$_{i_2}{(t_6)}_{i_1}$};
                \node(t_5) at (10,-5){$_{i_1}{(t_5)}_{i_1}$};
                \node(t_4) at (6,-7){$_{i_3}{(t_4)}_{i_0}$};
                \node(t_3) at (10,-7){$_{i_1}{(t_3)}_{i_0}$};
                \node(t_2) at (6,-9){$_{i_2}{(t_2)}_{i_1}$};
                \node(t_1) at (10,-9){$_{i_1}{(t_1)}_{i_1}$};

                \draw[arrow] (x_2) to[loop above,looseness=15] node[above]{$\rho_{56} \o \rho_{56}$} (x_2);
                \draw[arrow] (x_1) to[loop,out=120,in=150,looseness=15] node[above,sloped]{$\rho_{23} \o \rho_{23}$} (x_1);
                \draw[arrow] (x_3) to[bend right=30] node[above]{$\rho_6 \o \rho_6$} (x_2);
                \draw[arrow] (x_2) to[bend left=10] node[above]{$\rho_5 \o \rho_5$} (x_3);
                \draw[arrow] (x_0) to[bend right=18] node[below right,sloped]{$\rho_3 \o \rho_3$} (x_1);
                \draw[arrow] (x_1) to[bend left=10] node[above]{$\rho_2 \o \rho_2$} (x_0);
                \draw[arrow] (x_1) to node[above,sloped]{$\rho_{2345}\o \rho_{2345}$}  node[below,sloped]{$\rho_{45}\o \rho_{45}$} (x_3);
                \draw[arrow] (x_0) to node[above,sloped]{$\rho_{3456}\o \rho_{3456}$} node[below,sloped]{$\rho_{34}\o \rho_{34}$} (x_2);
                \draw[arrow] (x_0) to[bend left] node[above,sloped]{$\rho_{345}\o \rho_{345}$}  (x_3);
                \draw[arrow] (x_1) to[bend left,in=163] node[above,sloped,text width=4cm,text centered,pos=0.35]{$\rho_{234}\o \rho_{234}$\\  $+ \rho_{4}\o\rho_{4}+ \rho_{456}\o\rho_{456}$\\  $+ \rho_{23456}\o \rho_{23456}$} (x_2);

                \draw[arrow] (t_1) to node[below]{$\rho_4 \o 1$} (t_2);
                \draw[arrow] (t_3) to node[above]{$\rho_4 \o 1$} (t_4);
                \draw[arrow] (t_5) to node[above]{$\rho_4 \o 1$} (t_6);
                \draw[arrow] (t_7) to node[above]{$\rho_4 \o 1$} (t_8);
                \draw[arrow] (t_9) to node[above]{$\rho_4 \o 1$} (t_10);
                \draw[arrow] (t_11) to node[below]{$\rho_4 \o 1$} (t_12);
                \draw[arrow] (t_1) to node[left]{$1 \o \rho_2$} (t_3);
                \draw[arrow] (t_3) to node[left]{$1 \o \rho_3$} (t_5);
                \draw[arrow] (t_5) to node[left]{$1 \o \rho_4$} (t_7);
                \draw[arrow] (t_7) to node[left]{$1 \o \rho_5$} (t_9);
                \draw[arrow] (t_9) to node[left]{$1 \o \rho_6$} (t_11);
                \draw[arrow] (t_2) to node[right]{$1 \o \rho_2$} (t_4);
                \draw[arrow] (t_4) to node[right]{$1 \o \rho_3$} (t_6);
                \draw[arrow] (t_6) to node[right]{$1 \o \rho_4$} (t_8);
                \draw[arrow] (t_8) to node[right]{$1 \o \rho_5$} (t_10);
                \draw[arrow] (t_10) to node[right]{$1 \o \rho_6$} (t_12);
                \draw[arrow] (t_1) to[bend right] node[right]{$1 \o \rho_{23}$} (t_5);
                \draw[arrow,opacity=0.2] (t_3) to[bend right] (t_7);
                \draw[arrow,opacity=0.2] (t_5) to[bend right] (t_9);
                \draw[arrow,opacity=0.2] (t_7) to[bend right] (t_11);
                \draw[arrow] (t_2) to[bend left] node[left]{$1 \o \rho_{23}$} (t_6);
                \draw[arrow,opacity=0.2] (t_4) to[bend left] (t_8);
                \draw[arrow,opacity=0.2] (t_6) to[bend left] (t_10);
                \draw[arrow,opacity=0.2] (t_8) to[bend left] (t_12);
                \draw[arrow,opacity=0.2] (t_1) to[bend right=40] (t_7);
                \draw[arrow,opacity=0.2] (t_3) to[bend right=40] (t_9);
                \draw[arrow,opacity=0.2] (t_5) to[bend right=40] (t_11);
                \draw[arrow,opacity=0.2] (t_2) to[bend left=40] (t_8);
                \draw[arrow,opacity=0.2] (t_4) to[bend left=40] (t_10);
                \draw[arrow,opacity=0.2] (t_6) to[bend left=40] (t_12);
                \draw[arrow,opacity=0.2] (t_1) to[bend right=50] (t_9);
                \draw[arrow,opacity=0.2] (t_3) to[bend right=50] (t_11);
                \draw[arrow,opacity=0.2] (t_2) to[bend left=50] (t_10);
                \draw[arrow,opacity=0.2] (t_4) to[bend left=50] (t_12);
                \draw[arrow,opacity=0.2] (t_1) to[bend right=60] (t_11);
                \draw[arrow,opacity=0.2] (t_2) to[bend left=60] (t_12);

                \draw[arrow] (x_1) to[] node[below,sloped]{$1 \o \rho_{4} + \rho_{23} \o \rho_{234} $} (t_11);
                \draw[arrow] (x_1) to[bend left=60] node[above,sloped,pos=0.9]{$\rho_{23456} \o 1$} (t_2);
                \draw[arrow] (x_1) to[bend left=60] node[below,sloped]{$\rho_{456} \o 1$} (t_6);
                
                \draw[arrow] (x_0) to[] node[below,sloped]{$\rho_{3} \o \rho_{34} + \rho_1 \o \rho_{123456} $} (t_11);
                \draw[arrow] (x_0) to[bend right] node[below,sloped]{$\rho_{3456} \o 1$} (t_4);

                \draw[arrow] (x_2) to[bend right=60] node[below,sloped]{$\rho_{56} \o 1$} (t_8);
                \draw[arrow] (x_2) to[bend right=60] node[below,sloped]{$1$} (t_12);

                \draw[arrow] (x_3) to[bend left=45] node[below,sloped,pos=0.8]{$\rho_6 \o 1$} (t_10);

                \draw[arrow] (x_0) to[bend right] node[below,sloped]{$\rho_1 \o \rho_1$} (t_1);
                \draw[arrow] (x_0) to[bend right] node[below,sloped]{$\rho_1 \o \rho_{12}$} (t_3);
                \draw[arrow] (x_0) to[bend right] node[below,sloped]{$\rho_1 \o \rho_{123}$} (t_5);
                \draw[arrow] (x_0) to[bend right] node[below,sloped]{$\rho_1 \o \rho_{1234}$} (t_7);
                \draw[arrow] (x_0) to[bend right] node[below,sloped]{$\rho_1 \o \rho_{12345}$} (t_9);
                % \draw[arrow] (x_0) to[] node[below,sloped]{$\rho_1 \o \rho_{123456}$} (t_11);

                \draw[arrow] (t_11) to[bend right] node[below,sloped,pos=0.3]{$\rho_{234567} \o \rho_7$} (x_3);
                \draw[arrow] (t_9) to[bend right=50] node[above,sloped,pos=0.2]{$\rho_{234567} \o \rho_{67}$} (x_3);
                \draw[arrow] (t_7) to[bend right=50] node[below,sloped,pos=0.1]{$\rho_{234567} \o \rho_{567}$} (x_3);
                \draw[arrow] (t_5) to[bend right=50] node[below,sloped,pos=0.1]{$\rho_{234567} \o \rho_{4567}$} (x_3);
                \draw[arrow] (t_3) to[bend right=50] node[below,sloped,pos=0.1]{$\rho_{234567} \o \rho_{34567}$} (x_3);
                \draw[arrow] (t_1) to[bend right=50] node[below,sloped,pos=0.1]{$\rho_{234567} \o \rho_{234567}$} (x_3);

                \draw[arrow] (x_0) to[] node[above,sloped, pos=0.15]{$\rho_{1234567} \o \rho_{1234567}$} (x_3);
                \draw[arrow] (t_2) to[] node[above,sloped,pos=0.7]{$\rho_7 \o (\rho_4,\rho_7)$} (x_3);
            \end{tikzpicture}
            \caption{Bimodule ${}^{\B(\Zpmc_2)}N(\tau_C^{-1})_{\B(\Zpmc_2)}$.}
            \label{fig:C_inv_bimodule}
        \end{figure}

        \newpage
        
        \begin{figure}[H]
            \tikzstyle{arrow} = [thick,->,>=stealth]
            \begin{tikzpicture}[scale=1]
                \node[](x_2) at (-4,7){$_{i2}{(x2)}_{|(0, 2),(1, 3),(5, 7)|}$};
                \node[](x_3) at (6,7){$_{i3}{(x3)}_{|(0, 2),(1, 3),(4, 6)|}$};
                \node[](x_0) at (-4,0){$_{i0}{(x0)}_{|(1, 3),(4, 6),(5, 7)|}$};
                \node[](x_1) at (6,0){$_{i1}{(x1)}_{|(0, 2),(4, 6),(5, 7)|}$};
                \node[](r) at (1,-7){$_{i1}{(r)}_{|(1, 3),(4, 6),(5, 7)|}$};
                {\scriptsize
                \draw[arrow] (x_2) to[bend left] node[above,text width=4.5cm,text centered]{
                $\rho_{567} \o |(0, 2),(1, 3),(4 \! \rightarrow \!  7)|^\text{U-1}$  \\  
                $+ \rho_5 \o |(0, 2),(1, 3),(4 \! \rightarrow \!  5)| ^\text{U-1}$ \\  
                $+ \rho_{7} \o |(0, 2),(1, 3),(6 \! \rightarrow \!  7)|^\text{U-1}$ } (x_3);

                \draw[arrow] (x_3) to[bend right=10] node[above]{$(\rho_6 \o |(0, 2),(1, 3),(5 \! \rightarrow \!  6)| ^\text{U-1}$} (x_2);

                \draw[arrow] (x_0) to[bend right=10] node[below,text width=5cm,text centered]{
                $\rho_1 \o |(4, 6),(5, 7),(0 \! \rightarrow \!  1)| ^\text{U-1}$ \\ 
                $+ \rho_{123} \o |(4, 6),(5, 7),(0 \! \rightarrow \!  3)|^\text{U-1}$ } (x_1);

                \draw[arrow] (x_0) to node[above,sloped,text width=5cm,text centered]{
                $\rho_{1234}\o |(1, 3),(5, 7),(0 \! \rightarrow \!  4)|^\text{U-1}$ 
                \\
                $\rho_{34}\o |(5, 7),(1 \! \rightarrow \!  4),(2 \! \rightarrow \!  3)|^\text{U-6}$
                } 
                node[below,sloped,text width=5cm,text centered]{
                $\rho_{123456}\o|(1, 3),(5, 7),(0 \! \rightarrow \!  6)|^\text{U-1}$ 
                \\
                $\rho_{3456}\o |(5, 7),(1 \! \rightarrow \!  6),(2 \! \rightarrow \!  3)|^\text{U-6} $
                } (x_2);

                \draw[arrow] (x_0) to[bend left, out=15] node[above,sloped,text width=6cm,text centered]{
                 $\rho_{12345}\o |(1, 3),(4, 6),(0 \! \rightarrow \!  5)|^\text{U-1}$\\  
                 $+ \rho_{345}\o |(4, 6),(1 \! \rightarrow \!  5),(2 \! \rightarrow \!  3)|^\text{U-6}$\\ 
                  $\ + \rho_{34567}\o |(4, 6),(1 \! \rightarrow \!  7),(2 \! \rightarrow \!  3)|^\text{U-6}$\\ 
                   $+ \rho_{1234567}\o |(1, 3),(4, 6),(0 \! \rightarrow \!  7)|^\text{U-1}$} (x_3);

                \draw[arrow] (x_1) to[bend left=13, in=165] node[above,sloped,text width=6cm,text centered, near start]{
                $\rho_{234}\o |(0, 2),(5, 7),(1 \! \rightarrow \!  4)|^\text{U-1}$\\  
                $+ \rho_{4}\o|(0, 2),(5, 7),(3 \! \rightarrow \!  4)|^\text{U-1} $\\
                $+ \rho_{456}\o|(0, 2),(5, 7),(3 \! \rightarrow \!  6)|^\text{U-1}$\\  
                $+ \rho_{23456}\o |(0, 2),(5, 7),(1 \! \rightarrow \!  6)|^\text{U-1} $ } (x_2);
                
                \draw[arrow] (x_1) to node[above,sloped,text width=6cm,text centered]{
                $\rho_{4567}\o |(0, 2),(4, 6),(3 \! \rightarrow \!  7)|^\text{U-1}$ \\ 
                $+\rho_{2345}\o |(0, 2),(4, 6),(1 \! \rightarrow \!  5)|^\text{U-1}$}  
                node[below,sloped,text width=6cm,text centered]{
                $\rho_{45}\o |(0, 2),(4, 6),(3 \! \rightarrow \!  5)|^\text{U-1} $ \\
                $+ \rho_{234567}\o |(0, 2),(4, 6),(1 \! \rightarrow \!  7)|^\text{U-1}$} (x_3);
                
                \draw[arrow] (r) to[bend left] node[above,sloped]{$ \rho_{2} \o 1^\text{U-2}$} (x_0);
                \draw[arrow] (r) to[bend left, in=180] node[above,sloped]{$1\o |(4, 6),(5, 7),(2 \! \rightarrow \!  3)|^\text{U-2}$} (x_1);
                
                \draw[arrow] (r) to[bend left=55] node[below,sloped,text width=6cm,text centered]{$\rho_{4}\o |(1, 3),(5, 7),(2 \! \rightarrow \!  4)|^\text{U-3}$ 
                \\ $\rho_{456}\o |(1, 3),(5, 7),(2 \! \rightarrow \!  6)|^\text{U-3}$} (x_2);

                \draw[arrow] (r) to[bend right=54] node[near start, below,sloped,text width=6cm,text centered]{$\rho_{45}\o |(1, 3),(4, 6),(2 \! \rightarrow \!  5)|^\text{U-3}$\\ $+ \rho_{4567}\o |(1, 3),(4, 6),(2 \! \rightarrow \!  7)|^\text{U-3}$} (x_3);

                \draw[arrow] (x_1) to[bend left=10] node[above,sloped]{$\rho_{23}\o |(4, 6),(5, 7),(1 \! \rightarrow \!  2)|^\text{U-4}$} (r);
                \draw[arrow] (x_0) to[bend left=10] node[below,sloped]{$\rho_{3}\o |(4, 6),(5, 7),(1 \! \rightarrow \!  3)|^\text{U-4}$} (r);
                }
            \end{tikzpicture}
            \caption{Bimodule ${}^{\B(\Zpmc)}\widehat{CFDD}(\tau_E)^{\A(\Zpmc, g -1)}$.}
            \label{fig:arc-slide_DD_bimodule}
        \end{figure}

    \newpage
    \subsection{Method to compute Hochschild homology}\label{alg_hoh}
        \begin{wrapfigure}{r}{0.4\textwidth}
        \includegraphics[width=0.4\textwidth]{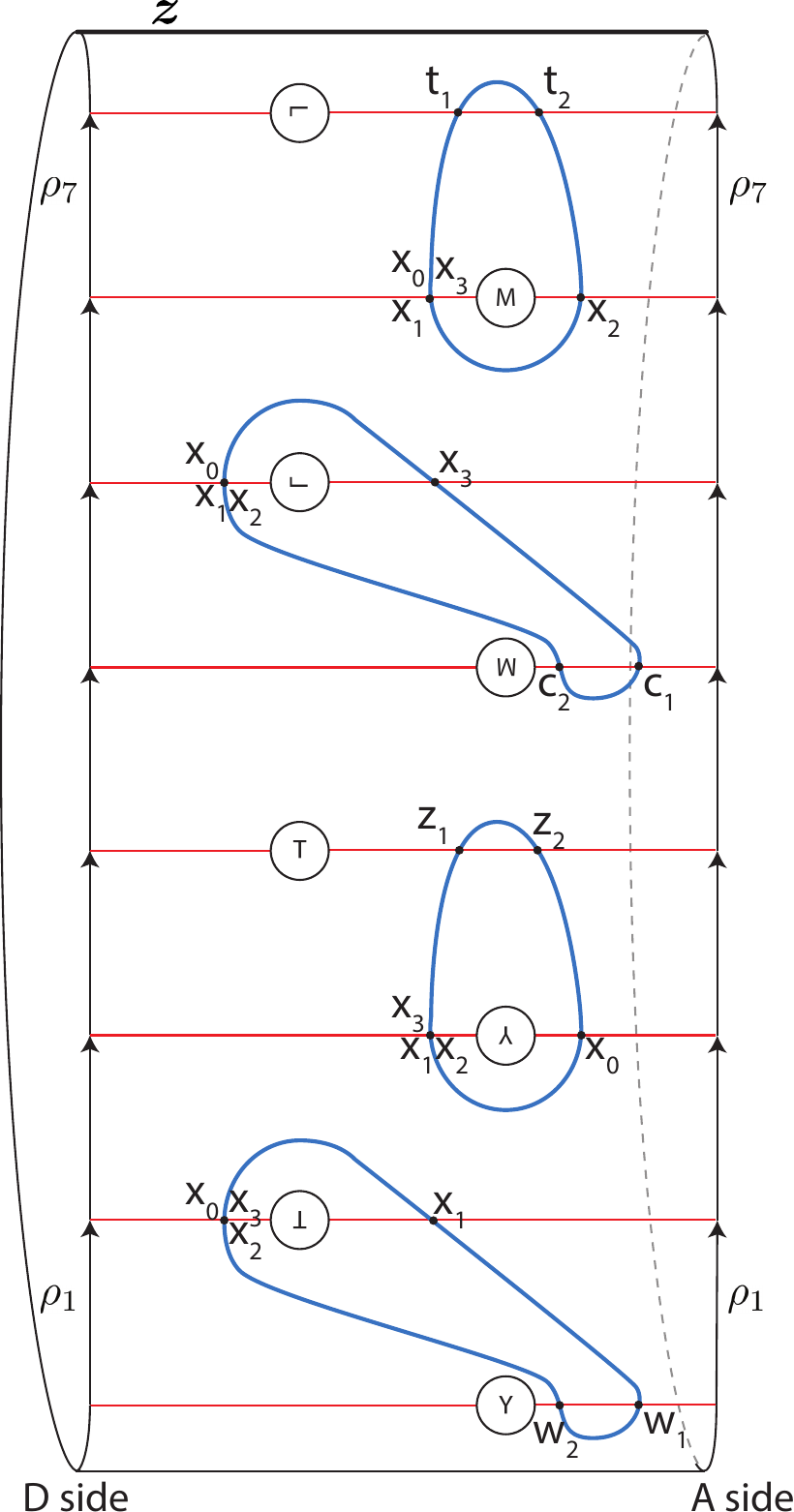}
        \caption{Heegaard diagram for $M_{\id}$ with no periodic domains, where $\id\co F^\circ(\Zpmc_2) \rightarrow F^\circ(\Zpmc_2)$ is the identity mapping class.}
        \label{fig:heegard_diagram_id_bounded}
        \end{wrapfigure}

        It is the Hochschild homology of a bimodule that we are going to equate with a version of fixed point Floer cohomology. Thus we would like to be able to compute it. The method from~\cite[Section~2.3.5]{LOT-bim} for computing Hochschild homology for type $DA$ bimodules works well, as long as the bimodule is bounded (see~\cite[Definition 2.2.46]{LOT-bim}). We will be computing Hochschild homology for many bimodules that are \emph{not bounded}, including all bimodules~\eqref{ten_bimodules}, and so the method does not apply off the shelf. To address this we multiply a given bimodule by a certain bounded bimodule on the left and on the right such that the $A_\infty$ homotopy equivalence class does not change: 

        \vspace{0.2cm}
        {\color{white} this is a} $\displaystyle  [\mathbb{I}]^b \boxtimes {} N(\phi)\boxtimes [\mathbb{I}]^b \simeq  N(\phi).$

        \vspace{0.2cm}
        \noindent We now describe the construction of ${}^{\B(\Zpmc)}[\mathbb{I}]^b_{ \B(\Zpmc)}$. 

        Tensoring with the identity bimodule ${}^{\B(\Zpmc)}N(\id)_{\B(\Zpmc)} \simeq {}^{\B(\Zpmc)}[\mathbb{I}]_{\B(\Zpmc)}$ (see~\cite[Definition~2.2.48]{LOT-bim}) does not change the $A_\infty$ homotopy equivalence class. Thus to obtain the bimodule $[\mathbb{I}]^b$ it is enough to find a way to change ${}^{\B(\Zpmc)}[\mathbb{I}]_{\B(\Zpmc)}$ such that it becomes bounded, but does not change its $A_\infty$ homotopy equivalence class.

        In the genus two case, we claim that the needed bimodule ${}^{\B(\Zpmc_2)}[\mathbb{I}]^b_{ \B(\Zpmc_2)}$ is depicted in Figure~\ref{fig:id_bounded_bimodule}. The graph on that figure does not have any cycles, thus the bimodule is bounded. Program~\cite[Showcase~11]{Pyt} finds that canceling four differentials $c_1 \rightarrow c_2$, $t_1 \rightarrow t_2$, $z_1 \rightarrow z_2$, and $w_1 \rightarrow w_2$ in ${}^{\B(\Zpmc_2)}[\mathbb{I}]^b_{ \B(\Zpmc_2)}$ gives ${}^{\B(\Zpmc_2)}[\mathbb{I}]_{\B(\Zpmc_2)}$, hence ${}^{\B(\Zpmc_2)}[\mathbb{I}]^b_{\B(\Zpmc_2)} \simeq {}^{\B(\Zpmc_2)}[\mathbb{I}]_{\B(\Zpmc_2)}$. Using this bounded identity bimodule, as well as ~\cite[Proposition~2.3.54]{LOT-bim}, in~\cite{Pyt} we implemented an algorithm to find Hochschild homology of any type $DA$ bimodule over $\B(\Zpmc_2)$. Resulting computations of Hochschild homology became the basis of Conjecture~\ref{conj:in_intro}.

        The bounded model was guessed based on the diagram in Figure~\ref{fig:heegard_diagram_id_bounded} (three intersections on the left side of the diagram are omitted for generators $z_1,z_2,c_1,c_2,t_1,t_2,w_1,w_2$), which is the Heegaard diagram from Figure~\ref{fig:heegard_diagram_id} but with perturbed blue curves so that there are no periodic domains.

        \begin{figure}[H]
            \tikzstyle{arrow} = [thick,->,>=stealth]
            \begin{tikzpicture}[scale=0.9]
                \node(w_2) at (-6,0){$_{i_1}{(w_2)}_{i_0}$};
                \node(w_1) at (-3,0){$_{i_1}{(w_1)}_{i_0}$};
                \node(x_0) at (0,0){$_{i_0}{(x_0)}_{i_0}$};
                \node(x_1) at (3,0){$_{i_1}{(x_1)}_{i_1}$};
                \node(z_1) at (6,0){$_{i_0}{(z_1)}_{i_1}$};
                \node(z_2) at (9,0){$_{i_0}{(z_2)}_{i_1}$};

                \node(c_2) at (-6,16){$_{i_3}{(c_2)}_{i_2}$};
                \node(c_1) at (-3,16){$_{i_3}{(c_1)}_{i_2}$};
                \node(x_2) at (0,16){$_{i_2}{(x_2)}_{i_2}$};
                \node(x_3) at (3,16){$_{i_3}{(x_3)}_{i_3}$};
                \node(t_1) at (6,16){$_{i_2}{(t_1)}_{i_3}$};
                \node(t_2) at (9,16){$_{i_2}{(t_2)}_{i_3}$};

                \draw[arrow] (c_1) to node[above]{$1$} (c_2);
                \draw[arrow] (c_1) to node[above]{$\rho_{6} \o \rho_{56}$} (x_2);
                \draw[arrow] (c_1) to[bend left] node[above]{$1 \o \rho_{5}$} (x_3);
                \draw[arrow] (c_1) to[bend left=40] node[above]{$\rho_6 \o \rho_{567}$} (t_2);
                \draw[arrow] (x_2) to[bend right] node[above]{$\rho_5 \o 1$} (c_2);
                \draw[arrow] (x_2) to[bend left=36] node[above]{$1 \o \rho_7$} (t_2);
                \draw[arrow] (x_3) to node[above]{$\rho_6 \o \rho_6$} (x_2);
                \draw[arrow] (x_3) to[bend left] node[above]{$\rho_6  \o \rho_{67}$} (t_2);
                \draw[arrow] (t_1) to node[above]{$\rho_7 \o 1$} (x_3);
                \draw[arrow] (t_1) to node[above]{$1$} (t_2);

                \draw[arrow] (w_1) to node[below]{$1$} (w_2);
                \draw[arrow] (w_1) to node[below]{$\rho_{2} \o \rho_{12}$} (x_0);
                \draw[arrow] (w_1) to[bend right] node[below]{$1 \o \rho_{1}$} (x_1);
                \draw[arrow] (w_1) to[bend right=40] node[below]{$\rho_2 \o \rho_{123}$} (z_2);
                \draw[arrow] (x_0) to[bend left] node[below]{$\rho_1 \o 1$} (w_2);
                \draw[arrow] (x_0) to[bend right=36] node[below]{$1 \o \rho_3$} (z_2);
                \draw[arrow] (x_1) to node[below]{$\rho_2 \o \rho_2$} (x_0);
                \draw[arrow] (x_1) to[bend right] node[below]{$\rho_2  \o \rho_{23}$} (z_2);
                \draw[arrow] (z_1) to node[below]{$\rho_3 \o 1$} (x_1);
                \draw[arrow] (z_1) to node[below]{$1$} (z_2);

                \draw[arrow] (w_2) to[bend left] node[sloped,below,pos=0.1]{$\rho_{456} \o (\rho_1,\rho_{4567})$} (t_2);
                \draw[arrow] (w_2) to[bend left] node[sloped,above]{$\rho_{4} \o (\rho_1,\rho_4,\rho_7)$} (t_1);
                \draw[arrow] (w_2) to[bend left] node[sloped,above]{$\rho_{45} \o (\rho_1,\rho_{45})$} (x_3);
                \draw[arrow] (w_2) to[bend left] node[sloped,above]{$\rho_{4} \o (\rho_1,\rho_{4}) + \rho_{456} \o (\rho_1,\rho_{456})$} (x_2);
                \draw[arrow] (w_2) to[bend left] node[sloped,below,pos=0.1]{$\rho_{45} \o (\rho_1,\rho_{4})$} (c_1);

                \draw[arrow] (w_1) to node[sloped,above,pos=0.9]{$\rho_{23456} \o \rho_{1234567}$} (t_2);
                \draw[arrow] (w_1) to node[sloped,above,pos=0.8]{$\rho_{234} \o (\rho_{1234},\rho_{7})$} (t_1);
                \draw[arrow] (w_1) to node[sloped,above, pos=0.3]{$\rho_{2345} \o \rho_{12345}$} (x_3);
                \draw[arrow] (w_1) to node[sloped,above,pos=0.3]{$\rho_{23456} \o \rho_{123456} + \rho_{234} \o \rho_{1234}$} (x_2);
                \draw[arrow] (w_1) to node[sloped,above,pos=0.1]{$\rho_{2345} \o \rho_{1234}$} (c_1);

                \draw[arrow] (x_0) to node[sloped,below,pos=0.1]{$\rho_{3456} \o \rho_{34567}$} (t_2);
                \draw[arrow] (x_0) to node[sloped,above,pos=0.22]{$\rho_{34} \o (\rho_{34},\rho_7)$} (t_1);
                \draw[arrow] (x_0) to node[sloped,above,pos=0.8]{$\rho_{345} \o \rho_{345}$} (x_3);
                \draw[arrow] (x_0) to node[sloped,above,pos=0.2]{$\rho_{34} \o \rho_{34} + \rho_{3456} \o \rho_{3456}$} (x_2);
                \draw[arrow] (x_0) to node[sloped,below,pos=0.8]{$\rho_{345} \o \rho_{34}$} (c_1);

                \draw[arrow] (x_1) to node[sloped,below,pos=0.2]{$\rho_{23456} \o \rho_{234567} + \rho_{456} \o \rho_{4567}$} (t_2);
                \draw[arrow] (x_1) to node[sloped,below,pos=0.7]{$\rho_{4} \o (\rho_{4},\rho_7) + \rho_{234} \o (\rho_{234},\rho_7)$} (t_1);
                \draw[arrow] (x_1) to node[sloped,below]{$\rho_{45} \o \rho_{45} + \rho_{2345} \o \rho_{2345}$} (x_3);
                \draw[arrow] (x_1) to node[sloped,above]{$\rho_{234} \o \rho_{234} + \rho_{23456} \o \rho_{23456} + \rho_{456} \o \rho_{456}+\rho_{4} \o \rho_{4}$} (x_2);
                \draw[arrow] (x_1) to node[sloped,above,pos=0.8]{$\rho_{45} \o \rho_{4} + \rho_{2345} \o \rho_{234}$} (c_1);

                \draw[arrow] (z_2) to node[sloped,below,pos=0.1]{$\rho_{3456} \o \rho_{4567}$} (t_2);
                \draw[arrow] (z_2) to node[sloped,above,pos=0.3]{$\rho_{34} \o (\rho_{4},\rho_7)$} (t_1);
                \draw[arrow] (z_2) to node[sloped,above,pos=0.3]{$\rho_{345} \o \rho_{45}$} (x_3);
                \draw[arrow] (z_2) to node[sloped,above,pos=0.3]{$\rho_{34} \o \rho_{4} + \rho_{3456} \o \rho_{456}$} (x_2);
                \draw[arrow] (z_2) to node[sloped,below,pos=0.2]{$\rho_{345} \o \rho_{4}$} (c_1);

                % \draw[arrow] (x_2) to[bend right=60] node[above,sloped,pos=0.5,xshift=1cm]{$\rho_{56} \o 1$} (t_8);
            \end{tikzpicture}
            \caption{Bimodule ${}^{\B(\Zpmc_2)}[\mathbb{I}]^b_{ \B(\Zpmc_2)}$.}
            \label{fig:id_bounded_bimodule}
        \end{figure}

\Needspace{8\baselineskip}
\section{Background on fixed point Floer cohomology}\label{sec:fixed point Floer}
    In this section, we sketch the definition of \emph{fixed point Floer homology}, and describe the existing computational methods.

    Fixed point Floer cohomology was initially defined by Floer in~\cite{Flo} for symplectomorphisms which are Hamiltonian isotopic to the identity. It was extended to other symplectomorphisms by Dostoglou and Salamon in~\cite{DS}. Seidel in~\cite{Sei1} studied fixed point Floer cohomology of Dehn twists on surfaces. In~\cite{Sei2} Seidel defined fixed point Floer cohomology for any mapping class $\phi$ of a surface, proving that the choice of the symplectic representative of $\phi$ does not matter. 

    \Needspace{8\baselineskip}
    \subsection{Case of a closed surface}
        Consider a closed oriented surface $\Sigma$ with genus $g>1$. The construction of fixed point Floer cohomology $HF(\phi)$ for orientation preserving mapping classes $\phi \in MCG(\Sigma)$ works as follows:

        \begin{enumerate}

        \item We first choose a symplectic area form $\omega$ on $\Sigma$, and an area-preserving monotone representative $\phi$ of a mapping class. See~\cite{Sei2} for the definition of $monotone$ representative, and for the technical discussion on why Floer cohomology does not depend on the choice of such representative. 

        \item Next, we want to make sure that $\phi$ has only non-degenerate fixed points, i.e. at the fixed points we want $\det(d\phi - \id )\neq 0$. This is usually achieved by perturbing $\phi$ by the time-one isotopy $\psi^1_{X_{H_t}}$ along the Hamiltonian vector field $X_{H_t}$, where $H_t:\Sigma \rightarrow {\mathbb R}$ is a time-dependent generic Hamiltonian.

        \item Now, the chain complex $CF_*(\phi)$ is generated over $\mathbb{F}_2$ by the fixed points of $\phi$. Non-degeneracy of fixed points implies that they are isolated, and so $CF_*(\phi)$ is finitely generated. Note that the fixed points of $\phi$ correspond to the constant sections of the fiber bundle $T_{\phi}\xrightarrow{\Sigma} S^1 $, where $T_{\phi}$ is the mapping torus $T_{\phi}=\Sigma \times [0,1] / (\phi(p),0)\sim (p,1)$.  

        \item Next we need to pick an almost complex structure $J$ on $T_{\phi} \times {\mathbb R}$. First, we pick a generic time-dependent almost complex structure $j$ on $\Sigma$, and then extend it to the rest of the tangent space of $T_{\phi} \times {\mathbb R}$ naturally, i.e. the direction of the circle inside $T_{\phi}$ and the direction of $\mathbb{R}$ are interchanged by $J$. 

        \item The differential $\partial: CF_*(\phi) \rightarrow CF_{*-1}(\phi)$ is now defined by counting the number of points in the moduli space  $\mathcal M_1(x,y)$ of pseudo-holomorphic cylinder sections of the fiber bundle $T_\phi \times {\mathbb R} \xrightarrow{\Sigma} S^1 \times  {\mathbb R} $, which limit to the constant sections $x$ and $y$ of $T_{\phi}\xrightarrow{\Sigma} S^1 $ at $+\infty$ and $-\infty$ respectively. Namely, if $x_i$ are all constant sections of $T_{\phi}\xrightarrow{\Sigma} S^1 $ (remember that these correspond to fixed points of $\phi$), then
            $$\partial (x_j) = \sum_{x_i} \#( \mathcal M_1(x_j,x_i) /{\mathbb R} ) \cdot x_i.$$
        Note that only the index one cylinders are counted (i.e. those which come in $1$-dimensional family), up to translation along $\R$. The differential goes from $+\infty$ to $-\infty$, as in Morse homology. 

        \item This differential satisfies $\partial^2=0$, and passing to the homology of the dual complex $HF^*(\phi)=H(CF^*(\phi),d)$ gives \emph{fixed point Floer cohomology} --- an invariant, which depends only on the mapping class $\phi \in MCG(\Sigma)$. The $\mathbb{Z}_2$-grading on this invariant is provided by the sign of $\det(d\phi - \id )$ at fixed points.
        \end{enumerate}

        In the construction above we omitted a lot of deep and hard pseudo-holomorphic  theory: one has to study the moduli space $\mathcal M(x_j,x_i)$, prove that the sum $\partial (x_j) = \sum_{x_i} \#( \mathcal M(x_j,x_i) /{\mathbb R} ) \cdot x_i$ is finite, prove that $\partial^2=0$, and also prove that the homotopy equivalence class of $CF^*(\phi)$ does not depend on all of the choices made: generic area-preserving monotone representative $\phi$, and the complex structure $J$. We refer the reader to papers~\cite{Gau,CC,Sei3II,Ul} and references therein for discussion of these issues, and further details.
    \subsection{Case of a surface with boundary}
\begin{figure}[!ht]
        \centering
            \includegraphics[width=0.5\textwidth]{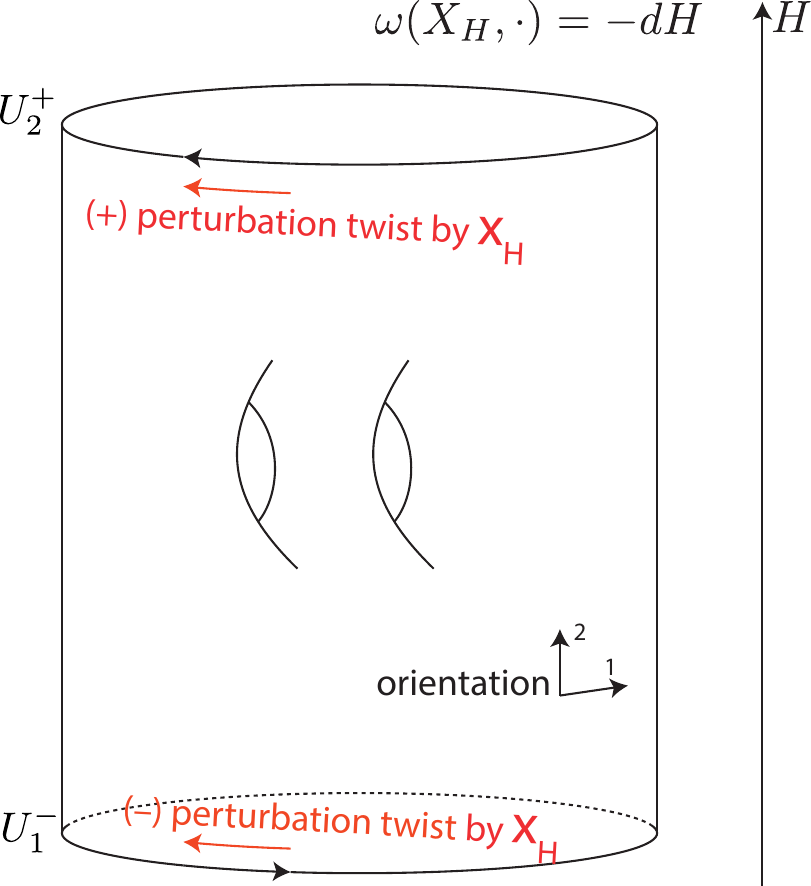}
            \caption{Perturbation twists near the boundary.}
            \label{fig:perturbation_conventions}
        \end{figure}

        There is a natural generalization of the above construction to surfaces with boundary. Suppose $\Sigma$ is an oriented surface of any genus with $n\neq 0$ boundary components $U_1 \cup U_2 \cup \dots \cup U_n$. 
        %rmk: genus g=1 is ok if there is boundary, because maps S^1xS^1 to Sigma_g have zero degree. Also, one can just use for example Ulrajevic work for full generality.
        We will consider orientation preserving mapping classes $\phi \in MCG_0(\Sigma)$, pointwise fixing the boundary. 
        \begin{enumerate}
            \item We first choose a symplectic area form $\omega$ on $\Sigma$, and an exact area-preserving representative $\phi$ of a mapping class. See~\cite[Appendix C]{Gau}, or~\cite[Lemma~3.3]{Ul} for an explanation of why the construction does not depend on this choice. 

            \item As in the closed case, we want every fixed point to be non-degenerate, and therefore isolated. Here comes the key difference from the closed surface case: because $\phi|_{\partial \Sigma}=\id|_{\partial \Sigma}$, we will need to perturb $\phi$ near the boundary. Thus, in order to specify a perturbation, as an input we will also take decorations of every boundary component with a sign, which will tell us how the perturbation behaves near the boundary components. If $U_i$ is decorated by $(+)$, then the perturbation in the neighborhood of $U_i$ should be a twist in the direction of the natural orientation of $U_i$ (or, equivalently, along the Reeb flow on the boundary, in the terminology of contact geometry). This corresponds to perturbation along a Hamiltonian vector field with the Hamiltonian $H$ having a time-independent local maximum on $U_i$, see Figure~\ref{fig:perturbation_conventions}. If $U_j$ is decorated by $(-)$, then the perturbation should be a twist in the opposite direction of the natural orientation, i.e. the Hamiltonian should have a local minimum on $C_j$. These twists near the boundary should be small enough, i.e. $\le 2\pi$ if one full twist is $2\pi$. The unions of positively and negatively decorated components will be denoted
            \begin{equation}\label{dec}
            U^+ = \bigcup_{\substack{\text{positively}\\\text{decorated}}} U_i^+ \qquad \text{and} \qquad U^- = \bigcup_{\substack{\text{negatively}\\\text{decorated}}} U_j^-
            \end{equation}

            \item The rest of the construction is the same as in the closed surface case, and we denote the resulting fixed point Floer cohomology by 
            $$HF^*(\phi\:;\:  U^+,U^-).$$ 
        \end{enumerate}

        \begin{remark}
        Notice that the naming of the twists comes from comparing the direction of the twist with the orientation on the boundary. It is not related to positive or negative Dehn twists. In fact, positive $(+)$ direction of twisting corresponds to the left handed twisting, which appears in negative (left handed) Dehn twists. $(-)$ direction of twisting corresponds to the right handed twisting, which appears in positive (right handed) Dehn twists. 
        \end{remark}

        % \begin{remark}
        % Note that paths of exact symplectomorphisms are Hamiltonian isotopies, see~\cite[Lemma~2.33]{Ul}.
        % \end{remark}

        % \begin{remark}
        % In case of closed surfaces, one needs to require genus to be more than one, so that all loops of free loops are null homotopic, i.e. maps $S^1 \time S^1 \rightarrow \Sigma$ are homotopic to a map $S^1 \time S^1 $ into some embedded curve $C \subset \Sigma$. This is used in the definition of Floer homology.
        % \end{remark}

        % see also the following papers:
        % Every area preserving self-diffeomorphism of $\hat{\Sigma}$ can be deformed to an exact one, see~\cite[Lemma~1.1]{BEE}. 
        % Every path of exact symplectomorphism is a Hamiltonian deformation, see~\cite[Lemma~3.1.1]{ER}.
        % Space of exact symplectomorphisms isotopic to identity is path-connected, see~\cite[Lemma~3.5.1]{ER}

    \Needspace{8\baselineskip}
    \subsection{Existing computational methods}\label{subsec:exist_comp_methods}
        In the case of the identity mapping class the Floer cohomology is the same as Morse cohomology with respect to the Hamiltonian we use for perturbing $\id$. See~\cite[Lemma~3.9]{Sei3II} for a proof and references to the original works. Thus we have 
        $$HF^*(\id\:;\:  U^+,U^-)=MH^*(\Sigma \:;\ H \! \! \uparrow^{U^+}_{U^-}),$$where by $MH^*$ we denote Morse cohomology. It is equal to the relative homology 
        $$H^*(\Sigma, \ U^-),$$ 
        because the Hamiltonian has local minimum on $U^-$ and local maximum on $U^+$.

        First computations for non-trivial mapping classes were done by Seidel in~\cite{Sei1}. Suppose $\phi$ is a composition of right handed Dehn twists along curves $R=\{R_1,R_2,\dots,R_l\}$ and left handed Dehn twists along curves $L=\{L_1,L_2,\dots,L_k\}$. Suppose that all the curves are disjoint, their complement has no disc components, and that no $L_i$ is homotopic to $R_j$. Then
        $$HF^*(\phi\:;\:  U^+,U^-)=H^*(\Sigma - L, \ R \cup U^-)$$ 
        for arbitrary decorations of boundary components, as in~\eqref{dec}.
        This result is again achieved via reducing the computation to the Morse cohomology, where the Hamiltonian has local minimum on the curves $R_i$ and local maximum on the curves $L_j$.

        Then Gautschi in~\cite{Gau} computed fixed point Floer cohomology for algebraically finite mapping classes (those include periodic mapping classes, and also reducible mapping classes which restrict to periodic classes to each component). Eftekhary in~\cite{Eft} then generalized Seidel's work to Dehn twists along curves $\{R_1,R_2,\dots,R_l\}$ which form a forest, i.e. the intersection graph (with $R_i$ being vertices and intersection points $R_i\cap R_j$ being edges) does not contain cycles; this is the most useful result for us, and so we quote it in detail below. The last computations were done by Cotton-Clay in~\cite{CC}, where he showed how to compute fixed point Floer cohomology  for all pseudo-Anosov mapping classes and for all reducible mapping classes (including those with pseudo-Anosov components). Thus there is a way to compute fixed point Floer cohomology for any mapping class. 

        The following is a generalization of Eftekhary's work to the case with boundary.
        \begin{theorem}[Eftekhary]\label{Eftekhary's theorem}
        Suppose $\Sigma$ is a surface, whose boundary, if non-empty, is divided into positively and negatively decorated components $U=U^+ \sqcup U^-$, as in~\eqref{dec}. Suppose $\phi\co MCG_0(\Sigma, \partial \Sigma = U)$ is a mapping class equal to a composition of right handed Dehn twists along the curves $R=\{R_1,R_2,\dots,R_l\}$, and left handed Dehn twists along the curves $L=\{L_1,L_2,\dots,L_m\}$. Assume that 
        \begin{itemize}
        \item $R$ is a forest;
        \item $L$ is a forest;
        \item $L \cap R = \emptyset $;
        \item no $L_i$ is homotopic to $R_j$; 
        \item all the curves $L_i$  and  $R_j$ are homologically essential.
        \end{itemize}   
        Then
        $$
        HF^*(\phi\:;\:  U^+,U^-)=H^*(\Sigma - L, \ R \cup U^-).
        $$
        \end{theorem} 
        We will use this theorem to perform computations of fixed point Floer cohomology in Section~\ref{supporting_comps}.
        
\Needspace{8\baselineskip}
\section{Conjectural isomorphism}\label{sec:conj_iso}
    \Needspace{8\baselineskip}
    \subsection{Statement}\label{sec:statement}
        Here we state our main conjecture, around which the paper is revolving. As in Section~\ref{sec:bimodule from bordered}, we consider the strongly based mapping class group of the genus $g$ surface with one boundary component. Because we want to be able to take fixed point Floer cohomology of $\phi$, we assume the genus to be greater than one.

        Having a mapping class $\phi\co  (\Sigma,\partial \Sigma=U_1=S^1) \rightarrow (\Sigma,\partial \Sigma=U_1=S^1)$, let us consider the induced mapping class $\tilde{\phi}\co (\tilde{\Sigma},\partial \tilde{\Sigma}=U_1 \cup U_2) \rightarrow (\tilde{\Sigma},\partial \tilde{\Sigma}=U_1 \cup U_2)$, where $\tilde{\Sigma}=\Sigma \setminus D^2$ is obtained by removing a disc in the small enough neighborhood of the boundary, such that $\phi$ is identity on that neighborhood.

        \begin{conjecture}\label{conj}
        For every mapping class $\phi \in MCG_0(\Sigma,\partial\Sigma=S^1=U_1)$ there is an isomorphism of $\mathbb{Z}_2$-graded vector spaces
        $$
        HH_*(N(\phi^{-1})) \cong HF^{*+1}(\tilde{\phi}\:;\:  U_2^+,U_1^-).$$
        \end{conjecture}

        In the rest of the paper we will be investigating the conjecture from different angles. First, from the practical and concrete point of view, by performing lots of computations in Section~\ref{supporting_comps} we will obtain a strong evidence that the conjecture is true. Then, in Section~\ref{sec:bordered<->Fuk}, we will outline the symplectic geometric interpretation of bordered Heegaard Floer theory. Based on that, in Section~\ref{sec:LF}, we will reinterpret Conjecture~\ref{conj} in the context of the Fukaya category: we will see that it is a special case of a more general conjecture of Seidel ~\cite[Conjecture 7.18]{Sei3II}, which states that the open-closed map in the context of the Fukaya-Seidel category of Lefschetz fibration is an isomorphism.

        % \begin{remark}
        % The reason that on the left side of the isomorphism we have $\phi^{-1}$ is because the bimodule coming from the bordered theory is homotopy equivalent to $\bigoplus\limits_{i,j} hom_{\F_z(\Sigma)} (\alpha_i,\phi(\alpha_j))$ (see Section~\ref{sec:bordered<->Fuk} for the explanation of this), and the bimodule defined in~\cite{Sei3II} is $\bigoplus\limits_{i,j} hom_{\F_z(\Sigma)} (\phi(\alpha_i),\alpha_j)\simeq \bigoplus\limits_{i,j} hom_{\F_z(\Sigma)} (\alpha_i,\phi^{-1}(\alpha_j))$, see Section~\ref{subsec:O-C map} for more on this.
        % \end{remark}

        % \begin{remark}
        % Let us also point out that the ranks of $HF^*(\tilde{\phi}\:;\:  U_2^+,U_1^-)$ and $HF^*(\tilde{\phi}^{-1}\:;\:  U_2^+,U_1^-)$ should be the same by Poincare duality for Floer homology.
        % \end{remark}
    \Needspace{8\baselineskip}
    \subsection{Supporting computations}\label{supporting_comps}
        To support Conjecture~\ref{conj}, we perform computations in the genus two case. From the bimodule side, the key results and tools we use are the computation of the Dehn twist bimodules in Section~\ref{subsec:bim_comps}, and the computer program~\cite{Pyt} we wrote to tensor bimodules and compute their Hochschild homologies. From the fixed point Floer cohomology side, we heavily exploit Theorem~\ref{Eftekhary's theorem}. We did lots of computations, and all of them confirmed the conjecture; below we describe five of them, in each case showing that the invariants $HH_*(N(\phi^{-1}))$ and $HF^{*+1}(\tilde{\phi}\:;\:  U_2^+,U_1^-)$ are isomorphic.

        As in Section~\ref{subsec:bim_comps}, we fix a set of curves generating the mapping class group as in Figure~\ref{fig:Sigma_2}, and use a parameterization $\Sigma_2 \cong F^\circ(\Zpmc_2)$ as in Figure~\ref{fig:Sigma_2_parameterized}.

        \begin{computation}[$\phi=\id$]\label{comp:id}
            Taking the identity bimodule $N(\id)=[\mathbb I]$ as an input, the program~\cite[Showcase~4]{Pyt} finds that Hochschild homology is 
            $$HH_*(N(\id))=(\mathbb{F}_2)^4_{(0)},$$ 
            concentrated in grading 0. Applying Theorem~\ref{Eftekhary's theorem}, we obtain that fixed point Floer cohomology in this case is also four-dimensional: 
            $$HF^*(\tilde{\id}\:;\:   U_2^+,U_1^-)=MH^*( \tilde{\Sigma}_2;H \! \! \uparrow^{U_2^+}_{U_1^-})=H^*(\tilde{\Sigma}_2, \ U_1)=(\mathbb{F}_2)^4_{(1)},$$ 
            concentrated in grading 1, see Figure~\ref{fig:HF-id-} for an illustration. These computations confirm Conjecture~\ref{conj} for $\phi=\id\co(\Sigma_2,\partial \Sigma_2) \rightarrow (\Sigma_2,\partial \Sigma_2)$.
            \begin{figure}[H]
            \centering
                \includegraphics[width=1\textwidth]{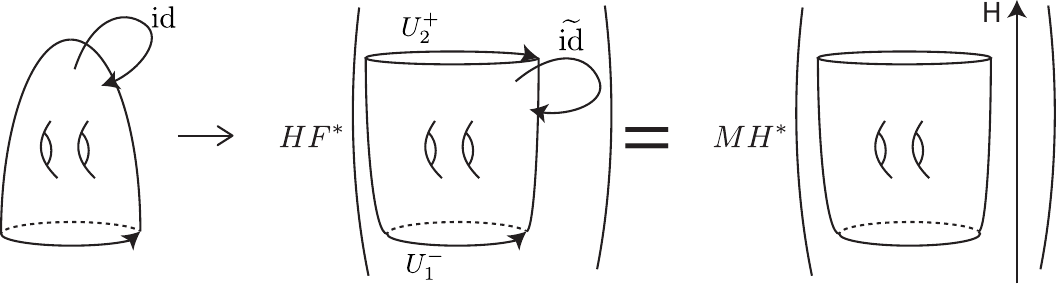}
                \caption{Computation of $HF^*(\tilde{\id}\:;\:   U_2^+,U_1^-)=MH^*( \tilde{\Sigma}_2;H \! \! \uparrow^{U_2^+}_{U_1^-})=H^*(\tilde{\Sigma}_2, \ U_1)=(\mathbb{F}_2)^4$ for $\id\co(\Sigma_2,\partial \Sigma_2) \rightarrow (\Sigma_2,\partial \Sigma_2)$.}
                \label{fig:HF-id-}
            \end{figure}
        \end{computation}

        \begin{computation}[$\phi=\tau_l$]\label{comp:DehnTwist}
            Suppose $\phi = \tau_l$ is a right handed Dehn twist along any of the curves $A,~B,~C,~D,$ or $E$. Then the program~\cite[Showcase~5]{Pyt} finds that Hochschild homology is 
            $$HH_*(N(\tau_l^{-1}))=(\mathbb{F}_2)^4_{(0)}.$$
            Applying Theorem~\ref{Eftekhary's theorem}, we obtain that fixed point Floer cohomology in this case is 
            $$HF^*(\tilde{\tau_l}\:;\:  U_2^+,U_1^-)\overset{(1)}{=}H^*(\tilde{\Sigma}_2,l \cup U_1)=(\mathbb{F}_2)^4_{(1)}.$$ 
            Equality $(1)$ is obtained by cutting $\tilde{\Sigma}_2$ along $l$ and computing Morse cohomology with respect to Hamiltonian which looks like in Figure~\ref{fig:HF-one_right_twist-}. These computations confirm Conjecture~\ref{conj} for $\phi=\tau_l:(\Sigma_2,\partial \Sigma_2) \rightarrow (\Sigma_2,\partial \Sigma_2)$.

            \begin{figure}[H]
            \centering
                \includegraphics[width=0.88\textwidth]{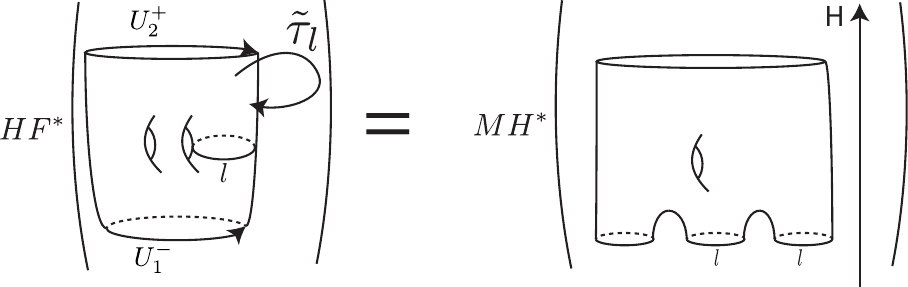}
                \caption{Computation of $HF^*(\tilde{\tau_l}\:;\:  U_2^+,U_1^-)=MH^*(\cdots)=H^*(\tilde{\Sigma}_2,l \cup U_1)=(\mathbb{F}_2)^4,$ for the right handed Dehn twist $\tau_l:(\Sigma_2,\partial \Sigma_2) \rightarrow (\Sigma_2,\partial \Sigma_2)$.}
                \label{fig:HF-one_right_twist-}
            \end{figure}
        \end{computation}

        \begin{computation}[$\phi=\tau_l^{-1}$]\label{comp:invDehnTwist}
            For left handed Dehn twists we have the same answers: the program~\cite[Showcase~6]{Pyt} finds 
            $$HH_*(N(\tau_l))=(\mathbb{F}_2)^4_{(0)},$$ 
            and  Theorem~\ref{Eftekhary's theorem} implies 
            $$HF^*(\tilde{\tau_l}^{-1}\:;\:  U_2^+,U_1^-)\overset{(1)}{=}H^*(\tilde{\Sigma}_2-l,U_1)=(\mathbb{F}_2)^4_{(1)}.$$
            Equality $(1)$ is obtained by cutting $\tilde{\Sigma}_2$ along $l$ and computing Morse cohomology with respect to Hamiltonian which looks like in Figure~\ref{fig:HF-one_left_twist-}. These computations confirm Conjecture~\ref{conj} for $\phi=\tau_l^{-1}\co(\Sigma_2,\partial \Sigma_2) \rightarrow (\Sigma_2,\partial \Sigma_2)$.

            \begin{figure}[H]
            \centering
                \includegraphics[width=0.88\textwidth]{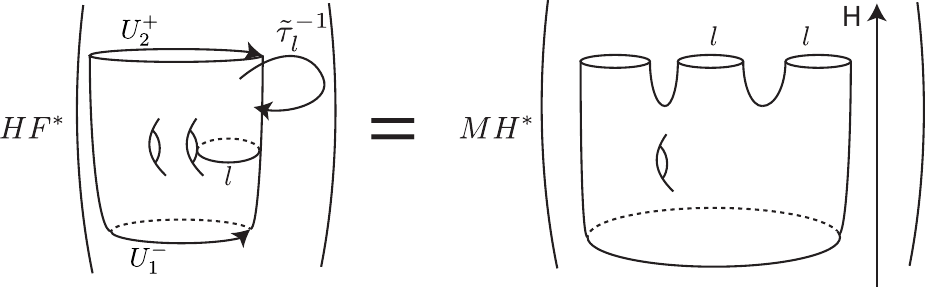}
                \caption{Computation of $HF^*(\tilde{\tau_l}^{-1}\:;\:  U_2^+,U_1^-)=MH^*(\cdots)=H^*(\tilde{\Sigma}_2-l,U_1)=(\mathbb{F}_2)^4$ for the left handed Dehn twist $\tau_l^{-1}\co(\Sigma_2,\partial \Sigma_2) \rightarrow (\Sigma_2,\partial \Sigma_2)$.}
                \label{fig:HF-one_left_twist-}
            \end{figure}
        \end{computation}

        \begin{proof}[Proof of Theorem~\ref{thm:int}]
        Based on the computations above, we prove the isomorphism $HH_*(N(\phi^{-1})) \cong HF^{*+1}(\tilde{\phi}\:;\:  U_2^+,U_1^-)$ for the identity mapping class $\phi=\id \lefttorightarrow \Sigma_2$, and also any Dehn twist $\phi=\tau  \lefttorightarrow \Sigma_2$. 

        Computation~\ref{comp:id} proves the conjecture for $\phi=\id$, while Computations~\ref{comp:DehnTwist} and~\ref{comp:invDehnTwist} prove the conjecture for Dehn twists along curves $A,~B,~C,~D,~E$ from Figure~\ref{fig:Sigma_2}. 

        We now prove the isomorphism for any Dehn twist along a \emph{non-separating} curve $m$. Any such curve $m$ can be mapped to the curve A by a mapping class $f$, and so any Dehn twist $\tau_m$ is conjugate to the Dehn twist $\tau_A$: $f\tau_A f^{-1}=\tau_{f(A)} =\tau_{m} $. Both fixed point Floer homology (it follows from the definition) and Hochschild homology (see Remark~\ref{rmk:inv_conj}) are invariant with respect to conjugation, and thus it follows that the isomorphism $HH_*(N(\tau_A^{-1})) \cong HF^{*+1}(\tilde{\tau_A}\:;\:  U_2^+,U_1^-)$ implies the isomorphism for any Dehn twist $\tau_m$.

        The case of a Dehn twist along a non-trivial \emph{separating} curve $m$ is done similarly. Any such curve $m$ either separates $\Sigma_2$ into two genus one surfaces, or into a genus two surface and an annulus. Therefore $m$ can be mapped to either the curve $G_1=\partial \big( Nbd_\epsilon (A\cup B) )\big)$, or the boundary curve $ G_2=\partial \big( Nbd_\epsilon(A\cup B \cup C \cup C) \big) \sim \partial \Sigma_2$. Thus it is enough to prove the conjecture for Dehn twists along $G_1$ and $G_2$. Seidel's result~\cite{Sei1} implies 
        \[ HF^{*}(\tilde{\tau^{\pm 1}_{G_1}}\:;\:  U_2^+,U_1^-) = (\mathbb F_2)^5_{(1)}\oplus (\mathbb F_2)_{(0)}, \qquad   HF^{*}(\tilde{\tau^{\pm 1}_{G_2}}\:;\:  U_2^+,U_1^-) =  (\mathbb F_2)^5_{(1)}\oplus (\mathbb F_2)_{(0)}. \]
        By leveraging the chain relations in the mapping class group $\tau_{G_1}=(AB)^6, \ \tau_{G_2}=(ABCD)^{10}$, the program~\cite[Showcase~7]{Pyt} finds
        \[ HH_*(N(\tau^{\pm 1}_{G_1}))=  (\mathbb F_2)^5_{(0)}\oplus (\mathbb F_2)_{(1)}, \qquad   HH_*(N(\tau^{\pm 1}_{G_2})) = (\mathbb F_2)^5_{(0)}\oplus (\mathbb F_2)_{(1)}, \]
        which finishes the proof of Theorem~\ref{thm:int} in case of a Dehn twist along a separating curve.
        \end{proof}

        Experimenting with mapping classes arising from two disjoint forests of curves (i.e. those where we can use Eftekhary's Theorem~\ref{Eftekhary's theorem}), we always got equal ranks of homologies. Let us highlight two more examples.

        \begin{computation}[$\phi=\tau_A \tau_B \tau_C \tau_D$]
            This mapping class is a monodromy of an open book on $S^3$ with the binding being the $(5,2)$ torus  knot, and the page being the genus two surface with boundary. Either by remembering that $\widehat{HFK}(S^3, T_{(5,2)})=(\mathbb{F}_2)^5$ contributing $\mathbb{F}_2$ to each of the five Alexander gradings $-2,-1,0,1,2$, or by invoking the program~\cite[Showcase~8]{Pyt}, we get 
            $$\widehat{HFK}(T_{(5,2)} \:;\:  -g+1)=HH_*(N(\tau_D^{-1} \tau_C^{-1} \tau_B^{-1} \tau_A^{-1}))=(\mathbb{F}_2)_{(0)}.$$ 
            Applying Theorem~\ref{Eftekhary's theorem} gives 
            $$HF^*(\widetilde{\tau_A \tau_B \tau_C \tau_D}\:;\:  U_2^+,U_1^-)=H^*(\tilde{\Sigma}_2,U_1 \cup A \cup B \cup C\cup D)=(\mathbb{F}_2)_{(1)},$$ confirming Conjecture~\ref{conj}. It is the lowest rank that we observed in our computations. 
        \end{computation}
            \begin{remark}
            Interestingly, it is known that $\rk(\widehat{HFK}(\text{fibered knot} \:;\:  -g+1))=\rk(HH_*(N(\phi)))>0$
            is always true, see~\cite{BV}. 
            \end{remark}

        \begin{computation}[$\phi=\tau_A^5 \tau_B \tau_C \tau_D  \tau_E^5$]\label{ps-a}
            This mapping class is pseudo-Anosov if viewed as a mapping class of a closed genus two surface (see~\cite{Eft}). In this case the program~\cite[Showcase~9]{Pyt} finds 
            $$HH_*(N(\phi^{-1}))=(\mathbb{F}_2)^{9}_{(1)}\oplus (\mathbb{F}_2)_{(0)},$$ 
            and Theorem~\ref{Eftekhary's theorem} implies 
            $$HF^*(\tilde{\phi}\:;\:  U_2^+,U_1^-)=H^*(\tilde{\Sigma}_2,U_1 \cup 5A \cup B \cup C \cup D \cup 5E)=(\mathbb{F}_2)^{9}_{(0)}\oplus (\mathbb{F}_2)_{(1)},$$
            confirming Conjecture~\ref{conj}.
        \end{computation}

\Needspace{8\baselineskip}
\section{Background on bordered theory vs Fukaya categories}\label{sec:bordered<->Fuk}
    This section covers background material, which is needed for the subsequent Section~\ref{sec:LF}. 

    Let us recall how we associate a bimodule to a mapping class in Section~\ref{sec:hf_bim}. To a surface we associate a $dg$ algebra. To a mapping class we associate an $A_\infty$ bimodule over that algebra. Composition of mapping classes corresponds to the box tensor product of bimodules (Theorem~\ref{DA pairing}). In Section~\ref{alt_bim} below, we describe another geometric construction of the same mapping class invariant. The basis for that construction is reinterpretation of bordered Heegaard Floer theory in terms of the partially wrapped Fukaya category of the surface, which is due to Auroux~\cite{Aur1,Aur2}, and which we recall in Section~\ref{aur_constr}.
    \Needspace{8\baselineskip}
    \subsection{Auroux's construction}\label{aur_constr}
\begin{figure}[!ht]
        \centering
            \includegraphics[width=0.5\textwidth]{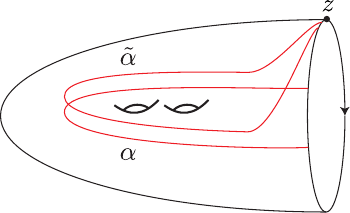}
            \caption{Perturbation near infinity for the partially wrapped Fukaya category of a surface.}
            \label{fig:partially_wrapped_perturbation}
        \end{figure}

        Consider a surface $\Sigma$ with one or more boundary components. Fix a set of points (\say{stops}) $Z$  on the boundary $ \partial \Sigma$, and fix a Liouville domain structure on $\Sigma$; the latter consists of an exact symplectic form $\omega = d \theta$, such that the Liouville vector field $X_\theta$ (defined by equation $\omega(X_\theta,\cdot)=\theta$) points outwards the boundary. Assume that every boundary component has at least one point from $Z$. To all this data we can associate the partially wrapped Fukaya category $\F_Z(\Sigma)$ (Auroux also considers Fukaya categories of $Sym^k(\Sigma)$, but we only need the $k=1$ case). We break down the construction of  $\F_Z(\Sigma)$ into several steps:
        \begin{enumerate}
        \item We first need to pass from $\Sigma$ to $\hat{\Sigma}$, which is a Liouville manifold obtained by completion of the surface $\Sigma$ by a cylindrical end. In more details, $\hat{\Sigma}$ is constructed by taking the symplectization of the boundary $([0,+\infty) \times \partial \Sigma, d(r \cdot \theta) )$, and gluing its $0 \le r \leq 1$ part to the Liouville-flow-collar neighborhood of $\partial \Sigma$, by $i:((0,1] \times \partial \Sigma, d(r \cdot \theta) ) \rightarrow ((-\infty,0]\times  \partial \Sigma, \omega) \subset \Sigma $, s.t. $i((e^r,x))=(r,x)$.  

        \item Objects of the partially wrapped Fukaya category $\F_Z(\Sigma)$ are curves of two types:
            \begin{itemize} 
            \item Closed exact embedded Lagrangian curves in $\hat{\Sigma}$; 
            \item Non-compact properly embedded curves (arcs for short) such that the ends at infinity of $\hat{\Sigma}$ stabilize to be rays $(r,0)$ in the cylindrical end.   
            \end{itemize}
        \item Morphism spaces are Lagrangian Floer cochain complexes $hom_{\F_Z(\Sigma)}(L_1,L_2)=CF_\text{Lagr}^*(\widetilde{L_1},L_2)$, where $\widetilde{L_1}$ is a Lagrangian submanifolds perturbed by a generic Hamiltonian. Because we consider non-compact Lagrangians, the behavior of the Hamiltonian perturbation at infinity of $\hat{\Sigma}$ affects $hom_{\F_Z(\Sigma)}(L_1,L_2)$ in an essential way. Auroux constructs a specific Hamiltonian, which wraps the ray of the arc around the cylindrical end until it reaches the stop, i.e. one of the rays in $Z \times [1,+\infty)$; the best way to understand the resulting perturbation is to look at Figure~\ref{fig:partially_wrapped_perturbation}, where $Z=\{z\}$. See~\cite{Aur1} for the details of the construction of this Hamiltonian.
        \item $A_\infty$ operations 
        \begin{equation}\label{a_inf_op}
        \mu_d: hom_{\F_Z(\Sigma)}(L_0,L_1)\o \dots \o
        hom_{\F_Z(\Sigma)}(L_{d-1},L_d) \rightarrow
        hom_{\F_Z(\Sigma)}(L_0,L_d)
        \end{equation}
        are given by counting holomorphic discs with $d+1$ marked points on the boundary (as usual in the Fukaya categories). In this step consistent perturbations and complex structures are needed to be picked, see~\cite{Aur1} for the details.
        \item The final step is to prove that $A_\infty$ operations satisfy $A_\infty$ relations, and that the quasi-equivalence class of $\F_Z(\Sigma)$ does not depend on all the choices made in the definition (perturbations and complex structures). Again, for this we refer the reader to the original article~\cite{Aur1}.
        \end{enumerate} 

        As shown in~\cite[Theorem~1]{Aur2}, if $\alpha_1,\dots,\alpha_k \in Ob(\F_Z(\Sigma))$ are non-intersecting arcs, and their complement is a set of discs each containing one point from $Z$, then they $generate$ $\F_Z(\Sigma)$. This means that every object $L\in Ob(\F_Z(\Sigma))$, considered as an object in certain enlarged category of twisted complexes $L\in Ob(Tw \F_Z(\Sigma))$, is quasi-isomorphic to a direct summand of iterated mapping cones between the generating objects $\alpha_1,\dots,\alpha_k$, see~\cite[Section~3]{Aur3} for the details. The deep theorem behind this fact  is that Lefschetz thimbles generate the Fukaya-Seidel category associated to a Lefschetz fibration~\cite[Theorem~18.24]{Sei5}. The Fukaya-Seidel category is closely related to the partially wrapped Fukaya category, see the next Section~\ref{sec:LF}.

        Notice that surface $F^\circ(\Zpmc)$ from Section~\ref{sec:bimodule from bordered}, associated to a genus $g$ pointed matched circle $\Zpmc$, contains a basepoint $z$ in it boundary, which can serve as a stop, giving rise to the Fukaya category $\F_z(F^\circ(\Zpmc))$. Moreover,  $F^\circ(\Zpmc)$ contains a distinguished set of generators of $\F_z(F^\circ(\Zpmc))$, which corresponds to the matched pairs of points in $\Zpmc$, see the surface on the right of Figure~\ref{fig:pointed_matched_circle}.

        \begin{theorem}[Auroux]\label{thm:aur}
        Suppose $\alpha_1,\dots,\alpha_{2g}$ are arcs in $F^\circ(\Zpmc)=\Sigma$ corresponding to the matched pairs in $\Zpmc$. Then the bordered one-strand-moving $dg$ algebra is quasi-isomorphic to the $A_\infty$ $hom$-algebra of the partially wrapped Fukaya category with respect to the generating set $\alpha_1,\dots,\alpha_{2g}$, i.e. 
        $$
        \B(\Zpmc) \simeq \A(\F_z(\Sigma)) \coloneqq  \bigoplus\limits_{1\leq i,j\leq 2g} hom_{\F_z(\Sigma)} (\alpha_i,\alpha_j).
        $$
        \end{theorem}

        \begin{remark}
        The full statement of Auroux's theorem involves all the summands of the bordered algebra and the Fukaya categories of symmetric products: for $0\leq k \leq 2g$ one has 
        $$ \A(\Zpmc,-g+k) \simeq \A(\F_z(Sym^k(\Sigma))) = \bigoplus\limits_{1\leq i,j\leq C_{2g}^k} hom_{\F_z(Sym^k(\Sigma))} (\lambda_i,\lambda_j),$$ where $\{ \lambda_i \}$ is a set of generators coming from the products of $k$ arcs in $\alpha_1,\dots,\alpha_{2g}$.
        \end{remark}

        Although we will not need this, let us mention that Auroux reinterpreted in terms of the Fukaya categories not only the bordered algebras, but also the type $A$ bordered Heegaard Floer modules, via the Yoneda embedding construction.

        % \begin{remark}
        % Having a finite set of generating objects $\alpha_1,\dots,\alpha_k$ is very useful, because the Yoneda embedding construction (see~\cite[Section~3.4]{Aur3}) then gives a fully faithful embedding of $A_\infty$ category into the category of $A_\infty$ modules over the $hom$-algebra $\bigoplus\limits_{1\leq i,j\leq k} hom_{\F_z(\Sigma)} (\alpha_i,\alpha_j)$. In particular, having a bordered Heegaard diagram for a bordered $3$-manifold, Auroux identified not only the algebras, but also the type A bordered Heegaard Floer module with a module coming from the Fukaya category via the Yoneda embedding construction.
        % \end{remark} 

        \begin{example}[Torus algebra]

        Let us illustrate the above theorem on a torus. The way to get the torus bordered algebra $\B(\Zpmc_1)$ from a pointed matched circle $\Zpmc_1$ is pictured in Figure~\ref{fig:Z_1--B-Z_1-}. The Fukaya categorical interpretation of the torus bordered algebra $\B(\Zpmc_1)$ is pictured in Figure~\ref{fig:bordered=partially_wrapped}. The elements of the algebra appear as generators of Lagrangian Floer complexes. Note that we cannot see the product structure (i.e. holomorphic triangles) on this picture, because for that we need to consistently pick perturbations for the three Lagrangians involved in the product operation.

        \begin{figure}[!ht]
        \centering
            \includegraphics[width=0.45\textwidth]{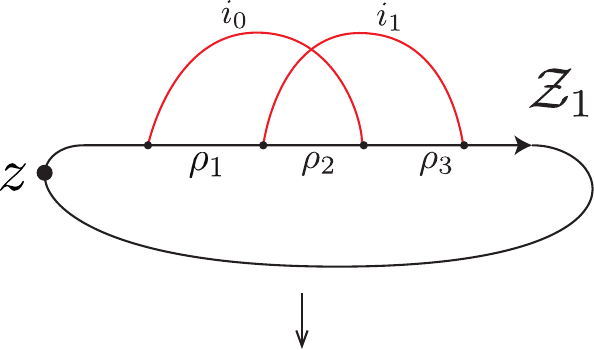}
            \tikzstyle{arrow} = [thick,->,>=stealth]

            \begin{tikzpicture}[scale=1.1,baseline=1.5cm]
            \node[circle](i_0) at (0,0){$i_0$};
            \node[circle](i_1) at (3,0){$i_1$};

            \draw[arrow] (i_0) to[bend left] node[above, sloped]{$\rho_1$} (i_1);
            \draw[arrow] (i_1) to node[above, sloped]{$\rho_2$} (i_0);
            \draw[arrow] (i_0) to[bend right] node[above, sloped]{$\rho_3$} (i_1);
            
            \draw (0,1) to[bend right] node[left,text width=3cm, text centered]{Path algebra \\ over $\mathbb{F}_2$} (0,-1);
            \draw (3,1) to[bend left] (3,-1);
            \draw (-2.5,-1.4) to (4,-1.4);
            \node at (1,-2){Relations $\rho_2 \rho_1=\rho_3 \rho_2=0$};
            \node[scale=1.2] at (-4.5,-1.4){$\B(\Zpmc_1)=$};

            \end{tikzpicture}
            \caption{Torus bordered algebra, constructed from a pointed matched circle.}
            \label{fig:Z_1--B-Z_1-}
        \end{figure}

        \begin{figure}[h]
        \centering
            \includegraphics[width=0.8\textwidth]{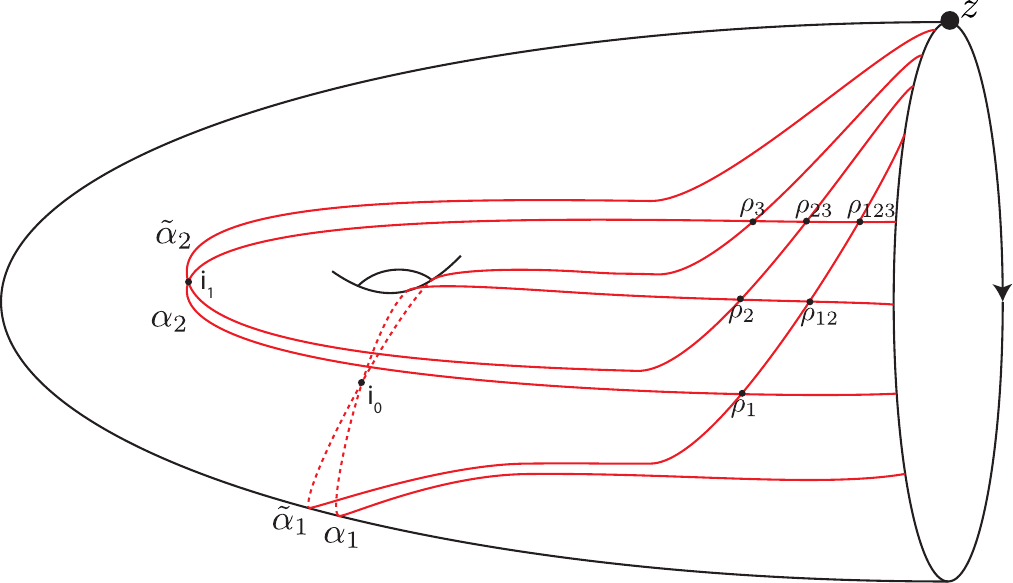}
            \caption{Elements of the torus bordered algebra $\B(\Zpmc_1)$, viewed as elements of morphism spaces between generating objects $\alpha_1,\alpha_2$ of the partially wrapped Fukaya category of the torus.}
            \label{fig:bordered=partially_wrapped}
        \end{figure}

        \end{example}
    \Needspace{8\baselineskip}
    \subsection{Alternative bimodule construction via Fukaya categories}\label{sec:bimodule_via_Fukaya}\label{alt_bim}
        Suppose we are given an exact self-diffeomorphism $\phi$ of $(\Sigma, \partial \Sigma = S^1)$, pointwise fixing the boundary. Suppose that the partially wrapped Fukaya category with one stop $\F_z(\Sigma)$ is generated by $\alpha_1,\ldots,\alpha_{2g}$. Then there is a standard way to associate to $\phi$ a $graph$ bimodule: $A_\infty$ bimodule of $AA$ type $_{\A(\F_z(\Sigma))}N_{\F_z}(\phi)_{\A(\F_z(\Sigma))}$ over the $hom$-algebra $\A(\F_z(\Sigma)) = \bigoplus\limits_{1\leq i,j\leq 2g} hom_{\F_z(\Sigma)} (\alpha_i,\alpha_j)$. 
        \begin{enumerate}
            \item The bimodule as a vector space is equal to 
            \begin{equation} \label{eq:bimodule from Fuk}
            N_{\F_z}(\phi)=\bigoplus\limits_{1\leq i,j\leq {2g}} hom_{\F_z(\Sigma)} (\alpha_i,\hat{\phi}(\alpha_j)),
            \end{equation}
            where  $\hat{\phi}$ is  an exact compactly supported self-diffeomorphism of the completion  $\hat{\Sigma}$ induced by $\phi$. From now on we will abuse notation and use $\phi$ for $\hat{\phi}$.

            \item The higher actions are given using $A_\infty$ operations~\ref{a_inf_op}. For example, the action $m^{1|1|1}\co\A(\F_z(\Sigma)) \o N_{\F_z}(\phi) \o \A(\F_z(\Sigma)) \rightarrow N_{\F_z}(\phi)$ is given via the following operation counting holomorphic discs with four marked points:  
            \begin{align*}
            {\scriptstyle
            \mu_3: \left( \bigoplus\limits_{1\leq i,j\leq {2g}} hom_{\F_z(\Sigma)} (\alpha_i,\alpha_j) \right) }
              & {\scriptstyle \o \left( \bigoplus\limits_{1\leq i,j\leq {2g}} hom_{\F_z(\Sigma)} (\alpha_i,\phi(\alpha_j)) \right) \o \left( \bigoplus\limits_{1\leq i,j\leq {2g}} hom_{\F_z(\Sigma)} (\phi(\alpha_i),\phi(\alpha_j)) \right)  \rightarrow  } \\ 
             & {\scriptstyle \rightarrow  \bigoplus\limits_{1\leq i,j\leq {2g}} hom_{\F_z(\Sigma)} (\alpha_i,\phi(\alpha_j)) }.
            \end{align*}
        \end{enumerate}
        Note that the bimodule $_{\A(\F_z(\Sigma))}N_{\F_z}(\phi)_{\A(\F_z(\Sigma))}$ is of $AA$ type in the bordered theory terminology, whereas $^{\B(\Zpmc)}N(\phi)_{\B(\Zpmc)}$ from Section~\ref{sec:bimodule from bordered} was of $DA$ type. Remember that if $F^\circ(\Zpmc)=\Sigma$, by Theorem~\ref{thm:aur} we have $\B(\Zpmc)\simeq \A(\F_z(\Sigma))$, and so both bimodules are over the same algebra $\B(\Zpmc)$. We now unify the two constructions of bimodules. 

        \begin{proposition}\label{prop:bims_are_eq}
        Suppose $\alpha_1,\dots,\alpha_{2g}$ are arcs in $F^\circ(\Zpmc)=\Sigma$ corresponding to the matched pairs in the pointed matched circle $\Zpmc$. Suppose $\phi \in MCG_0(F^\circ(\Zpmc),F^\circ(\Zpmc))$. Then the modulification of type $DA$  bimodule $^{\B(\Zpmc)}N(\phi)_{\B(\Zpmc)}$ is homotopy equivalent to $_{\B(\Zpmc)}N_{\F_z}(\phi)_{\B(\Zpmc)}$:
        $$ {}_{\B(\Zpmc)}{\B(\Zpmc)}_{\B(\Zpmc)} \boxtimes {}^{\B(\Zpmc)}N(\phi)_{\B(\Zpmc)} \simeq {}_{\B(\Zpmc)}N_{\F_z}(\phi)_{\B(\Zpmc)}.$$
        \end{proposition}

        For the proof we refer the reader to~\cite[Lemma~4.2]{AGW14}. The main idea is to use $\alpha$-$\beta$-bordered Heegaard diagrams introduced in~\cite{LOT-mor}. 

        \begin{corollary}\label{cor}
        Hochschild homologies of the two bimodules are isomorphic:
        $$HH_*({}^{\B(\Zpmc)}N(\phi)_{\B(\Zpmc)})\cong HH_*( {}_{\A(\F_z(\Sigma))}N_{\F_z}(\phi)_{\A(\F_z(\Sigma))} ).$$
        \end{corollary}
        This follows from Proposition~\ref{prop:bims_are_eq} and~\cite[Proposition~2.3.54]{LOT-bim}.

\Needspace{8\baselineskip}
\section{Theoretical evidence for the conjecture}\label{sec:LF}
    Let us step back and see what we have accomplished thus far. To a surface associated to a pointed matched circle, $(\Sigma, \partial \Sigma )=F^\circ(\Zpmc)$, we associated two quasi-isomorphic algebras:
    $$
    \begin{tikzcd}[column sep=15pt,row sep=45pt]
        & (\Sigma, \partial \Sigma )=F^\circ(\Zpmc) \arrow[dr, rightsquigarrow,  "\text{Section~\ref{sec:hf_bim}}" right ]  \arrow[dl, rightsquigarrow, "\text{Section~\ref{aur_constr}}" left]&   
        \\
         \A(\F_z(\Sigma)) = \bigoplus\limits_{1\leq i,j\leq 2g} hom_{\F_z(\Sigma)} (\alpha_i,\alpha_j)&&\B(\Zpmc) \arrow[ll, leftrightarrow, "\simeq" above, "\text{Theorem~\ref{thm:aur}}" below] 
    \end{tikzcd}
    $$
    To the mapping class $\phi \lefttorightarrow (\Sigma, \partial \Sigma )=F^\circ(\Zpmc)$ we associated two homotopy equivalent bimodules, and a version of fixed point Floer cohomology. We conjectured that Hochschild homology of either of the bimodules is isomorphic to fixed point Floer cohomology, and supported Conjecture~\ref{conj} by computations in Section~\ref{supporting_comps}.
    $$
    \begin{tikzcd}[row sep=50pt,column sep=15pt]
        &  
        \text{Mapping class }\phi \lefttorightarrow (\Sigma, \partial \Sigma )=F^\circ(\Zpmc) 
        \arrow[dl, rightsquigarrow, "\text{Section~\ref{alt_bim}}" left]
        \arrow[d, rightsquigarrow, "\text{Section~\ref{sec:bimodule from bordered}}" left]
        \arrow[ddr, rightsquigarrow, "\text{Section~\ref{sec:fixed point Floer}}"]
        &
        \\
        _{\A(\F_z(\Sigma))}N_{\F_z}(\phi)_{\A(\F_z(\Sigma))}
        \arrow[d, rightsquigarrow, "\text{Section~\ref{sec:hoch}}" right] 
        &
        ^{\B(\Zpmc)}N(\phi)_{\B(\Zpmc)}
        \arrow[l, leftrightarrow, "\simeq" above, "\text{Prop.~\ref{prop:bims_are_eq}}" below] 
        \arrow[d, rightsquigarrow, "\text{Section~\ref{sec:hoch}}" left] 
        &
        \\
        HH_*(N_{\F_z}(\phi))
        \arrow[rr, bend right=25, "\text{Section~\ref{subsec:O-C map}}" above, "\text{(in the double basepoint case)}" below]
        & 
        HH_*(N(\phi))
        \arrow[l, leftrightarrow, "\cong" above, "\text{Corollary~\ref{cor}}" below]
        \arrow[r, leftrightarrow, "\overset{?}{\cong}" above, "\text{Conjecture~\ref{conj}}" below] 
        & 
        HF^*(\tilde{\phi}\:;\ U_2^+,U_1^-)
    \end{tikzcd}
    $$
    In contrast to Section~\ref{supporting_comps}, where we performed concrete computations supporting Conjecture~\ref{conj}, in this section we describe some more theoretical evidence. We are going to construct the map corresponding to the lowest arrow in the diagram above. For that we will need mild generalizations of all the invariants, namely we will need to use two basepoints instead of one in all the constructions. In Section~\ref{1bp_to_2bp_bim} we describe the known two-basepoints-generalizations $N_{\F_{\{z_1,z_2\}}}(\phi)$ and $N^\text{2bp}(\phi)$ of bimodules $N_{\F_z}(\phi)$ and $N(\phi)$. In Section~\ref{1bp_to_2bp_hf} we first describe some heuristics, which explain why $HF^*(\tilde{\phi}\:;\:U_2^+,U_1^-)$ corresponds to the one basepoint case. Based on this we introduce the double basepoint version $HF^*(\dtilde{\phi}\:;\ U_2^+,U_3^+, U_1^-)$. After that we state a double basepoint version of Conjecture~\ref{conj}:
    $$HH_*(N^\text{2bp}(\phi)) \ \overset{?}{\cong} \ HF^*(\dtilde{\phi}\:;\:  U_2^+,U_3^+, U_1^-).$$ 
    In Section~\ref{LF} we cover the background material on how the Lefschetz fibration structure on a surface gives rise to a special type of the Fukaya category, the Fukaya-Seidel category. Finally, in Section~\ref{subsec:O-C map} we describe the so called open-closed map, which results in a map $HH_*(N^\text{2bp}(\phi)) \ \rightarrow \ HF^*(\dtilde{\phi}\:;\:U_2^+,U_3^+, U_1^-)$. This map is widely believed to be an isomorphism, see ~\cite[Conjecture 7.18]{Sei3II}, and if true, it would prove the double basepoint version of Conjecture~\ref{conj}.

    \Needspace{8\baselineskip}
    \subsection{From one basepoint to two: bimodule}\label{1bp_to_2bp_bim}
\begin{figure}[ht]
        \centering
            \includegraphics[width=0.5\textwidth]{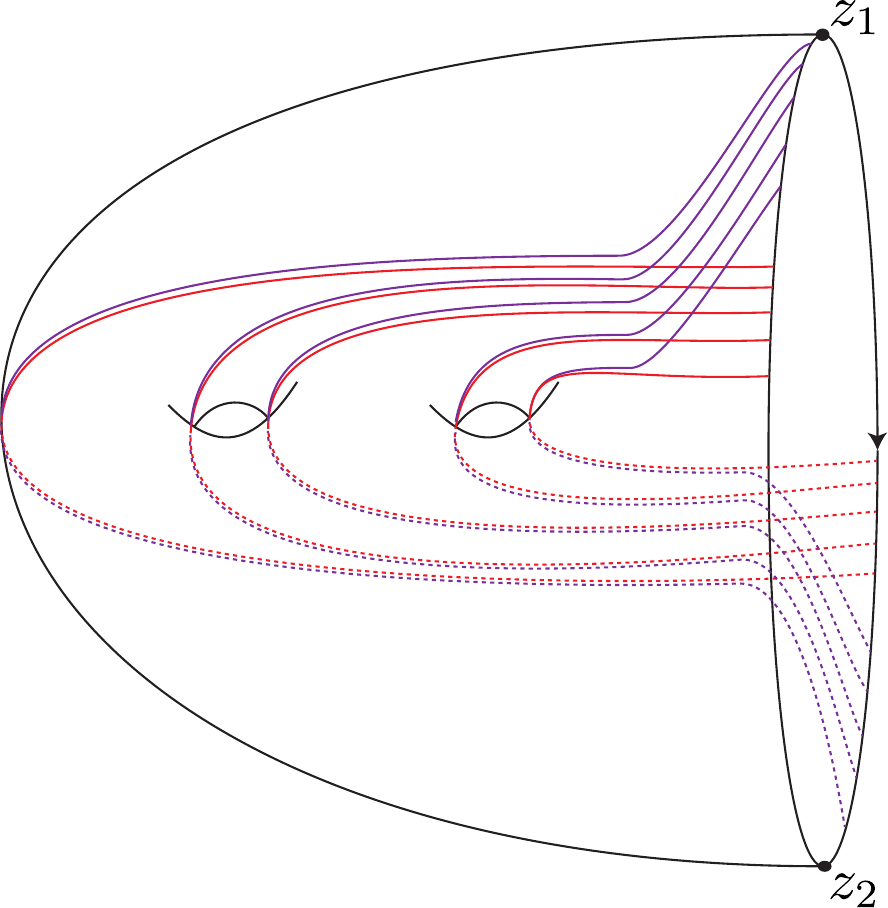}
            \caption{Generators of the partially wrapped Fukaya category $\F_{\{z_1,z_2\}}(\Sigma_2)$.}
            \label{fig:partially_wrapped_2bpts}
        \end{figure}

        Let us explain how to modify our previous constructions of bimodules, if we want to have two basepoints instead of one. First of all, in Section~\ref{aur_constr} the partially wrapped Fukaya category was defined for any number of basepoints. In Figure~\ref{fig:partially_wrapped_2bpts} we draw a generating set of Lagrangians (red curves) in the case of two basepoints on the genus two surface, together with their perturbations (purple curves).  Now we have five Lagrangian arcs as generators, instead of four in the one basepoint case (Figure~\ref{fig:pointed_matched_circle}). In general, for genus $g$ surface $(\Sigma,\partial \Sigma)$ the number of generating arcs will be $2g+1$ and $2g$ for two and one basepoint cases, respectively. We denote the double basepoint partially wrapped Fukaya category by $\F_{\{z_1,z_2\}}(\Sigma)$

\begin{figure}[ht]
        \centering
            \includegraphics[width=0.6\textwidth]{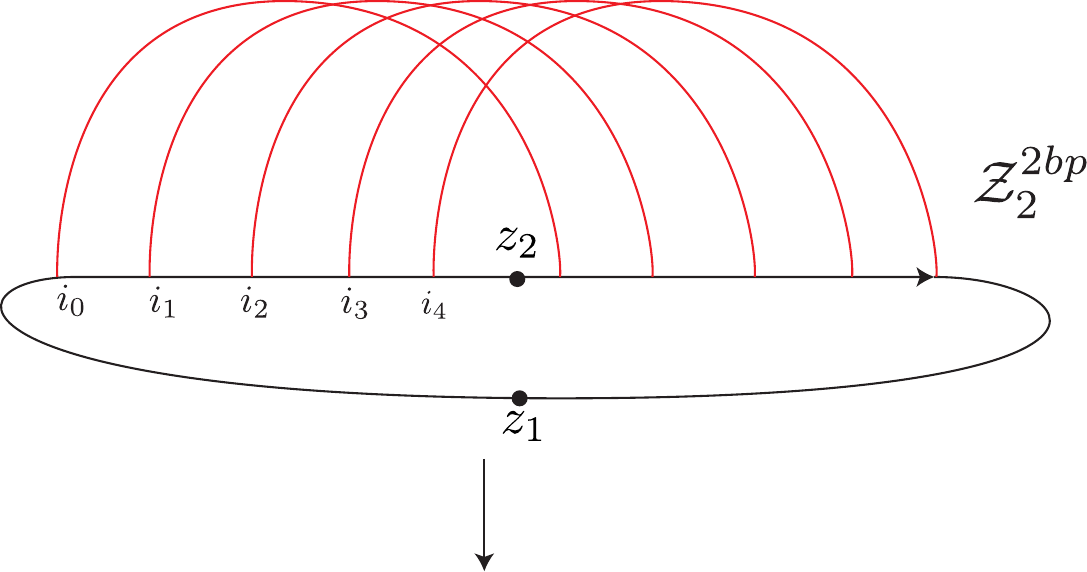}
            \tikzstyle{arrow} = [thick,->,>=stealth]
            \begin{tikzpicture}[scale=1,baseline=1.5cm]
            \node[circle](i_0) at (0,0){$i_0$};
            \node[circle](i_1) at (2,0){$i_1$};
            \node[circle](i_2) at (4,0){$i_2$};
            \node[circle](i_3) at (6,0){$i_3$};
            \node[circle](i_4) at (8,0){$i_4$};

            \draw[arrow, color=green] (i_0) to[bend left] (i_1);
            \draw[arrow, color=blue] (i_0) to[bend right] (i_1);
            \draw[arrow, color=green] (i_1) to[bend left] (i_2);
            \draw[arrow, color=blue] (i_1) to[bend right] (i_2);
            \draw[arrow, color=green] (i_2) to[bend left] (i_3);
            \draw[arrow, color=blue] (i_2) to[bend right] (i_3);
            \draw[arrow, color=green] (i_3) to[bend left] (i_4);
            \draw[arrow, color=blue] (i_3) to[bend right] (i_4);
            
            \draw (0,1) to[bend right] node[left,text width=3cm, text centered]{Path \\ algebra over $\mathbb{F}_2$ \\ (composing only the same color paths)} (0,-1);
            \draw (8,1) to[bend left] (8,-1);
            \node[scale=1.2] at (-4.5,0){$\B(\Zpmc^\text{2bp}_2)=$};
            \end{tikzpicture}
            \caption{A genus two double basepoint example of how to get a $dg$ algebra out of a pointed matched circle. Paths consisting of different color arrows are prohibited.}
            \label{fig:Z_2-2bp--B-Z_2-2bp-}
        \end{figure}

        For the corresponding double basepoint version of pointed matched circle, as well as the corresponding algebra, see Figure~\ref{fig:Z_2-2bp--B-Z_2-2bp-}.

        Just as in Section~\ref{alt_bim}, we can define a graph bimodule for a mapping class $\phi \in MCG_0(\Sigma, \partial \Sigma)$ via the double basepoint partially wrapped Fukaya category $\F_{\{z_1,z_2\}}(\Sigma_2)$: 
        $$_{\A(\F_{\{z_1,z_2\}}(\Sigma))}N_{\F_{\{z_1,z_2\}}}(\phi)_{\A(\F_{\{z_1,z_2\}}(\Sigma))}=\bigoplus\limits_{1\leq i,j\leq k} hom_{\F_{\{z_1,z_2\}}(\Sigma)} (\alpha_i,\phi(\alpha_j)).$$ 
        Auroux's Theorem~\ref{thm:aur} works for any number of basepoints, and so if $(\Sigma,\partial \Sigma)=F^\circ(\Zpmc^\text{2bp})$, then $\A(\F_{\{z_1,z_2\}}(\Sigma)) \simeq \B(\Zpmc^\text{2bp})$. 
        
\begin{figure}[ht]
        \centering
        \includegraphics[width=0.5\textwidth]{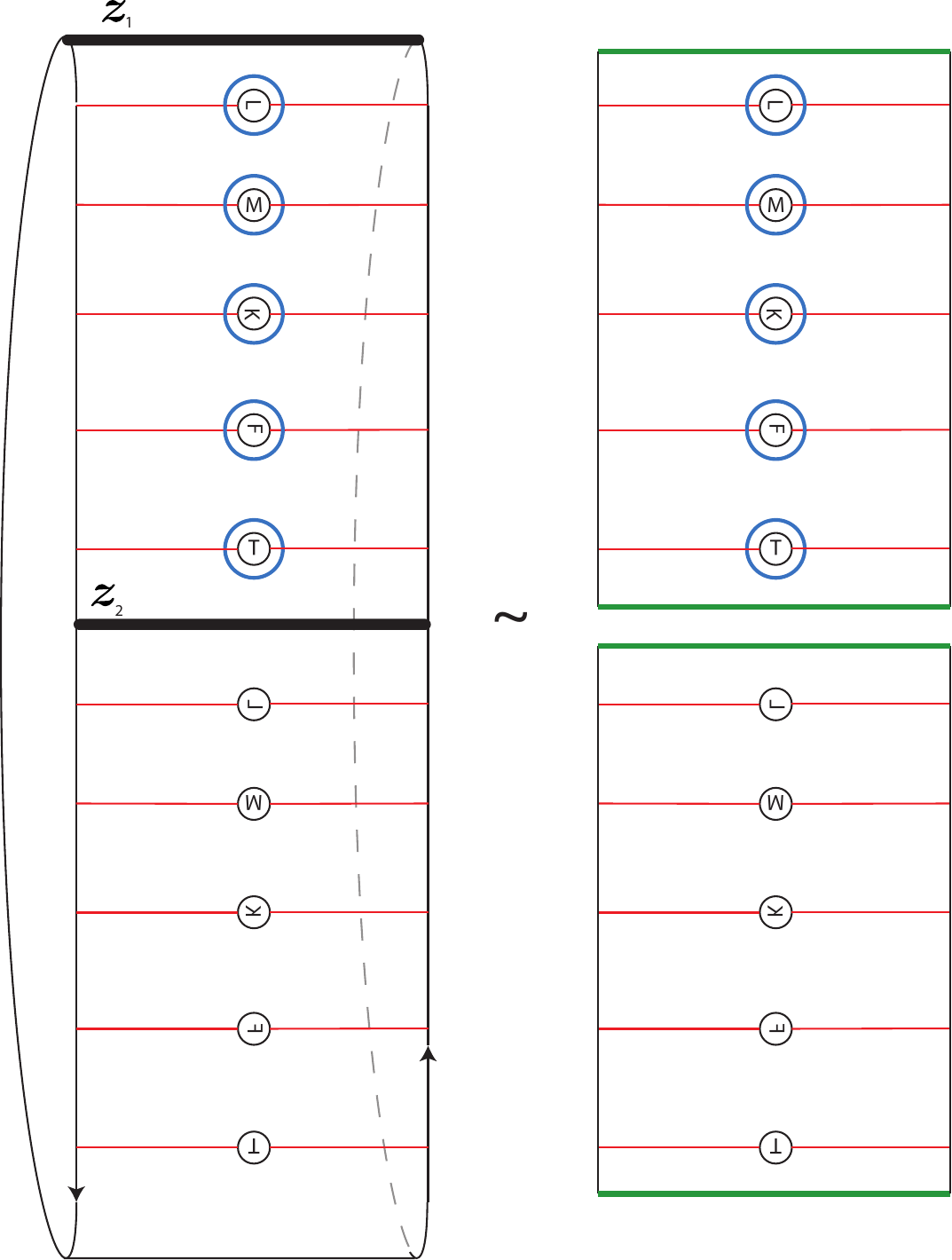}
        \caption{Double basepoint diagrams for mapping cylinder of $\id\co\Sigma_2 \rightarrow \Sigma_2$.}
        \label{fig:heegard_diagram_id_2bpts}
        \end{figure}

        Despite the fact that we constructed $_{\B(\Zpmc^\text{2bp})}N_{\F_{\{z_1,z_2\}}}(\phi)_{\B(\Zpmc^\text{2bp})}$, for computations we would prefer to have a type $DA$ bimodule $^{\B(\Zpmc^\text{2bp})}N^\text{2bp}(\phi)_{\B(\Zpmc^\text{2bp})}$, generalizing $N(\phi)=N^\text{1bp}(\phi)$. Such a bimodule, as in the one basepoint case, comes from a Heegaard diagram for the mapping cylinder, but equipped with two basepoints on each boundary, and two arcs connecting them. The necessary machinery of bordered Heegaard diagrams with multiple basepoints was invented by Zarev in~\cite{Zarev}; in Figure~\ref{fig:heegard_diagram_id_2bpts} we draw two diagrams for the identity mapping class in the genus two case: on the right there is Zarev's bordered sutured Heegaard diagram, and on the left we drew a double basepoint Heegaard diagram, which would be a natural generalization of the one basepoint diagram. The two diagrams carry the same holomorphic information, and the bimodules coming from them are the same. The reason why they are different is because Zarev used the language of sutured manifolds and their bordered versions. In order to go from the left diagram to the right, instead of drawing basepoints and basepoint arcs, one deletes their neighborhoods and then sets the boundary coming from these neighborhoods (drawn in green) to be forbidden for holomorphic discs.

        The generalization of Figure~\ref{fig:heegard_diagram_id_2bpts} to non-trivial mapping classes of genus $g$ surfaces is completely analogous to the one basepoint case. As in the one basepoint case, Proposition~\ref{prop:bims_are_eq} holds true:
        $$ {}_{\B(\Zpmc^\text{2bp})}{\B(\Zpmc^\text{2bp})}_{\B(\Zpmc^\text{2bp})} \boxtimes {}^{\B(\Zpmc^\text{2bp})}N^\text{2bp}(\phi)_{\B(\Zpmc^\text{2bp})} \simeq {}_{\B(\Zpmc^\text{2bp})}N_{\F_{\{z_1,z_2\}}}(\phi)_{\B(\Zpmc^\text{2bp})}.$$

        Let us now discuss the relationship between one and two basepoint cases. It turns out that the Hochschild homologies of one and two basepoint bimodules are related, namely
        \begin{equation}\label{HH_for_1bp_and_2bp}
        \rk(HH_*(N^\text{2bp}(\phi))=\rk(HH_*(N(\phi))+1.
        \end{equation} 
        The reason is that Hochschild homology of the double basepoint bimodule is equal to knot Floer homology of the binding of the corresponding open book in the second lowest Alexander grading, where knot in a Heegaard diagram is specified by four basepoints, instead of two: $HH_*(N^\text{2bp}(\phi))=\widehat{HFK}^\text{4bp}(M_\phi^\circ,K \:;\:  -g)$.  And the difference between the four basepoint and the usual two basepoint knot Floer homologies is known: $\widehat{HFK}^\text{4bp}(M_\phi^\circ,K)=(\mathbb{F}_2)^2 \o \widehat{HFK}(M_\phi^\circ,K)$, where the Alexander gradings of two generators of $(\mathbb{F}_2)^2$ are 0 and -1. Thus we have
        \begin{align*} 
        &\rk(HH_*(N^\text{2bp}(\phi))=\rk(\widehat{HFK}^\text{4bp}(M_\phi^\circ,K \:;\:  -g))= \\ 
        = \ &\rk \big((\mathbb{F}_2)_{(\text{A}=-1)} \o \widehat{HFK}(M_\phi^\circ,K \:;\:  -g+1) \oplus (\mathbb{F}_2)_{(\text{A}=0)} \o \widehat{HFK}(M_\phi^\circ,K \:;\:  -g) \big)= \\
        = \ &\rk \big( \widehat{HFK}(M_\phi^\circ,K \:;\:  -g+1) \oplus \widehat{HFK}(M_\phi^\circ,K \:;\:  -g) \big)=\\
        = \ & \rk(HH_*(N(\phi)) + 1,
        \end{align*}
        because the lowest $-g$ Alexander grading of knot Floer homology of a fibered knot (i.e. the binding of an open book) is always one.
    \Needspace{8\baselineskip}
    \subsection{From one basepoint to two: fixed point Floer cohomology} \label{sec:HF_2bpts}\label{1bp_to_2bp_hf}
\begin{figure}[H]
        \centering
            \includegraphics[width=1\textwidth]{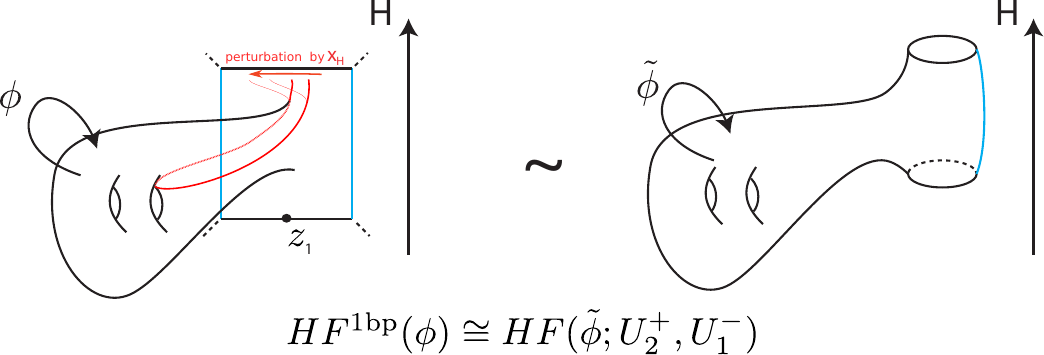}
            \caption{Left: the behavior of Hamiltonian perturbation one needs to consider in the one basepoint version $HF^\text{1bp}(\phi)$. Right: an equivalent to it theory $HF^{*}(\tilde{\phi}\:;\:  U_2^+,U_1^-)$ on the right, where the Hamiltonian is constant on the boundaries.}
            \label{fig:1bpt_perturbation}
        \end{figure}

        Let us first explain why the choice of fixed point Floer cohomology $HF^*(\tilde{\phi}\:;\:U_2^+,U_1^-)$  corresponds one basepoint. The point is that there is a version of fixed point Floer cohomology $HF^\text{1bp}(\phi)$, which is equal to $HF^{*}(\tilde{\phi}\:;\:  U_2^+,U_1^-)$, but defined without deleting a second disc from the surface. In Section~\ref{sec:fixed point Floer}, we decided not to give a rigorous definition of $HF^\text{1bp}(\phi)$,
        % (which would be analogous to~\cite[Section~6]{Sei3II})
        and chose to use existing methods and work with $HF^{*}(\tilde{\phi}\:;\:  U_2^+,U_1^-)$.
        % see~\cite[Section~6]{Sei3II} and~\cite[Proposition~5.12]{Sei3I} for compactness argument both in vertical and horizontal directions. 
        But, because we would like to have the double basepoint version, we now indicate the setup for $HF^\text{1bp}(\phi)$. It can be defined only for infinite area surfaces with a cylindrical end, rather than compact surfaces with boundary, and so we have to work with the induced compactly supported exact self-diffeomorphisms $\phi$ on the completion $\hat{\Sigma}$. Analogously to how we chose the specific Hamiltonian perturbations near $\partial \Sigma$ in the definition of $HF^{*}(\tilde{\phi}\:;\:  U_2^+,U_1^-)$, we now need to specify behavior of Hamiltonian perturbation near $\infty$ on $\hat{\Sigma}$. Inspired by~\cite[Section~6]{Sei3II}, we indicate the behavior of the Hamiltonian on $\infty$ in Figure~\ref{fig:1bpt_perturbation}, the left side. Upwards and downwards the Hamiltonian is linear with respect to the radial coordinate. Comparing the left and the the right side of the figure (on the right we glued the blue boundaries together), we can see why such Hamiltonian perturbation is equivalent to the one considered in the definition of $HF^{*}(\tilde{\phi}\:;\:  U_2^+,U_1^-)$ --- the generators (fixed points) and the differentials in Floer cohomology on the left side and on the right side are in $1-1$ correspondence.

        The reason why we call $HF^\text{1bp}(\phi)$ the one basepoint version is because the Hamiltonian on the left side of Figure~\ref{fig:1bpt_perturbation} could be used in the definition of the partially wrapped Fukaya category with one basepoint. As indicated by the orange curve and its perturbation in Figure~\ref{fig:1bpt_perturbation}, putting the basepoint in the bottom part of infinity, and allowing Lagrangian arcs to go only to the upper part of infinity, one obtains perturbations sending all the arcs to the left, i.e. to the basepoint\footnote{Technically, one has to take the limit of these perturbations, because the Hamiltonian is linear with respect to the radial coordinate, as opposed to quadratic one used in~\cite{Aur2}.}.
        % \begin{remark}
        % When defining perturbations in the partially wrapped Fukaya category, we have to make sure every Lagrangian arc will be wrapped enough to intersect the other ones to the left of it at infinity (see Figure~\ref{fig:bordered=partially_wrapped}). Thus we actually cannot consider a linear Hamiltonian with respect to radial coordinate $H=r$ on the upper half of infinity, as it is pictured in Figure~\ref{fig:1bpt_perturbation}. Instead, in order to make sure that Hamiltonian has big enough derivative for every pair of Lagrangian arcs, we need to take a limit $H=\delta r, \delta \rightarrow +\infty$, as it was done in~\cite[Definition 4.1]{Aur1}. But because we always consider only finite number of generating arcs going to infinity, it is enough to just consider linear Hamiltonian with big enough derivative $H=\delta r, \delta \gg 0$.
        % % one can also try to make this hamiltonian quadratic
        % \end{remark}

\begin{figure}[H]
        \centering
            \includegraphics[width=1\textwidth]{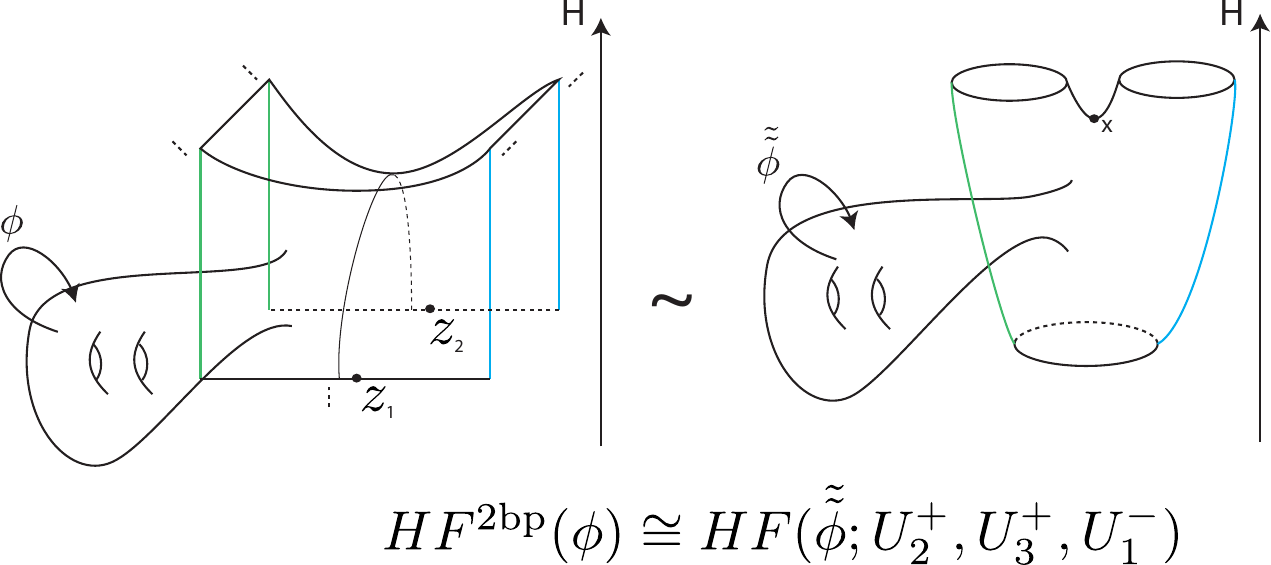}
            \caption{Left: the behavior of Hamiltonian perturbationone needs to consider in the double basepoint version $HF^\text{2bp}(\phi)$ (see~\cite[Section~6]{Sei3II} for a rigorous setup). Right: an equivalent to it theory $HF^*(\tilde{\tilde{\phi}}\:;\:  U_2^+,U_3^+,U_1^-)$, where Hamiltonian is constant on boundaries.}
            \label{fig:2bpt_perturbation}
        \end{figure}

        Now we consider the double basepoint counterpart $HF^\text{2bp}(\phi)$ of the above construction. The corresponding behavior of Hamiltonian near $\infty$ on $\hat{\Sigma}$ is pictured on the left of Figure~\ref{fig:2bpt_perturbation}. The cohomology theory  $HF^\text{2bp}(\phi)$ was developed in~\cite[Section~6]{Sei3II}, viewing $\hat{\Sigma}$ as double branched cover of $\C$, or equivalently, as a total space of a $0$-dimensional Lefschetz fibration over $\C$. A different, but equivalent version of Floer cohomology (where Hamiltonian is constant on boundary components) is depicted on the right --- instead of one disc we need to take out two discs this time, and we denote the resulting Floer cohomology by $HF^*(\tilde{\tilde{\phi}}\:;\:  U_2^+,U_3^+,U_1^-)$. 
     
        \begin{remark}
        The one and two basepoint versions of fixed point Floer cohomology should be related, and for all the cases we considered, as in the case of rank relationship equation~\ref{HH_for_1bp_and_2bp} between Hochschild homologies, we had $\rk(HF^*(\tilde{\tilde{\phi}}\:;\:  U_2^+,U_3^+,U_1^-))=\rk(HF^*(\tilde{\phi}\:;\:  U_2^+,U_1^-))+1.$ We did not find a general explanation for this. The reason might be that if one compares cochain complexes, then they are identical except $CF^*(\tilde{\tilde{\phi}}\:;\:  U_2^+,U_3^+,U_1^-)$ has one more generator $x$, depicted on the right of Figure~\ref{fig:2bpt_perturbation}.  This generator does not have any differentials going out of it, as they are gradient lines going up from $x$ for suitable Hamiltonian, and there are no generators above $x$. It is likely that it also does not have any differentials going in (or, rather, one can arrange the Hamiltonian in such a way).
        \end{remark}
        % which can be possibly proved via neck stretching argument analogous to one in paper~\cite{Sei1}. The stretching should be along such a circle, that energy bounds are constants, mb. See Gautschi_PhD, Propositino 2.10.

        Now we are ready to state a double basepoint version of Conjecture~\ref{conj}. Namely, the following should be true:
        \begin{conjecture}\label{conj_2bp}
        For every mapping class $\phi \in MCG_0(\Sigma,\partial\Sigma=S^1=U_1)$ there is an isomorphism of $\mathbb{Z}_2$-graded vector spaces
        $$HH_*(N^{\text{2bp}} (\phi^{-1})) \cong HF^{*+1}(\tilde{\tilde{\phi}}\:;\:  U_2^+,U_3^+,U_1^-).$$
        \end{conjecture}
        We did lots of computations using~\cite{Pyt}, and all of them support this conjecture. These computations are completely analogous to the ones in Section~\ref{supporting_comps}, and so we choose not  to describe them. Rather, below we describe a more theoretical evidence for Conjecture~\ref{conj_2bp}.
    \Needspace{8\baselineskip}
    \subsection{\texorpdfstring{$\Sigma$ as Lefschetz fibration and the Fukaya-Seidel category}{Σ as Lefschetz fibration and the Fukaya-Seidel category}}\label{LF}
        Take an area preserving double branched cover $f\co\hat{\Sigma}\rightarrow \C$ of an exact surface $\hat{\Sigma}$ with cylindrical end over the complex numbers. For instance, one may take a quotient by the hyperelliptic involution, which we drew below in Figure~\ref{fig:LF}. We can view this cover as an exact symplectic fibration with singularities, as in~\cite[Setup 5.1]{Sei3II}. This fibration is in fact a $0$-dimensional Lefschetz fibration, with 2g+1 critical points.  
        % Being a Lefschetz fibration corresponds to critical points not having more than order two branching (i.e. triple and more branch covers also can be Lefschetz fibrations). 
        We assume that the critical values $p_1,\dots,p_{2g+1}$ all satisfy $\Re(p_i)=0$ and $\Im(p_1)<\dots<\Im(p_{2g+1})$. The genus two case satisfying these properties is drawn in Figure~\ref{fig:LF}.

        \begin{figure}[H]
        \centering
            \includegraphics[width=0.5\textwidth]{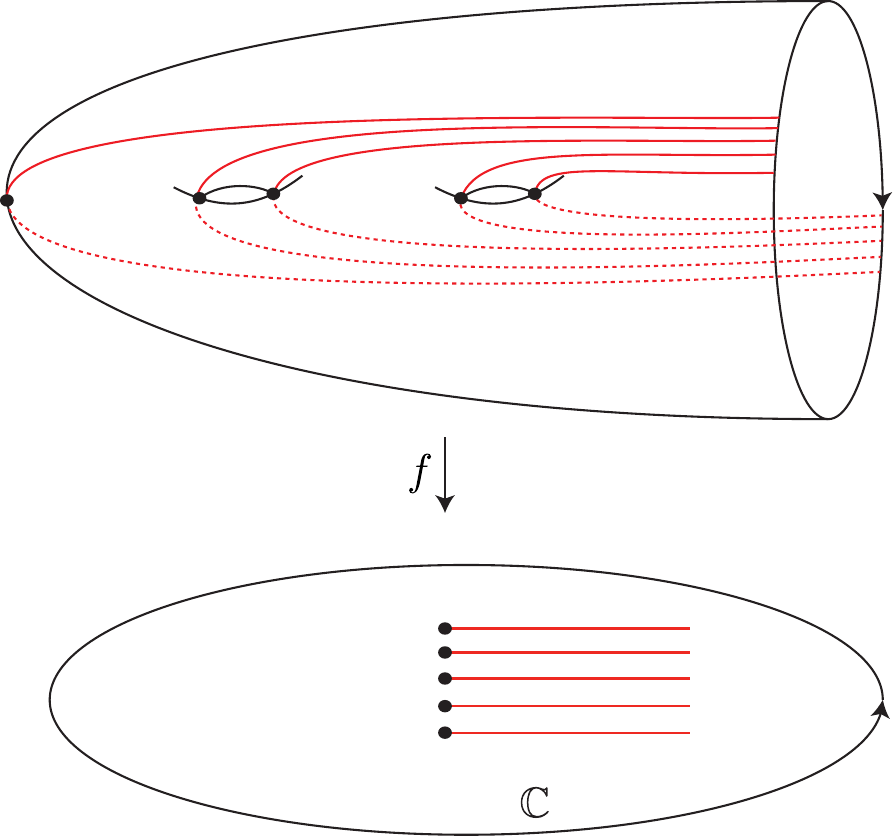}
            \caption{0-dimensional Lefschetz fibration structure on the genus two surface.}
            \label{fig:LF}
        \end{figure}

        It turns out that the map $f\co\hat{\Sigma}\rightarrow \C$ can produce a specific version of fixed point Floer cohomology of $\phi \lefttorightarrow \Sigma$, and a specific version of the Fukaya category of $\Sigma$. It turns out that those versions are equal to $HF^\text{2bp}(\phi)$ and $\F_{\{z_1,z_2\}}(\Sigma)$ respectively. We discuss this below.

        Following~\cite[Section~6]{Sei3II}, having an exact compactly supported self-diffeomorphism $\phi\co \hat{\Sigma} \rightarrow \hat{\Sigma}$ of $0$-dimensional Lefschetz fibration, we can consider fixed point Floer cohomology $HF^*(\phi,\delta > 0,\epsilon)$, where $\epsilon$ is not important to us because fibers of $f\co\hat{\Sigma}\rightarrow \C$ do not have boundary, and $\delta$ is responsible for  perturbation at infinity by the Hamiltonian 
        \begin{equation}\label{hamiltonian_2bp}
        H:\hat{\Sigma}\rightarrow \mathbb R, \ x \mapsto \delta \Re(f(x)).
        \end{equation} 
        This theory depends only on the sign of $\delta$, and from the definitions of Hamiltonians it follows that 
        $$HF^\text{2bp}(\phi)=HF^*(\phi,\delta > 0).$$
        
        The Lefschetz fibration structure over the complex plane can be also used to define a special type of $A_\infty$ category $\F_f(\hat{\Sigma})$, which is called the Fukaya-Seidel category; see~\cite{Sei4,Sei5}, and more recent articles~\cite{Sei3I,Sei3II} for a more relevant setup for us. In our case of the double branched cover $f\co\hat{\Sigma}\rightarrow \C$, the objects of the category are compact exact Lagrangians in $\hat{\Sigma}$, and also non-compact ones which are Lefschetz thimbles associated to admissible arcs\footnote{Admissible arcs in $\C$ are proper rays which start at a critical value of $f$, do not pass over other critical values, and at some point stabilize to be horizontal, oriented to the right rays. In our case of $f\co\hat{\Sigma}\rightarrow \C$ Lefschetz thimbles are just preimages of admissible arcs.} in $\C$. Perturbation at infinity is defined using the same Hamiltonian~\ref{hamiltonian_2bp}, but making sure that $\delta \gg 0$ is big enough. 
        % It turns out that it is exactly the same type of Hamiltonian perturbation as for $\F_{\{z_1,z_2\}}(\Sigma)$, see the left side of Figure~\ref{fig:2bpt_perturbation}. 

        % Fukaya-Seidel category $\F_f(\hat{\Sigma})$ is quasi-equivalent to partially wrapped Fukaya category $\F_{\{z_1,z_2\}}(\Sigma)$.
        It turns out that the Fukaya-Seidel category $\F_f(\hat{\Sigma})$ is quasi-equivalent to the partially wrapped Fukaya category $\F_{\{z_1,z_2\}}(\Sigma)$, despite the fact that the non-compact objects allowed are different. For a setup, which mediates between the partially wrapped category with two basepoints and the Fukaya-Seidel category, see~\cite[Section~3.2]{Aur1}.
        % (where the Lefschetz fibration structure is used to define $\F_{\{z_1,z_2\}}(\Sigma)$) 
        % For a generalization of that setup to the partially wrapped Fukaya category (which does not use fibration structure) see~\cite[Section~4.1]{Aur1}. 

        It was proved in~\cite{Sei5} that Lefschetz thimbles (one for each critical points) generate the Fukaya-Seidel category. The following are several consequences of the quasi-equivalence $\F_f(\hat{\Sigma}) \simeq \F_{\{z_1,z_2\}}(\Sigma)$. If we choose a generating set of thimbles for the category $\F_f(\hat{\Sigma})$, then these thimbles also generate $\F_{\{z_1,z_2\}}(\Sigma)$. An example of a generating set of Lefschetz thimbles in the genus two case is drawn in Figure~\ref{fig:LF}, and the same set of generators for $\F_{\{z_1,z_2\}}(\Sigma)$ was drawn in Figure~\ref{fig:partially_wrapped_2bpts}.
        The $hom$-algebras for the two categories are also the same: 
        $$\bigoplus\limits_{1\leq i,j\leq k} hom_{\F_f(\hat{\Sigma})} (\alpha_i,\alpha_j) \simeq \bigoplus\limits_{1\leq i,j\leq k} hom_{\F_{\{z_1,z_2\}}(\Sigma)} (\alpha_i,\alpha_j) \simeq \B (\Zpmc^\text{2bp}),$$ 
        and the bimodules corresponding to exact automorphisms 
        $\phi\co \hat{\Sigma} \rightarrow \hat{\Sigma}$ are also the same: 
        $$N_{\F_f}(\phi)\coloneq\bigoplus\limits_{1\leq i,j\leq k} hom_{\F_f(\hat{\Sigma})} (\alpha_i,\phi(\alpha_j)) \simeq \bigoplus\limits_{1\leq i,j\leq k} hom_{\F_{\{z_1,z_2\}}(\Sigma)} (\alpha_i,\phi(\alpha_j)) \simeq N^\text{2bp}(\phi).$$
    \Needspace{8\baselineskip}
    \subsection{Open-closed map}\label{subsec:O-C map}
        The open-closed map is a map between Hochschild homology of a graph bimodule and a fixed point Floer cohomology $OC:HH_*(N_{\text{Fuk}}(\phi)) \rightarrow HF(\phi)$. The key point is that this map can be defined only if Hamiltonian perturbations near $\infty$ on $\hat{\Sigma}$  are the same for $\text{Fuk}$ and $HF(\phi)$. Because Hamiltonian perturbations for $\F_f(\hat{\Sigma})$ and $HF^\text{2bp}(\phi)$ are the same, one expects to have a map 
        $$OC:HH_*(N_{\F_f}(\phi)) \rightarrow HF^\text{2bp}(\phi),$$
        which would define a map in one direction of Conjecture~\ref{conj_2bp} (remember that $HH_*(N_{\F_f}(\phi)) \cong HH_*(N^\text{2bp}(\phi))$ and $HF^\text{2bp}(\phi) \cong HF^{*+1}(\tilde{\tilde{\phi}}\:;\:  U_2^+,U_3^+,U_1^-)$).
        Indeed, such a map was constructed by Seidel in~\cite{Sei3II}, and we describe the construction details below.

\begin{figure}[h]
        \centering
            \includegraphics[width=0.4\textwidth]{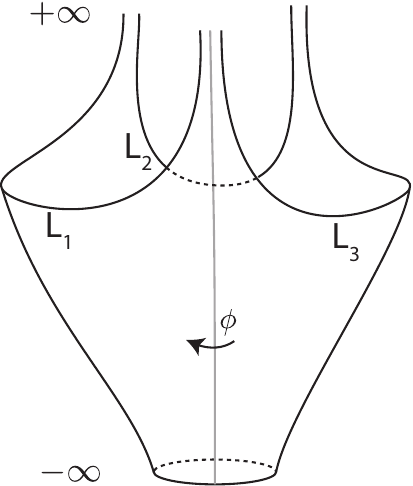}
            \caption{The open-closed map counts such holomorphic objects inside $\hat{\Sigma}$. Compare this to~\cite[Figure 3]{Sei3II}.}
            \label{fig:o-c}
        \end{figure}

        In~\cite[Section~7]{Sei3II} Seidel constructed open-closed map in the case where the symplectic manifold is an exact symplectic fibration with singularities over $\C$, which includes Lefschetz fibrations, and in particular double branched covers $f\co\hat{\Sigma} \rightarrow \C$. In our case of $f\co\hat{\Sigma} \rightarrow \C$ the open-closed map counts isolated points in the moduli space of holomorphic maps from a Riemann surface drawn in Figure~\ref{fig:o-c} to $\hat{\Sigma}$, with a twist $\phi$ along the gray line (compare with~\cite[Figure 3]{Sei3II}). These maps have the following boundary conditions: a twisted orbit of Hamiltonian vector field $X_H$ on one end, which is equivalent to a constant section of the mapping torus $T_{\phi \circ \psi^1_{X_{H}}}$, and a chain of Lagrangians on the other, with consistent perturbations. Along the gray line the map has a twist $\phi$. So the strip end with the gray line limits to an intersection point of $\phi \circ \psi^1_{X_{H}}(L_3) \cap L_1$ to the left of the gray line, and to the intersection point $L_3 \cap (\phi \circ \psi^1_{X_{H}})^{-1}(L_1)$ to the right of the gray line.

        % In our case $\phi\co \hat{\Sigma}\rightarrow \hat{\Sigma}$ is an exact symplectomorphism, which outside of the compact set is a Hamiltonian vector field flow, where Hamiltonian is pictured on the left of figures~\ref{fig:1bpt_perturbation} and~\ref{fig:2bpt_perturbation}, depending on the version we are working in. Hamiltonian needs to have big enough derivative. 

        % Hamiltonian for perturbations should be chosen to be the same as in definition of the bimodule and fixed point Floer cohomology.  

        In this setting Seidel in~\cite[Equation 7.15]{Sei3II} defines a bimodule $\mathcal P(\phi,\delta, \epsilon)=\bigoplus\limits_{1\leq i,j\leq k} hom (\psi^1_{X_H}(\phi(\alpha_i)),\alpha_j)$. The $\epsilon$ does not play any role for us, because in our case of $0$-dimensional Lefschetz fibration $f\co\hat{\Sigma} \rightarrow \C$ there is no boundary in a fiber. The $\delta$ is responsible for the Hamiltonian $H(x)=\delta \Re(f(x))$ which is used to perturb $\phi$ at $\infty$ of $\hat{\Sigma}$.
        If we assume $\delta \gg 0$ so that the generating Lagrangians are wrapped enough to intersect each other at $\infty$, 
        then this bimodule is 
        \begin{align*}
        &\mathcal P(\phi,\delta \gg 0)=\bigoplus\limits_{1\leq i,j\leq k} hom (\psi^1_{X_H}(\phi(\alpha_i)),\alpha_j)= \bigoplus\limits_{1\leq i,j\leq k} hom_{\F_f(\hat{\Sigma})} (\phi(\alpha_i),\alpha_j)=\\
        &=\bigoplus\limits_{1\leq i,j\leq k} hom_{\F_{\{z_1,z_2\}}(\Sigma)} (\phi(\alpha_i),\alpha_j)= 
        \bigoplus\limits_{1\leq i,j\leq k} hom_{\F_{\{z_1,z_2\}}(\Sigma)} (\alpha_i,\phi^{-1}(\alpha_j))=N^\text{2bp}(\phi^{-1}).
        \end{align*}

        Now we turn our attention to~\cite[Conjecture 7.18]{Sei3II}. In our case this amounts to the following --- there is an open-closed map which gives an isomorphism:
        $$
        OC\co HH_*(N^\text{2bp}(\phi^{-1})) \  \xrightarrow{\overset{?}{\cong}} \ HF_\text{2bp}^{*+1}(\phi).
        $$ 

        As a consequence, our double basepoint Conjecture~\ref{conj_2bp} is a special case of Seidel's conjecture. The one basepoint version, i.e. Conjecture~\ref{conj}, most likely fits in a similar framework, where instead of the Fukaya-Seidel category $\F_f(\Sigma)\simeq \F_{\{z_1,z_2\}}(\Sigma)$ one should work with the one basepoint partially wrapped Fukaya category $\F_z(\Sigma)$, and construct there an appropriate version of a twisted open-closed map.

\newcommand*{\arxivPreprint}[1]{ArXiv preprint \href{http://arxiv.org/abs/#1}{#1}}
\newcommand*{\arxiv}[1]{(ArXiv:\ \href{http://arxiv.org/abs/#1}{#1})}
\bibliographystyle{alpha}
\bibliography{references}
\end{document}